\documentclass[10pt,english]{article}
\usepackage[T1]{fontenc}
\usepackage[latin1]{inputenc}
\usepackage{amssymb}

\makeatletter

\newcommand{\noun}[1]{\textsc{#1}}



\usepackage{babel}
\makeatother
\begin{document}

\title{Complex Spaces and Nonstandard Schemes}

\author{Adel \noun{Khalfallah}, Siegmund \noun{Kosarew}}

\maketitle
\textbf{Summary}. Nonstandard mathematics furnishes a remarkable connexion
between analytic and algebraic geometry. We describe this interplay
for the most basic notions like complex spaces/algebraic schemes,
generic points, differential forms etc. We obtain - by this point
of view - in particular new results on the prime spectrum of a Stein
algebra.

\section*{Introduction}

\addcontentsline{toc}{section}{Introduction} The methods of nonstandard
mathematics are in general ignored in analytic and algebraic geometry.
Only some very specific applications of model theory are used to be
known as for instance the Lefschetz principle, the theorem of Tarski-Seidenberg
or some simple proofs of Hilbert's Nullstellensatz. Our Leitmotif
is quite different and can be summarized by the following statements

~

$\bullet$ holomorphic functions (or convergent power series) should
be (the standard part of) polynomials of infinite (or hyperfinite)
degree

~

$\bullet$ complex spaces should be seen as hyperalgebraic schemes 

~

$\bullet$ generic points of irreducible complex spaces/schemes should
be certain nonstandard points

~

$\bullet$ describe differential forms as functions taking infinitesimal
values (Leibniz' vision)

~

$\bullet$ replace always {}``countable'' by {}``hyperfinite''.

~\\
This program is achieved in our paper. Clearly, other more specific
problems and conjectures in analytic/algebraic geometry can be reformulated
by our {}``nonstandard'' point of view, especially those which involve
the word {}``infinity''.

~

We describe briefly the essential content of this paper. One of our
fundamental constructions is that of a category of certain ringed
spaces, called bounded schemes, which contains the category of algebraic
$\mathbb{C}$- schemes and which admits an essential \emph{surjective}
functor, called the \emph{standard} \emph{part} \emph{functor}, to
the category of complex spaces. We generalize in this way the usual
passage from a nonstandard (bounded) real number to its standard part
(a non constructive and non trivial procedure). The advantage of this
new more algebraic category is that it allows to apply many constructions
of standard algebraic geometry which are not evident in the analytic
context. We obtain analytic results just by taking the standard part
functor. The essential surjectivity of the last one means, roughly
speaking, that we replace holomorphic functions by nonstandard polynomials
(also called \emph{internal} polynomials) of some hyperfinite degree
(which we can fix a priorily). The polynomials obtained in this way
are called \emph{bounded} since they map bounded (nonstandard points)
into themselves. Analogously to the classical algebraic case, we can
extend sometimes holomorphic maps to compactifications in our nonstandard
algebraic setting. 

By a new version of what should be the meaning of a point in a complex
space, i.e. a nonstandard one, we are able to prove that any prime
ideal (not necessarily closed) of a Stein algebra which satisfies
a Nullstellensatz, involving nonstandard points, is the zero set of
some (non unique) nonstandard point. Especially, we can describe geometrically
\emph{all} maximal ideals of a Stein algebra and determine their residue
fields which are $\mathbb{C}$ or a nonstandard complex number field
$^{*}\mathbb{C}$

~

$\bullet$ \emph{every maximal ideal in a Stein algebra is the vanishing
ideal of a (nonstandard) point.}

~

Another fact which we obtain is a geometric visualisation of generic
points in irreducible complex spaces. They are found to be certain
nonstandard (but bounded) points. The evaluation map in such a point
defines an inclusion 

~

$\mathcal{M}(X)\hookrightarrow{}^{*}\mathbb{C}$

~ \\
of the field $\mathcal{M}(X)$ of meromorphic functions on $X$ to
the field $^{*}\mathbb{C}$ of nonstandard complex numbers.

Finally, in the tradition of G.W.Leibniz, we propose a {}``real''
infinitesimal interpretation of the symbol $dx$ in our context and,
more generally, the notion of a differential form. Such a $dx$ is
an \emph{infinitesimal} \emph{variable} i.e. (an equivalence class
of) a function which maps near standard points to infinitesimal ones.
This notion is rigorously defined in section 6. We show that our differential
forms on $X$ identify naturally with the holomorphic ones on the
associated complex space $X^{an}$. 

In many of our construction (in fact, in the most interesting ones),
we work with so-called external sets (like for instance, the set of
infinitesimal internal polynomials). This means in particular that
the transfer principle of nonstandard mathematics - its most powerful
tool - does not work directly here. But it is still important for
intermediate steps in our arguments. The reader show notice that,
in contrast to the language used in some books on nonstandard analysis,
we use in this paper always the expression {}``bounded'' point instead
of a {}``limited'' one and similarily also for other objects, as
internal polynomials for example.\\

The classical text on nonstandard mathematics is Robinson's book \cite{Robinson1}
and also that of Stroyan and Luxemburg \cite{Stroyan Luxemburg}.
For a more recent introduction into hyperreals see \cite{Goldblatt}
and for an account to hyper categories \cite{Brunjes1}. Other interesting
rings of nonstandard numbers, motivated by asymptotic expansions,
are considered in \cite{Lightstone}.

\tableofcontents{}

\section{Algebras of internal polynomials}

Let $I$ be an infinite set and $\mathcal{U}$ be a fixed non-principal
ultrafilter on $I$.

\subsection{Ultraproduct of rings}

Let $(A_{i})_{i\in I}$ be a family of sets indexed by $I$. We write
$^{*}A_{I}=\Pi{}_{\mathcal{U}}A_{i}$ for the ultraproduct of the
$A_{i}$'s with respect to the ultrafilter $\mathcal{U}$. An element
$x$ of $^{*}A_{I}$ is an equivalence class of $(x_{i})$. We write
$x=[x_{i}]$. If $A_{i}=A$ for each $i\in I$, then $^{*}A=\Pi{}_{\mathcal{U}}A_{i}$
is the ultrapower of $A$ with respect to $\mathcal{U}.$ Let $f_{i}:A_{i}\longrightarrow B_{i}$
be a family of maps between $A_{i}$ and $B_{i}$. It induces a map
$^{*}f:\,^{*}A\longrightarrow\,^{*}B$, defined by $^{*}f([a_{i}])=[f_{i}(a_{i})]$.
Such map is called an \emph{internal map} and we write $^{*}f=[f_{i}]$.\\

Now, we consider algebraic structures on the $A_{i}$'s. If $(A_{i})_{i}$
is a family of rings, then $^{*}A_{I}$ is again a ring. In fact $^{*}A_{I}=\Pi{}_{i\in I}A_{i}/\mathfrak{I}$,
where $\mathfrak{I}$ is the ideal in $\Pi{}_{i\in I}A_{i}$ defined
by $\mathfrak{I}=\{ x\in\Pi{}_{i\in I}A_{i},\,\mathrm{V(}x)\in\mathcal{U}\}$
where $\mathrm{V}(x):=\{ i\in I\mid x_{i}=0\}$. Rings of this form
are called \emph{{*}-rings} or \emph{internal rings}. Let $^{*}A$
and $^{*}B$ be two internal rings. An \emph{internal homomorphism}
\emph{of {*}-rings} between $^{*}A$ and $^{*}B$ is given by $^{*}f=[f_{i}]$
where $f_{i}:A_{i}\longrightarrow B_{i}$ are morphisms of rings.
{*}-rings and internal morphisms of {*}-rings constitute a category
of {*}-rings. The operator {}``{*}'' gives a functor between the
category of rings, parametrized by the indexset $I$, and the category
of {*}-rings. If $A$ is an internal ring, $(\mathcal{\textrm{Alg}}{}_{A}^{int})$
denotes the category of internal $A$-algebras.

\subsection{Internal algebras}

Let $^{*}S_{I}=\Pi{}_{\mathcal{U}}S_{i}$ be an internal ring and
$n=[n_{i}]\in{}^{*}\mathbb{N}$. Let $S_{i}[X_{1},\ldots,X_{n_{i}}]$
represent the ring of polynomials in indeterminates $X_{1},\ldots,X_{n_{i}}$
over the ring $S_{i}$ and $S[X_{1},\ldots,X_{n}]_{int}$ denote the
ultraproduct of the polynomial rings $S_{i}[X_{1},\ldots,X_{n_{i}}]$
with respect to the ultrafilter $\mathcal{U}$, that is,

\[
S[X_{1},\ldots,X_{n}]_{int}=\prod_{\mathcal{U}}S_{i}[X_{1},\ldots,X_{n_{i}}].\]

Given $B\in(\mathcal{\textrm{Alg}}{}_{S}^{int}$), we say that $B$
is an \emph{{*}-algebra of hyperfinite type over $S$} (resp. \emph{{*}-algebra
of finite type over $S$})\emph{,} if there exist an integer $n\in\,^{*}\mathbb{N}$
(resp. a finite $n\in\,\mathbb{N}$) and a surjective internal morphism
$u:\, S[X_{1},\ldots,X_{n}]_{int}\longrightarrow B$. Let $\mathfrak{I}=\mathrm{Ker}\, u$,
then $\mathfrak{I}$ is an internal ideal of $S[X_{1},\ldots,X_{n}]_{int}$
and $B$ is internally isomorphic to $S[X_{1},\ldots,X_{n}]_{int}/\mathfrak{I}$.
Furthermore, the morphism $u$ is completely determined by the internal
sequence $(u(X_{1}),\ldots,u(X_{n}))$. Conversely, any hyperfinite
family $(t_{1},\ldots,t_{n})\in B^{n}$ determines a unique internal
morphism $u:\, S[X_{1},\ldots,X_{n}]_{int}\longrightarrow B$ such
that $u(X_{i})=t_{i}$ for each $i=1,\ldots,n$.\\
 If the ring $^{*}S_{I}$ is {*}-noetherian (i.e. $\left\{ i\in I,\, S_{i}\,\mbox{is noetherian }\right\} \in\mathcal{U}$),
then the ideal $\mathfrak{I}$ is hyperfinitely generated by a family
of nonstandard polynomials.

Let $B$ and $C$ be two {*}-algebras of hyperfinite type over $^{*}S$.
Then there exist $n,m\in\,^{*}\mathbb{N}$ and surjective morphisms
$u:\, S[X_{1},\ldots,X_{n}]_{int}\longrightarrow B$, $v:\, S[Y_{1},\ldots,Y_{m}]_{int}\longrightarrow C$,
such that we have internal isomorphisms 

$B\simeq S[X_{1},\ldots,X_{n}]_{int}/\mathfrak{I}$ , $C\simeq S[Y_{1},\ldots,Y_{m}]_{int}/\mathfrak{J}$.

Let $w:\, B\longrightarrow C$ be a morphism of {*}-algebras. Then
the morphism $w$ lifts to a morphism of {*}-algebras $\tilde{w}:\, S[X_{1},\ldots,X_{n}]_{int}\longrightarrow S[Y_{1},\ldots,Y_{m}]_{int}$,
sending the ideal $\mathfrak{I}$ into the ideal $\mathfrak{J}$.\\

In the sequel, let $n$ be a fixed \emph{finite} integer.

\subsection{The ring of internal polynomials}

Let $^{*}S_{I}=[S_{i}]_{I}$ be an internal ring. Elements of $S[X_{1},\ldots,X_{n}]_{int}$
are called \emph{internal polynomials} or nonstandard polynomials.
Suppose now we are given $d\in{}^{*}\mathbb{N}$. An internal polynomial
$P\in S[X_{1},\ldots,X_{n}]_{int}$ of degree at most $d$, will be
written in the uniqueform

\begin{eqnarray*}
P & = & \sum_{|\nu|\leq d}a_{\nu}\, X^{\nu},\,\,\textrm{where }a_{\nu}\in\,^{*}S\textrm{ and }d\in{}^{*}\mathbb{N}.\end{eqnarray*}

Here the sum ranges over all multi-indices $\nu=(\nu_{1},\ldots,\nu_{n})\in\,^{*}\mathbb{N}^{n}$
with $|\nu|=\nu_{1}+\ldots+\nu_{n}\leq d$ and as usual $X^{\nu}$
is stands for $X_{1}^{\nu_{1}}\ldots X_{n}^{\nu_{n}}$. By transfer,
the nonstandard notion of degree is defined. We have a canonical \emph{injective}
morphism of $({}^{*}S_{I})$-algebra $({}^{*}S_{I})[X_{1},\ldots,X_{n}]\longrightarrow S[X_{1},\ldots,X_{n}]_{int}$,
where $({}^{*}S_{I})[X_{1},\ldots,X_{n}]$ denotes the ring of polynomials
in indeterminates $X_{1},\ldots,X_{n}$ over the internal ring $^{*}S_{I}$.
Furthermore, let $P\in S[X_{1},\ldots,X_{n}]_{int}$ be an internal
polynomial of degree $d$, then $P\in{}(^{*}S_{I})[X_{1},\ldots,X_{n}]$
if and only if, $d$ is \emph{finite}. \\
 \\

We denote by $({}^{*}S_{I})[\![X_{1},\ldots,X_{n}]\!]$ the ring of
of power series over the internal ring $^{*}S_{I}$ in indeterminates
$X_{1},\ldots,X_{n}$. \\

\textbf{Proposition 1.3.1} \emph{\label{Proposition1.3.1}There is
a canonical surjective ring homomorphism}\[
\theta:\, S[X_{1},\ldots,X_{n}]_{int}\longrightarrow(^{*}S_{I})[\![X,\ldots,X_{n}]\!]\]
 \emph{given by {}``forgetting monomials of infinite degree''.}

\emph{Proof}\textbf{.} Let $P\in S[X_{1},\ldots,X_{n}]_{int}$. Then
$P=\sum_{|\nu|\leq d}a_{\nu}\, X^{\nu},\,\,\textrm{where }a_{\nu}\in\,^{*}S_{I}\textrm{ and }d\in\,^{*}\mathbb{N}$.
Let $\theta$ be the restriction map\[
\theta(P)=\sum_{\nu\in\mathbb{N}^{n}}a_{\nu}\, X^{\nu}\in\,(^{*}S_{I})[\![X,\ldots,X_{n}]\!]\]
 i.e. $\theta(P)$ is the standard power series with coefficients
in $^{*}S_{I}$. It is clear that the map $\theta$ is a ring homomorphism
and $\theta$ maps \emph{onto} $(^{*}S_{I})[\![X,\ldots,X_{n}]\!]$.
In fact, given $Q\in(^{*}S_{I})[\![X,\ldots,X_{n}]\!]$, $Q=\Sigma{}_{\nu\in\mathbb{N}^{n}}a_{\nu}\, X^{\nu}$,
by comprehensiveness (see for example \cite{Lightstone} Ch.2 §6),
the sequence $(a_{\nu})_{\nu\in\mathbb{N}^{n}}$ extends to an internal
sequence, also denoted by $(a_{\nu})_{\nu\in\,^{*}\mathbb{N}^{n}}$.
Fix any $d\in\,^{*}\mathbb{N}_{\infty}$ and put $P=\sum_{|\nu|\leq d}a_{\nu}\, X^{\nu}$.
Then clearly we have $\theta(P)=Q$. -

\subsection{The ring of bounded internal polynomials}

Let $\,^{*}\mathbb{C}$ be an enlargement of the field of complex
numbers. Elements of the ring $\mathbb{C}[X_{1},\ldots,X_{n}]_{int}$
can be considered as internal functions between $^{*}\mathbb{C}^{n}$
and $^{*}\mathbb{C}$. In the sequel, we will not distinguish between
internal polynomials and their associated internal functions. As above,
each $P\in\mathbb{C}[X_{1},\ldots,X_{n}]_{int}$ can be written in
the form $P=\Sigma{}_{|\nu|\leq d}a_{\nu}\, X^{\nu},$ where $a_{\nu}\in\,^{*}\mathbb{C}$
$\textrm{ and }d\in\,^{*}\mathbb{N}$. We associate to $P$ the internal
polynomial $\mid P\mid$ defined by

\[
\mid P\mid=\sum_{|\nu|\leq d}|a_{\nu}|\, X^{\nu}.\]

Let us fix some notations: $^{b}\mathbb{C}$ stands for the set of
bounded points of $^{*}\mathbb{C}$ and $^{i}\mathbb{C}$ for the
set of infinitesimal points of $^{*}\mathbb{C}$ . We denote by $^{*}\mathbb{N}_{\infty}$
the set of infinite integers. If $n,m{\in\,}^{*}\mathbb{N}$ with
$n\leq m$, ~$[\![n\ldots m]\!]:=\{ k\in{}^{*}\mathbb{N},n\leq k\leq m\}$.
It is evident that $^{b}(\mathbb{C}^{n})=({}^{b}\mathbb{C})^{n}$
and $^{i}(\mathbb{C}^{n})=({}^{i}\mathbb{C})^{n}$. \\

\textbf{Definition 1.4.1} Let \emph{$P\in\mathbb{C}[X_{1},\ldots,X_{n}]_{int}$}
be \emph{}an \emph{}internal \emph{}polynomial\emph{.} We \emph{}call
\emph{$P$ }

i) a \emph{bounded polynomial} if $P({}^{b}\mathbb{C}^{n})\subset{}{}^{b}\mathbb{C}$
\emph{,} that \emph{}is\emph{, $P$} sends \emph{}bounded \emph{}points
\emph{}of \emph{$^{*}\mathbb{C^{\textrm{n}}}$} to \emph{}bounded
\emph{}points \emph{}of \emph{$^{*}\mathbb{C}$,}

ii) an \emph{absolutely bounded polynomial} if $\mid P\mid({}^{b}\mathbb{C}^{n})\subset{}{}^{b}\mathbb{C}$
\emph{,} that \emph{}is\emph{, $\mid P\mid$} sends \emph{}bounded
\emph{}points \emph{}of \emph{$^{*}\mathbb{C^{\textrm{n}}}$} to \emph{}bounded
\emph{}points \emph{}of \emph{$^{*}\mathbb{C}$,}

iii) an \emph{infinitesimal polynomial} (resp\emph{. absolutely infinitesimal
polynomial}) if $P({}^{b}\mathbb{C}^{n})\subset{}{}^{i}\mathbb{C}$
(resp\emph{.} $\mid P\mid({}^{b}\mathbb{C}{}^{n})\subset{}{}^{i}\mathbb{C})$
\emph{.}

\emph{~}

Let $\,^{b}\mathbb{C}[X_{1},\ldots,X_{n}]$ denote the set of bounded
internal polynomials. It is a subring of $\mathbb{C}[X_{1},\ldots,X_{n}]_{int}$.
The subset $\,^{i}\mathbb{C}[X_{1},\ldots,X_{n}]$ of infinitesimal
internal polynomials is an ideal of $\,^{b}\mathbb{C}[X_{1},\ldots X_{n}]$.\\

\noun{Notation}. Let$\,^{bs}\mathbb{C}[X_{1},\ldots,X_{n}]:=\{ P=\Sigma a_{\nu}\, X^{\nu}\in\,^{b}\mathbb{C}[X_{1},\ldots,X_{n}]\,\mid\, a_{\nu}\in\mathbb{C}\mbox{ for every \emph{standard} }\nu\}$.
Trivially it is a subring of $\,^{b}\mathbb{C}[X_{1},\ldots,X_{n}]$.\\

If the degree of an internal polynomial is \emph{finite}, we can immediately
characterize bounded or infinitesimal polynomials by their coefficients
and the notions of boundedness and absolutely boundedness coincide.\\

\noun{Remark}. Let $P=\Sigma{}_{|\nu|\leq d}a_{\nu}\, X^{\nu}$
be an internal polynomial of \emph{bounded} degree $d$, where $a_{\nu}\in\,^{*}\mathbb{C}$.
Then

i) $P$ is a bounded polynomial if and only if $a_{\nu}\in\,{}^{b}\mathbb{C}$
for each $\nu\in\mathbb{N}^{n}$, such that $|\nu|\leq d$,

ii) $P$ is an infinitesimal polynomial if and only if $a_{\nu}\in\,{}^{i}\mathbb{C}$
for each $\nu\in\mathbb{N}^{n}$, such that $|\nu|\leq d$.

We have the inclusions

\[
({}^{b}\mathbb{C})[X_{1},\ldots,X_{n}]\subsetneq\,^{b}\mathbb{C}[X_{1},\ldots,X_{n}]\subsetneq\mathbb{C}[X_{1},\ldots,X_{n}]_{int}.\]
\\
Now, we consider the general case \\

\textbf{Proposition 1.4.2\label{Proposition-1.4.2}} \emph{Let $P\in\mathbb{C}[X_{1},\ldots,X_{n}]_{int}$
be an internal polynomial, i.e.}

\[
P=\sum_{|\nu|\leq d}a_{\nu}\, X^{\nu},\,\textrm{ where }a_{\nu}\in{}^{*}\mathbb{C}\textrm{ and }d\in{}^{*}\mathbb{N}_{\infty}.\]

\emph{Then $P$ is absolutely bounded if and only if the following
two conditions are satisfied}

\emph{i) $a_{\nu}\in\,^{b}\mathbb{C}$ for bounded $|\nu|$ (i.e.
$\nu\in\mathbb{N}^{n})$,}

\emph{ii) $|a_{\nu}|^{\frac{1}{|\nu|}}\in\,^{i}\mathbb{C}$ for infinite
$|\nu|$ such that $|\nu|\leq d$.}

\emph{~}

\emph{Proof}\textbf{.} Let $P=\Sigma{}_{|\nu|\leq d}a_{\nu}\, X^{\nu}$
be an absolutely bounded internal polynomial, then for each $\nu\in\mathbb{N}^{n}$,
we have $|a_{\nu}|\leq\mid P\mid(1,\ldots1)$ so we conclude that
\emph{$a_{\nu}\in\,^{b}\mathbb{C}$} for finite $|\nu|$. Suppose
that there exists a standard positive real $m$, such that\[
|a_{\nu}|^{\frac{1}{|\nu|}}\geq m\,\mathcal{\textrm{ for some infinite }}|\nu|\textrm{ such that }|\nu|\leq d\]

Let $q$ be a standard real such that $q>1$ and $\xi$ be a standard
positive real such that $\xi>\frac{q}{m}$. Then,

\[
\mid P\mid(\xi,\ldots,\xi)\geq|a_{\nu}|\,\xi^{|\nu|}\geq(m\,\xi)^{|\nu|}\geq q^{|\nu|}\]
\\
hence, $\mid P\mid(\xi,\ldots,\xi)$ will be infinite which is a contradiction.

Now let $\varepsilon>0$ be a standard positive real. We set

\[
A_{\varepsilon}=\{|\nu|\in\,^{*}\mathbb{N}\:,\,|a_{\nu}|^{\frac{1}{|\nu|}}\leq\varepsilon\}.\]
 Then $A_{\varepsilon}$ is an internal subset of $^{*}\mathbb{N}$
which contains $\{ n\in{}^{*}\mathbb{N}_{\infty},\, n\leq d\}$. By
the permanence principle, there exists a finite integer $n_{0}\in\mathbb{N}$,
such that $[\![n_{0}..d]\!]\subset A_{\varepsilon}$. This means that,

\[
|a_{\nu}|^{\frac{1}{|\nu|}}\leq\varepsilon\textrm{ for every }\nu\in\,^{*}\mathbb{N}^{n}\,:\,|\nu|\in[\![n_{0}..d]\!].\]
 Let $x=(x_{1},\ldots,x_{n})$ be a bounded point in $^{*}\mathbb{C}^{n}$,
such that $|x_{i}|\leq\frac{1}{2\varepsilon}$ for each $i\in[\![1\ldots n]\!]$.
Then, we have

\[
|a_{\nu}\, x^{\nu}|\leq\varepsilon^{|\nu|}|x^{\nu}|\leq\left(\frac{1}{2}\right)^{|\nu|}\textrm{ for each }\nu\in\,^{*}\mathbb{N}^{n}\,:\,|\nu|\in[\![n_{0}..d]\!].\]
 The two internal polynomials ${\displaystyle \Sigma{}_{n_{0}\leq|\nu|\leq d}|a_{\nu}|\, X^{\nu}}$
and ${\displaystyle \Sigma{}_{|\nu|<n_{0}}|a_{\nu}|\, X^{\nu}}$ are
bounded for each $x=(x_{1},\ldots,x_{n})$, such that $|x_{i}|\leq\frac{1}{2\varepsilon}$.

Now, let $x=(x_{1},\ldots,x_{n})$ be a bounded point in $^{*}\mathbb{C}^{n}$.
Then there exists a standard positive real $\varepsilon$ such that
$|x_{i}|\leq\frac{1}{2\varepsilon}$ for each $i\in[\![1\ldots n]\!].$
Hence $\tilde{P}(x_{1},\ldots x_{n})$ is bounded, q.e.d.\\

\textbf{Proposition 1.4.3\label{Proposition-1.4.3}} \emph{Let $P\in\mathbb{C}[X_{1},\ldots,X_{n}]_{int}$
be an internal polynomial, so}

\[
P=\sum_{|\nu|\leq d}a_{\nu}\, X^{\nu}\,\textrm{ where }a_{\nu}\in{}^{*}\mathbb{C}\textrm{ and }d\in{}^{*}\mathbb{N}_{\infty}.\]

\emph{Then $P$ is absolutely infinitesimal if and only if $\,|a_{\nu}|^{\frac{1}{|\nu|}}\in\,^{i}\mathbb{C}$
for each $\nu\in\,^{*}\mathbb{N}^{n}\setminus\{0\}$, $|\nu|\leq d$,
and $a_{0}$ infinitesimal.} \\

The last condition is equivalent to $a_{\nu}\in\,^{i}\mathbb{C}$,
if $|\nu|$ is finite (i.e. $\nu\in\mathbb{N}^{n})$ and $|a_{\nu}|^{\frac{1}{|\nu|}}\in\,^{i}\mathbb{C}$,
if $|\nu|\leq d$ and $|\nu|$ is infinite.\\

\emph{Proof}\textbf{.} Let us define\[
A:=\{|a_{\nu}|^{\frac{1}{|\nu|}},\,0<|\nu|\leq d\}.\]
 Then $A$ is an hyperfinite set and so has a greatest element. Let
$M=\textrm{max }A$, so $M\in A$ and $M\in\,^{i}\mathbb{C}$.

Let $x=(x_{1},\ldots,x_{n})$ be a bounded point of $^{*}\mathbb{C}^{n}$.
Then there exists a standard positive real $R>0$ such that $|x_{i}|\leq R$
for each $i\in[\![1\ldots n]\!]$. We get\[
\mid\mid P\mid(x_{1},\ldots,x_{n})|\leq\mid a_{0}\mid+\sum_{1\leq|\nu|\leq d}(M\, R)^{|\nu|}\leq|a_{0}|+\frac{MR}{(1-MR)^{n}}\in\,^{i}\mathbb{C}.\]

In order to show the other implication, we note that the internal
polynomial $P$ is in particular absolutely bounded, so by proposition
1.4.2 , we have $|a_{\nu}|^{\frac{1}{|\nu|}}\in\,^{i}\mathbb{C}$
for infinite $|\nu|$. For each $\nu\in\mathbb{N}^{n}$, we have $|a_{\nu}|\leq\mid P\mid(1,\ldots1)$.
Then, we conclude that \emph{$a_{\nu}\in\,^{i}\mathbb{C}$} for finite
$|\nu|$, q.e.d.\\

In fact, we prove that the two notions of boundedness defined above
finally coincide :\\

\textbf{Proposition 1.4.4} \emph{Let $P\in\mathbb{C}[X_{1},\ldots,X_{n}]_{int}$
be an internal polynomial of degree $d\in\,^{*}\mathbb{N}$. We have }

\emph{i) $P$ is a bounded polynomial if and only if $P$ is an absolutely
bounded polynomial,}

\emph{ii) $P$ is an infinitesimal polynomial if and only if $P$
is an absolutely infinitesimal polynomial.}

\emph{Proof}\textbf{.} Let $x=(x_{1},\ldots,x_{n})$ be a point of
$^{*}\mathbb{C}^{n}$, we have $|P(x_{1},\ldots,x_{n})|\leq\mid P\mid(|x_{1}|,\ldots,|x_{n}|)$.
So, it is clear that every absolutely bounded (resp. infinitesimal)
polynomial is a bounded (resp. an infinitesimal) polynomial. Now we
verify the converse. Let $P=\Sigma{}_{|\nu|\leq d}a_{\nu}\, X^{\nu}\,$
be a bounded polynomial and $R$ be a standard positive real. By proposition
1.4.2, we have to prove that $a_{\nu}$ is bounded for finite $|\nu|$
and $|a_{\nu}|^{\frac{1}{|\nu|}}$ is infinitesimal for infinite $|\nu|$.
We denote by $T_{R}=\{(\xi_{1},\ldots\xi_{n})\in\,^{*}\mathbb{C}^{n},\,|\xi_{1}|=\ldots=|\xi_{n}|=R\}$.
Applying transfer to the Cauchy formula, we obtain\[
a_{\nu}=\frac{1}{(2\pi i)^{n}}\int_{T_{R}}\,\frac{P(\xi_{1},\ldots,\xi_{n})}{\xi_{1}^{\nu_{1}+1}\ldots\xi_{n}^{\nu_{n}+1}}\, d\xi_{1}\ldots d\xi_{n}\]
 Again, by transfer, the polynomial $P$ attains its maximum on the
Shilov boundary of the polydisc at some point $\xi_{R}\in T_{R}$.
But since $P$ is bounded, we have $|P(\xi_{R})|$ a bounded number.
Hence, there exists a standard positive real $M_{R}$ such that

\[
|a_{\nu}|\leq\frac{M_{R}}{R^{|\nu|}},\,\forall\textrm{ }\nu\in\,{}^{*}\mathbb{N}^{n},\,|\nu|\leq d\]
 So, if $|\nu|$ is bounded then $a_{\nu}$ is bounded and if $|\nu|$
is infinite, then $M_{R}^{\frac{1}{|\nu|}}\approx0$ and we have

\[
|a_{\nu}|^{\frac{1}{|\nu|}}\leq\frac{2}{R}\textrm{ for each }\nu\textrm{ such that }|\nu|\in\,{}^{*}\mathbb{N}_{\infty}\textrm{, }|\nu|\leq d.\]
 Since $R$ is an arbitrary standard positive real, we get $|a_{\nu}|^{\frac{1}{|\nu|}}\in\,^{i}\mathbb{C}$
for infinite $|\nu|$, q.e.d.\\

\textbf{Corollary 1.4.5} \emph{The rings} $\,^{b}\mathbb{C}[X_{1},\ldots,X_{n}]$
\emph{and} $\,^{i}\mathbb{C}[X_{1},\ldots,X_{n}]$ \emph{are invariant
by any partial derivative} $\partial^{\alpha}$ \emph{for $\alpha\in\mathbb{N}^{n}$.}

~

\noun{Remark} 1.4.6 There is a more general version of the preceding
results (including Cauchy's formula for the coefficients), assuming
only that the given internal polynomial takes bounded/infinitesimal
values on a polydisc with appreciable multiradius.

\section{Comparing holomorphic functions and internal polynomials }

We already have defined the ring of bounded internal polynomials over
$^{*}\mathbb{C}$. Let $P$ be a bounded internal polynomial. Then
$P$ can be seen as internal function from the set of bounded points
of $^{*}\mathbb{C}^{n}$ to the set of bounded points of $^{*}\mathbb{C}$.
First, we prove that $^{°}P$, the standard part of $P$, is an entire
holomorphic function on $\mathbb{C}^{n}$. This result is an easy
corollary of the theorem of Robinson-Callot which we extended to several
complex variables. \\

In this section, we define the so-called \emph{standard part functor}
from the category of bounded polynomial algebras to the category of
Stein algebras of finite embedding dimension. We prove that this functor
is essentially surjective. This fact can be regarded as a nonstandard
algebraization of Stein spaces.

\subsection{Stein Algebras}

We recall some known facts about Stein algebras. Let $(X,\mathcal{O}_{X})$
be a complex space. Then $X$ is called Stein space, if it is holomorphically
separable and holomorphically convex. The algebra $\mathrm{\Gamma}(X,\mathcal{O}_{X})$
of all holomorphic functions on $X$ is a Fréchet algebra. A topological
algebra $A$ (over $\mathbb{C}$) is called a \emph{Stein algebra,}
if there exists a Stein space $(X,\mathcal{O}_{X})$ such that $A$
is morphically isomorphic to $\mathrm{\Gamma}(X,\mathcal{O}_{X})$.
\\

Stein algebras form a category where morphisms are morphisms of topological
$\mathbb{C}-$ algebras. The functor of global sections defines a
contravariant functor from the category of Stein spaces to the category
of Stein algebras. In fact, this functor is an \emph{anti-equivalence}
between these categories, see Forster\cite{Forster} for the algebro-topological
theory of Stein algebras.\\

We recall some results due to H.Cartan. Let $A=\mathrm{\Gamma}(X,\mathcal{O}_{X})$
be a Stein algebra, and $\mathfrak{a}$ is a \emph{closed} ideal of
$A$. Then $\mathcal{O}_{X}.\mathfrak{a}$ is a coherent sheaf of
ideals and $\mathrm{\Gamma(}X,\mathcal{O}_{X}\mathfrak{.a})$ is isomorphic
to $\mathfrak{a}$. Conversely, let $\mathcal{M}$ be a coherent sheaf
of ideals of $\mathcal{O}_{X}$. Then $\mathrm{\Gamma(}X,\mathcal{M})$
is a closed ideal of $A$ and $\mathcal{O}_{X}\mathrm{.\Gamma(}X,\mathcal{M})$
identifies with $\mathcal{M}$. \\

Let $A$ be a Stein algebra of finite embedding dimension. Then $A=\mathrm{\Gamma}(X,\mathcal{O}_{X})$,
where $X$ is Stein space of finite embedding dimension. It follows
from the proper embedding theorem that $X$ can be embedded as a closed
complex subspace of $\mathbb{C}^{n}$ for some $n\in\mathbb{N}$.

\subsection{Internal holomorphic functions}

By the enlargement construction, one can define the set of internal
holomorphic functions on an internal open subset of $^{*}\mathbb{C}$.
Robinson studied internal holomorphic functions and proved some external
properties of these functions. \\

Let $\mathcal{U}$ be an ultrafilter defined on $\mathbb{N}$ and
$^{*}\mathbb{C}=\Pi{}_{\mathcal{U}}\mathbb{C}$ be an enlargement
of $\mathbb{C}$. Let $D=[D_{i}]$ be an internal subset of $^{*}\mathbb{C}^{n}$
and $f=[f_{i}]:\, D\longrightarrow^{*}\mathbb{C}$ be an internal
function. We say that $f$ is \emph{an internal holomorphic function}
on $D$ if \begin{eqnarray*}
\{ i & \in & \mathbb{N}\mid\, D_{i}\subset\mathbb{C}^{n}\textrm{open and }\, f_{i}\mbox{ holomorphic on }D_{i}\}\in\mathcal{U}.\end{eqnarray*}

It is well known, by Osgood's theorem, that a function $f$ is holomorphic
on $\Omega$, an open subset of $\mathbb{C}^{n}$, if and only if
$f$ is continuous on $\Omega$ and partially holomorphic on $\Omega$.\\

\textbf{Theorem 2.2.1} \emph{Let $B$ be a S-open subset of $^{*}\mathbb{C}^{n}$.
We fix a bounded point $a$ in $B$ and $f:\, B{\longrightarrow\,}^{*}\mathbb{C}$
an internal holomorphic function on $B$. We assume that $f$ is bounded
on $\mu(a)$, the halo of $a$. Then there exists $V$ a S-open neighbourhood
of $a$ such that}

\emph{i) f is S-continuous in $V$.}

\emph{ii) There exists ${}{}^{\circ}f:\,^{\circ}V\longrightarrow\mathbb{C}$
a holomorphic function in $^{\circ}V$, such that $^{\circ}f\approx f$
and $^{\circ}(\partial^{\alpha}f)=\partial^{\alpha}(^{\circ}f)$ for
each $\alpha\in\mathbb{N}^{n}$.}

\emph{iii) If ${}{}^{\circ}f$ is not constant in $^{\circ}V$, then
$f$ is S-open in $a$, that is, $f(\mu(a))=\mu(f(a))$.}

~\\
 This theorem is an easy generalization of the theorem of Robinson-Callot
known in single complex variable, see Robinson\cite{Robinson1}, Callot\cite{Callot},
Fruchard\cite{Fruchard} and Lutz-Goze(\cite{Lutz goze}: lesson 11,
p:123).\\

\noun{Notation}. For $\nu=(\nu_{1},\ldots,\nu_{n})\in\mathbb{N}^{n}$
and $z=(z_{1},\ldots,z_{n})\in\mathbb{C}^{n}$ define $|\nu|=\Sigma{}_{i=1}^{n}\nu_{i}$
and $z^{\nu}=z_{1}^{\nu_{1}}\ldots z_{n}^{\nu_{n}}$. Let $r=(r_{1},\ldots,r_{n})\in\,^{*}\mathbb{R}^{n}$
, all $r_{i}>0$, $z_{0}=(z_{1}^{(0)},\ldots,z_{n}^{(0)})\in\,^{*}\mathbb{C}^{n}$.
Then $\textrm{P}^{n}(z_{0},r)=\{ z{\in{}}^{*}\mathbb{C}^{n}\,\mid\,|z_{i}-z_{i}^{(0)}|<r_{i}\mathcal{\textrm{ for }}i=1,\ldots,n\}$
is called the \emph{{*}-open polydisk} with polyradius $r$ and center
$z_{0}$ and $\overline{\textrm{P}^{n}(a,r)}=\{ z{\in{}}^{*}\mathbb{C}^{n}\,\mid\,|z_{i}-z_{i}^{(0)}|\leq r_{i}\mathcal{\textrm{ for }}i=1,\ldots,n\}$
the \emph{{*}-closed polydisk}. The Shilov boundary of the {*}-polydisk
$\textrm{P}^{n}(z_{0},r)$ is the set $\textrm{T}^{n}(z_{0},r)=\{ z{\in{}}^{*}\mathbb{C}^{n}\,\mid\,|z_{i}-z_{i}^{(0)}|=r_{i}\mathcal{\textrm{ for }}i=1,\ldots,n\}$.
If $r>0$ and $r=(r,\ldots,r)$, we write $\textrm{P}_{r}^{n}(z_{0})$
(resp. $\textrm{T}_{r}^{n}(z_{0})$) instead of $\textrm{P}^{n}(z{}_{0},r)$
(resp. $\textrm{T}^{n}(z_{0},r)$).\\

\emph{Proof}\textbf{.} i) By permanence principle, there exists a
standard positive real $r$ such that $\overline{\textrm{P}^{n}(a,r)}\subset B$
and the internal holomorphic function $f$ is bounded on $\overline{\textrm{P}^{n}(a,r)}$.
Let $\rho$ be a standard positive such that $\rho<r$ and \[
V:=\{ z{\in{}}^{*}\mathbb{C}^{n}\,\mid\,^{\circ}|z_{i}-a_{i}|<\rho\mathcal{\textrm{ for }}i=1,\ldots,n\}\]
V is a S-open neighbourhood of $a$ and $^{\circ}V$ is \emph{}the
open polydisk with polyradius $\rho$ and center $°a$. Let $z,\xi\in V$
such that $z\approx\xi$. Applying transfer to Cauchy's formula, we
have\[
f(z)=\frac{1}{2\pi i}\int_{\textrm{T}_{r}^{n}(a)}\frac{f(w)}{w-z}\, dw\textrm{ and }f(\xi)=\frac{1}{2\pi i}\int_{\textrm{T}_{r}^{n}(a)}\frac{f(w)}{w-\xi}\, dw\]
 where $w-z=(w_{1}-z_{1})\ldots(w_{n}-z_{n})$ and $w-\xi=(w_{1}-\xi_{1})\ldots(w_{n}-\xi_{n})$.

Again, by transfer, the function $f$ attains its maximum on $\textrm{T}_{r}^{n}(a)$,
the Shilov boundary of the polydisk $\textrm{P}^{n}(a,r)$. Since
$f$ is bounded on $\overline{\textrm{P}^{n}(a,r)}$, then $M:=\textrm{sup}_{w\in\textrm{T}^{n}(a,r)}|f(w)|$
is bounded. So,

\[
|f(z)-f(\xi)|\leq\frac{M}{(2\pi)^{n}}\int_{\textrm{T}_{r}^{n}(a)}\mid\frac{1}{w-z}-\frac{1}{w-\xi}|\, dw.\]
 It is easy to see that $w-z\approx w-\xi$ and, since $w-\xi$ is
appreciable, $\frac{1}{w-z}\approx\frac{1}{w-\xi}$. Moreover, there
exists $\eta$ a positive infinitesimal such that $\mid\frac{1}{w-z}-\frac{1}{w-\xi}|\leq\eta$
for every $w\in\textrm{T}_{r}^{n}(a)$. It follows that

\[
|f(z)-f(\xi)|\leq Mr^{n}\eta.\]
 This shows that $f(z)\approx f(\xi)$ for $z\approx\xi$ which is
equivalent to S-continuity of $f$ in $z$, see Appendix C.\\

ii) Let $^{\circ}a={({}}^{\circ}a_{1},\ldots,^{\circ}a_{n})$ be the
standard part of $a$. The standard part $^{\circ}f$ of the function
$f$ exist on $^{\circ}V$, the standard polydisk with polyradius
$\rho$ and center $°a$. We deduce from the first assertion that
$^{\circ}f$ is continuous on $^{\circ}V$ (see Appendix C). By Osgood's
theorem, it suffices to prove that $f$ is partially holomorphic.
But this follows directly from Robinson \cite{Robinson1}( Theorem
6.2.3 p. 156) and for each $\alpha\in\mathbb{N}^{n}$, we have \emph{${}{}^{\circ}(\partial^{\alpha}f)=\partial^{\alpha}{({}}^{\circ}f)$}.\\

iii) The S-continuity of $f$ implies that $f(\mu(a))\subset\mu(f(a))$.
Let $w\approx f(a)$. By hypothesis, $^{\circ}f$, the standard part
of $f$, is not constant, hence, there exists $b\not=a$, such that
$^{\circ}f(b)\not\neq{}^{\circ}f(a)$. Let $D$ be the {*}-complex
line passing through $a$ and $b$. Let $f_{1}$ denotes the restriction
of $f$ on $D\cap\mu(a)$. $f_{1}$ is an internal holomorphic map
and not a constant function. Then, by Robinson's theorem \cite{Robinson1}
Theorem 6.2.8 p. 158, there exists $z\in D\cap\mu(a)$ such that $f(z)=w$,
q.e.d.\\

\noun{Remark 2.2.2} The condition that $^{\circ}f$ is not constant
is essential to deduce that $f$ is S-open. For a counterexample,
it suffices to consider the trivial example of constant infinitesimal
map. Moreover, this condition cannot be weakened by {}``$f$ is non
constant'': Let $\varepsilon$ be an infinitesimal and $f(z)=\varepsilon z$.
Then we have $f(\mu(0))\subsetneq\mu(0)$.

\subsection{Internal holomorphic maps into complex hyperbolic spaces}

In theorem (2.2.1), we note that the condition that $f$ is bounded
in every point of $\mu(a)$, the halo of $a$, is essential and cannot
be replaced by the condition that $f(a)$ is bounded as shows this
example. Let $\omega$ be an infinite real and $f(z)=\omega z$ be
the internal holomorphic function on $^{*}\mathbb{C}$. We have $f(0)=0$
but $f$ is not bounded on every point of $\mu(0)$. Hence, we could
think of imposing additional assumptions on the target space, but
garde a weak condition on $f$. In this direction A. Robinson proved
(see \cite{Robinson1} Theorem 6.3.2 p.160)

~

\textbf{Theorem 2.3.1} \cite{Robinson1} \emph{Let $B$ be a $S$-domain
in $^{*}\mathbb{C}$ and $f:B\rightarrow{}{}^{*}\mathbb{C}\setminus\{0,1\}$
be an internal holomorphic function. Suppose that $f$ takes a bounded
value at some point $z_{0}\in B$. Then $f$ is bounded and $S-$continuous
in $B.$}\\

We should note that this theorem is not true if $f$ takes values
in $^{*}\mathbb{C}\setminus\{0\}$ as shows the following example:
put $f(z):=\exp\omega z$, where $\omega$ is an infinite real; $f(0)$
is bounded, but $f$ is not $S$ continuous in $0$.

~

We want to generalize the above theorem of Robinson, motivated by
the fact that the space $\mathbb{C}\setminus\{0,1\}$ is complete
hyperbolic and, moreover, hyperbolically embedded in $\mathbb{P}_{\mathbb{C}}^{1}$.
\\

\noun{Notation. $\Delta_{r}$} stands for the disc with radius $r>0$
in the complex plane $\mathbb{C}$. If $r=1$, we simply denote $\Delta$
instead of $\Delta_{1}$. Let $X$ be a reduced complex space and
$\mathrm{Hol}(\Delta,X)$ the set of all holomorphic maps from $\Delta$
to $X$. We will denote by $d_{X}$ the Kobayashi pseudo-distance
of $X$.\\

\textbf{Theorem 2.3.2} \emph{Let $X$ be a hyperbolic complex space
and $f:\,^{*}\Delta\rightarrow\,^{*}X$ be an internal holomorphic
map. Then}

\emph{i) $f$ is S-continuous in $^{*}\Delta$.}

\emph{ii) Suppose that $f$ takes a bounded value at some point $z_{0}$
such $^{\circ}|z_{0}|<1$. Then $f$ is bounded in every point of
$\{ z\in\,^{*}\Delta\,|\,^{\circ}|z|<1\},$ the S-interior of $^{*}\Delta$.}

~

For the definition of bounded points in a metric space, see Appendix
C: metric spaces.

~

\emph{Proof.} For every \emph{$g\in\mathrm{Hol}(\Delta,X)$}, we have
\emph{$d_{X}(g(x),g(y))\leq d_{\Delta}(x,y)$} for any \emph{$x,y\in\Delta$}.
Hence by transfer, we obtain \emph{\[
^{*}d_{X}(f(x),f(y)){\leq{}}^{*}d_{\Delta}(x,y)\,\mathrm{\, for\, each\,}\, x,y{\in{}}^{*}\Delta.\]
 }

As a consequence, if $x,y\in^{*}\Delta$ such that $x\approx y$,
then $^{*}d_{X}(f(x),f(y))\approx0$ and $f$ is S-continuous on $^{*}\Delta$.
The assertion (ii) is an immediately deduced from the following lemma,
due to Robinson (see \cite{Robinson1}, Theorem 4.5.9 p:114)

~

\textbf{Lemma 2.3.3} \cite{Robinson1} \emph{Let $h$ be a $S$-continuous
map, defined on a $S$-connected set $D$. Then the points of $h(D)$
belong to the same galaxy.}

~

Next, we consider the case of complete hyperbolic complex spaces.

~

\textbf{Theorem 2.3.4} \emph{Let $f:\,^{*}\Delta\rightarrow\,^{*}X$
be an internal holomorphic map where $X$ is a complete hyperbolic
complex space. Then}

\emph{i) $f$ is S-continuous in $^{*}\Delta$.}

\emph{ii) Suppose that $f$ takes a near-standard value at some point
$z_{0}$ such $^{\circ}|z_{0}|<1$. Then $f$ is near-standard in
every point of $\{ z\in\,^{*}\Delta\,|\,^{\circ}|z|<1\},$ the S-interior
of $^{*}\Delta$ and there exists a holomorphic map $^{\circ}f:\Delta\rightarrow X$,
the standard part of $f$ verifying $^{\circ}f\approx f$.}

~

We recall the following definition

~

\textbf{Definition 2.3.5} Let \emph{$(X,d)$} be a metric \emph{}space\emph{.}
We \emph{}say \emph{}that \emph{$(X,d)$} is \emph{strongly complete}
or \emph{finitely compact} if \emph{}every \emph{}closed \emph{}ball
\emph{$B(x,r)=\{ y\in X\,|d(x,y)\leq r\}$} with \emph{$x\in X$}
and \emph{$r>0$} is \emph{}compact.\\

Comparing different notions of completeness, Kobayashi proved (see
\cite{Kobayashi}, proposition 1.1.9 p.4)\\

\textbf{Proposition 2.3.6} \cite{Kobayashi} \emph{Let $d$ be a distance
on a locally compact space $X$. Then}

\emph{i) If $(X,d)$} \emph{is strongly complete then $(X,d)$ is
complete.}

\emph{ii) If $d$ is inner distance then completeness implies strong
completeness}.\\

It is straightforward - using nonstandard characterization of compactness
- to check the following proposition\\

\textbf{Proposition 2.3.7} \emph{Let $(X,d)$ be a metric space. Then
$(X,d)$ is strongly complete if and only if $\mathrm{bd}{({}}^{*}X)=\mathrm{ns}{({}}^{*}X)$},
\emph{where} $\mathrm{bd}{({}}^{*}X)$ \emph{denotes} \emph{the} \emph{set}
\emph{of} \emph{bounded} \emph{points} \emph{in} $^{*}X$ \emph{(see}
\emph{Appendix C: metric} \emph{spaces}).\\

\emph{Proof (theorem 2.3.4)}. Using theorem 2.3.2, we only have to
prove assertion (ii). Since the Kobayashi distance $d_{X}$ is an
inner distance, we conclude that $X$ is in fact strongly complete
for $d_{X}$. As a consequence, bounded points of $^{*}X$ coincide
with near standard points. Now let $z{\in{}}^{*}\Delta$, such $^{\circ}|z|<1$.
By S-continuity, we get $f(\mu(z))\subset\mu(f(z))$. Let $W$ be
a S-neighborhood of $f(z)$ in $^{*}\mathbb{C}^{n}$. Applying the
permanence principle, there exists a standard positive real $r>0$
such that $f{({}}^{*}\Delta_{r})\subset W$. By theorem 2.2.1, we
conclude that the standard part $^{\circ}f$ exists, q.e.d.

~

Let $X$ be a complex manifold and $H$ be a metric on $TX$; for
simplicity, we will write $|v|$ instead of $H(v)$ for every $v\in TX$.

~

\textbf{Proposition 2.3.8} \emph{The complex manifold $X$ is hyperbolic
if and only if for every internal holomorphic map} $f:{}^{*}\Delta\rightarrow{}{}^{*}X$
\emph{with} $f(0)\in\mathrm{ns}{({}}^{*}X)$ \emph{verifies} $|f'(0)|$
\emph{is bounded.}

~

\emph{Proof}. Suppose that $X$ is hyperbolic and let $f:{}^{*}\Delta\rightarrow{}{}^{*}X$
be an internal holomorphic map. Then $f$ is S-continuous in $0$.
Hence for $W$ be a S-neighborhood of $f(0)$, there exists a standard
real $0<r<1$ such that $f{({}}^{*}\Delta_{r})\subset W$. By the
Cauchy formula, we deduce that $|f'(0)|$ is bounded. Conversely,
if $X$ is not hyperbolic then by the Royden infinitesimal criterion
of hyperbolicity, we conclude that there exist $x\in X$ and a sequence
of holomorphic maps $f_{n}:\Delta\rightarrow X$ such that the sequence
$(f_{n}(0))$ converges to $x$ and ${|f}_{n}'(0)|\rightarrow+\infty$.
Clearly the sequence $(f_{n})$ induces an internal holomorphic map
$F:{}^{*}\Delta\rightarrow{}{}^{*}X$ verifying $F(0)\in\mathrm{ns}{({}}^{*}X)$
and $|F'(0)|$ is not bounded. -

~

\noun{Remark}. Proposition 2.3.8 is a nonstandard translation of
the following characterization of hyperbolicity which asserts that
a manifold $X$ is hyperbolic if and only if $X$ satisfies the \emph{Landau
property} (see for instance \cite{Hahn Tim})\emph{,} that is, for
each $p\in X$ and each $W$ a relatively compact neighborhood of
$p$, there exists $R>0$ such that\[
\mathrm{sup}\{|f'(0)|:\, f\in\mathrm{Hol}(\Delta,X)\,\,\mathrm{with}\, f(0)\in W\}\leq R.\]

\textbf{Corollary} \textbf{2.3.9} \emph{Let $X$ be a compact manifold.
Then $X$ is hyperbolic if and only if, every internal holomorphic
map $f:\,^{*}\Delta\longrightarrow^{*}X$ satisfies $|f'(0)|$ is
bounded}.~\\
We close this section by giving a nonstandard characterization of
hyperbolicity in the compact case

~

\textbf{Theorem 2.3.10} \emph{i) Let $Y$ be a relatively compact
subspace of a complex space $Z$. Then $Y$ is hyperbolically embedded
in $Z$ if and only if every internal holomorphic map $f{:{}}^{*}\Delta\rightarrow{}{}^{*}Y$
has a holomorphic standard part $^{\circ}f:\Delta\rightarrow Z$.}

\emph{ii) Let $X$ be a compact complex space. Then $X$ is hyperbolic
if and only if every internal holomorphic map $f{:{}}^{*}\Delta{\rightarrow{}}^{*}X$
has a holomorphic standard part $^{\circ}f:\Delta\rightarrow X$.}

~\\
 This is a nonstandard interpretation of the following standard facts:

i) Let $Y$ be a relatively compact subspace of a complex space $Z$.
Then $Y$ is hyperbolically embedded in $Z$ if and only if \emph{$\mathrm{Hol}(\Delta,Y)$}
is relatively compact in \emph{$\mathrm{Hol}(\Delta,Z)$.}

ii) Let $X$ be a compact complex space. Then $X$ is hyperbolic if
and only if \emph{$\mathrm{Hol}(\Delta,X)$} is compact.

\subsection{Bounded internal holomorphic functions}

Similar as for the algebra of bounded internal polynomials, we define
$^{b}\mathcal{O}(\mathbb{C}^{n})$, the algebra (over $^{b}\mathbb{C}$)
of \emph{bounded} entire {*}- holomorphic functions on $^{*}\mathbb{C}^{n}$,
by

\[
^{b}\mathcal{O}(\mathbb{C}^{n})=\{ f\in\,{}^{*}\mathcal{O}(\mathbb{C}^{n})\,\mid\, f({}^{b}\mathbb{C}^{n})\subset\,{}^{b}\mathbb{C}\}.\]

\textbf{Proposition 2.4.1} \emph{Let $f\in{}^{b}\mathcal{O}(\mathbb{C}^{n})$
be a bounded {*}-holomorphic function on} \emph{$^{*}\mathbb{C}^{n}$.
Then} \emph{$^{\circ}f$, the standard of $f$ exists and is holomorphic
on $\mathbb{C}^{n}$. We have for the zero sets, ${}{}^{\circ}\mathrm{Z}(f)\subset\mathrm{Z}(\,^{\circ}f)$
and if $^{\circ}f$ is not constant, then \[
{}{}^{\circ}\mathrm{Z}(f)=\mathrm{Z}(\,^{\circ}f).\]
 Proof}\textbf{.} The first assertion is a direct consequence of theorem
2.2.1. For the second, the inclusion follows from the S-continuity
of $f$ and the S-openness of $f$ implies the equality. \\

Let now $X$ be a topological space. We denote by $\textrm{ns}(^{*}X)$
the set of \emph{near-standard} points of $^{*}X$. The knowledge
of $\textrm{ns}(^{*}X)$ is essential to define the standard part
map. \\

We endow $\mathbb{C}[X_{1},\ldots,X_{n}]$ and $\mathcal{O}(\mathbb{C}^{n})$
with the compact-open topology, that is the topology of uniform convergence
in each compact subset of $\mathbb{C}^{n}$.

For $g\in\mathcal{O}(\mathbb{C}^{n})$, let $\mu(g)$ denote the halo
of $g$ for the compact-open topology. \\

\textbf{Lemma 2.4.2} \emph{Let $g\in\mathcal{O}(\mathbb{C}^{n})$.
Then \[
\mu(g)=\{ f\in{}^{*}\mathcal{O}(\mathbb{C}^{n})\,\mid\, f(x)\approx\,^{*}g(x)\mathrm{\, for\, every\,}\, x{\in{}}^{b}\mathbb{C}^{n}\}.\]
 In addition, we have $\,^{b}\mathcal{O}(\mathbb{C}^{n})=\mathrm{ns}({}^{*}\mathcal{O}(\mathbb{C}^{n}))$
and $\,^{b}\mathbb{C}[X_{1},\ldots,X_{n}]=\mathrm{ns}({}^{*}\mathcal{O}(\mathbb{C}^{n}))\cap\,\mathbb{C}[X_{1},\ldots,X_{n}]_{int}$}.\\

\emph{Proof}\textbf{.} Let $\mathcal{C}(\mathbb{C}^{n},\mathbb{C})$
denote the set of continuous maps from $\mathbb{C}^{n}$ to $\mathbb{C}$
endowed with the compact-open topology. Let $g\in\mathcal{C}(\mathbb{C}^{n},\mathbb{C})$
and $\mu_{\mathcal{C}}(g)$ be the halo of $g$ in $^{*}\mathcal{C}(\mathbb{C}^{n},\mathbb{C})$.
Then we have that $\mu_{\mathcal{C}}(g)=\{ f\in\,{}^{*}\mathcal{C}(\mathbb{C}^{n},\mathbb{C})\,\mid\, f(x)\approx\,^{*}g(x)\textrm{ for every}\, x\in{}^{b}\mathbb{C}^{n}\}$,
see Appendix C: Standard part of a map. Since $\mathcal{O}(\mathbb{C}^{n})\subset\mathcal{C}(\mathbb{C}^{n},\mathbb{C})$
is equipped with the induced topology, then $\mu(g)=\mu_{\mathcal{C}}(g)\cap{}^{*}\mathcal{O}(\mathbb{C}^{n})$
for every $g\in\mathcal{O}(\mathbb{C}^{n})$.

It is straightforward to verify that $\textrm{ns}(^{*}\mathcal{O}(\mathbb{C}^{n}))\subset{}^{b}\mathcal{O}(\mathbb{C}^{n})$,
the converse inclusion is deduced from proposition 2.4.1. - \\

\noun{Remark} 2.4.3 If we put $^{b}\mathcal{C}(\mathbb{C}^{n},\mathbb{C}):=\{\, f\in{}{}^{*}\mathcal{C}(\mathbb{C}^{n},\mathbb{C})\,\mid\, f(\,^{b}\mathbb{C}^{n})\subset{}{}^{b}\mathbb{C}\}$
and $\mathcal{SC}(\mathbb{C}^{n},\mathbb{C}):=\{ f\in{}{}^{*}\mathcal{C}(\mathbb{C}^{n},\mathbb{C})\,\mbox{ S-continuous on }{}^{b}\mathbb{C}^{n}{}\}$.
Then we have, by Appendix C: Standard part of a map,\[
^{b}\mathcal{C}(\mathbb{C}^{n},\mathbb{C})\cap\mathcal{SC}(\mathbb{C}^{n},\mathbb{C})=\textrm{ns}({}^{*}\mathcal{C}(\mathbb{C}^{n},\mathbb{C})).\]

For $g\in\mathcal{O}(\mathbb{C}^{n})$ let $\mu_{s}(g)$ denote the
halo of $g$ in the topology of simple convergence. Then \[
\mu_{s}(g)=\{ f\in{}^{*}\mathcal{O}(\mathbb{C}^{n})\,\mid\, f(x)\approx\, g(x)\textrm{ for every}\, x\in\mathbb{C}^{n}\}.\]
~

\textbf{Proposition 2.4.4} \emph{The standard part map defines a ring
homomorphism}

\[
\textrm{st}:\,^{b}\mathcal{O}(\mathbb{C}^{n})\longrightarrow\mathcal{O}(\mathbb{C}^{n}).\]
 \emph{Its restriction to the subrings $^{b}\mathbb{C}[X_{1},\ldots,X_{n}]$
and $\,^{bs}\mathbb{C}[X_{1},\ldots,X_{n}]$ are surjective and we
have the following commutative diagram}

\[
\begin{array}{ccc}
 & st\\
^{bs}\mathbb{C}[X_{1},...,X_{n}] & \longrightarrow & \mathcal{O}(\mathbb{C}^{n})\\
\downarrow &  & \downarrow\\
\mathbb{C}[X_{1},...,X_{n}]_{int} & \longrightarrow & (^{*}\mathbb{C})[[X_{1},...,X_{n}]]\end{array}\]
\emph{where} \emph{the} \emph{vertical} \emph{arrows} \emph{are} \emph{the
natural} \emph{inclusions} \emph{and} \emph{the} \emph{lower} \emph{horizontal}
\emph{one} \emph{is} \emph{the} \emph{ring} \emph{homomorphism} \emph{$\theta$
defined} \emph{in} \emph{section} \ref{Proposition1.3.1}.\\

\emph{Proof}\textbf{.} The standard part map $"\mbox{st}"$ is well
defined since we know that $\,^{b}\mathcal{O}(\mathbb{C}^{n})=\textrm{ns}(^{*}\mathcal{O}(\mathbb{C}^{n}))$.
It defines a ring homomorphism, because the standard part map is compatible
with sums and products of complex numbers.

Now, we shall prove that the restriction of $"\mbox{st}"$ to \emph{$\mathbb{^{\mathcal{\textrm{bs}}}C}[X_{1},\ldots,X_{n}]$}
is surjective. Let $f\in\mathcal{O}(\mathbb{C}^{n})$, so\[
f=\sum_{\nu\in\mathbb{N}^{n}}a_{\nu}\, X^{\nu}.\]
 The sequence $(a_{\nu})_{\nu\in\mathbb{N}^{n}}$ extends to the internal
sequence $(a_{\nu})_{\nu\in^{*}\mathbb{N}^{n}}$. Let $N\in\,{}^{*}\mathbb{N}_{\infty}$.
Truncating $f$ at the order $N$, we get\[
f_{N}=\sum_{|\nu|\leq N}a_{\nu}\, X^{\nu}\]
 and $f_{N}$ is an internal polynomial. Since the partial sums of
$f$ converge on each compact of $\mathbb{C}^{n}$ to $f$, we have

\[
f_{N}(x)\approx\,{}^{*}f(x)\mbox{ for each }x{\in{}}^{b}\mathbb{C}^{n}\]
 that is, $f_{N}\in\mu(f)\subset{}{}^{{b}}\mathbb{C}[X_{1},\ldots X_{n}]$
. This shows $f_{N}\in{}{}^{{bs}}\mathbb{C}[X_{1},\ldots,X_{n}]$
and $\mbox{st}(f_{N})=f$.\\

The restriction of the map {}``$\mbox{st}$'' to $^{b}\mathbb{C}[X_{1},\ldots,X_{n}]$
can be described explicitly. Let $f\in{}{}^{{b}}\mathbb{C}[X_{1},\ldots,X_{n}]$,
of degree $d\in\,{}^{*}\mathbb{N}$ , and $f=\Sigma{}_{|\nu|\leq d}a_{\nu}\, X^{\nu}$.
Then \[
\mbox{st}(\sum_{|\nu|\leq d}a_{\nu}\, X^{\nu})=\sum_{\nu\in\mathbb{N}^{n}}\,^{°}(a_{\nu})\, X^{\nu}.\]
 \\
 First, it is evident that the power series $g(x):=\Sigma{}_{\nu\in\mathbb{N}^{n}}\,^{\circ}(a_{\nu})\, x^{\nu}$
defines an entire holomorphic function on $\mathbb{C}^{n}$. In fact,
let $\varepsilon>0$ be a standard positive real, there exists a finite
integer $n_{0}\in\mathbb{N}$, such that $|a_{\nu}|^{\frac{1}{|\nu|}}\leq\varepsilon\textrm{ for each }\nu\in\,^{*}\mathbb{N}^{n}\textrm{ such that }n_{0}\leq|\nu|\leq d.$
Hence, for each standard $\nu\in\,\mathbb{N}^{n}$ such that $|\nu|\geq n_{0}$,
we have ${}{}^{\circ}|a_{\nu}|^{\frac{1}{|\nu|}}\leq\varepsilon$
which implies that $\limsup_{|\nu\ |\rightarrow+\infty}{({}}^{\circ}\mid a_{\nu}|^{\frac{1}{|\nu|}})=0$.

Trivially the standard sequence of complex numbers ${({}}^{\circ}a_{\nu})_{\nu\in\mathbb{N}^{n}}$
extends to an internal sequence of complex numbers $(b_{n}){}_{n\in^{*}\mathbb{N}^{n}}$
such that for every \emph{standard} $\nu$, we have $b_{\nu}={}^{\circ}a_{\nu}$.
By transfer, we get $|b_{\nu}|^{\frac{1}{|\nu|}}\leq\varepsilon$
for every $|\nu|\geq n_{0}$, $\nu\in\,^{*}\mathbb{N}^{n}$. The uniform
convergence of $g$ on each compact gives us \[
^{*}g(x)\approx\sum_{|\nu|\leq d}b_{\nu}\, x^{\nu}\mbox{ for each }x\in{}{}^{b}\mathbb{C}^{n}.\]

We claim that $\Sigma{}_{|\nu|\leq d}a_{\nu}\, x^{\nu}\approx\Sigma{}_{|\nu|\leq d}b_{\nu}\, x^{\nu}\mbox{ for each }x{\in{}}^{b}\mathbb{C}^{n}$.
This is equivalent to $\Sigma{}_{|\nu|\leq d}(a_{\nu}-b_{\nu})\, x^{\nu}\in{}{}^{i}\mathbb{C}[X_{1},\ldots,X_{n}]$.
But this follows immediately from proposition 1.4.2 which gives a
characterization of infinitesimal bounded internal polynomials. Indeed,
for each standard $\nu$, we have $a_{\nu}\approx b_{\nu}$ and $|a_{\nu}-b_{\nu}|^{\frac{1}{|\nu|}}\approx0$
for infinite $\nu$, since $|a_{\nu}-b_{\nu}|^{\frac{1}{|\nu|}}\leq2^{\frac{1}{|\nu|}}\varepsilon\leq2\varepsilon$
for every $\nu\in{}^{*}\mathbb{N}^{n}$, $n_{0}\leq|\nu|\leq d$,
q.e.d.

\subsection{The standard part functor in the affine case}

Let $n,m$ be \emph{finite} positive integers and  $i_{b}:{}^{b}\mathbb{C}[X_{1},\ldots,X_{n}]\hookrightarrow\mathbb{C}[X_{1},\ldots,X_{n}]_{int}$
the inclusion of the algebra of bounded polynomials into the algebra
of internal polynomials.\\

\textbf{Proposition 2.5.1} \emph{Let $u$ be a {*}-homomorphism of
$^{*}\mathbb{C}$-algebras}

\emph{\[
u:\mathbb{C}[X_{1},\ldots,X_{n}]_{int}\longrightarrow\mathbb{C}[Y_{1},\ldots,Y_{m}]_{int}\]
 such that $u(X_{i})=g_{i}\in{}^{b}\mathbb{C}[Y_{1},\ldots,Y_{m}]$.
Then we have}

\emph{i) The homomorphism $u$ induces $^{b}u:\,{}^{b}\mathbb{C}[X_{1},\ldots,X_{n}]\longrightarrow\,{}^{b}\mathbb{C}[Y_{1},\ldots,Y_{m}]$,
a homomorphism of $^{b}\mathbb{C}$-algebras such that} $i_{b}\circ(^{b}u)=u\circ i_{b}$.

\emph{ii) The homomorphism} $^{b}u$ \emph{induces a map} ${\mathrm{st}}(^{b}u)$
\emph{which is a homomorphism of Stein algebras between $\mathcal{O}(\mathbb{C}^{n})$
and $\mathcal{O}(\mathbb{C}^{m})$ such that} $\mathrm{st}\circ({}^{b}u)=\mathrm{st}({}^{b}u)\circ\mathrm{st}$.

~

\emph{Proof}\textbf{.} Let $f$ be an internal polynomial of degree
at most $d{\in{}}^{*}\mathbb{N}$, so $f=\Sigma{}_{|\nu|\leq d}a_{\nu}\, X_{1}^{\nu_{1}}\ldots X_{n}^{\nu_{n}}$.
Then

\[
u(f)=\sum_{|\nu|\leq d}a_{\nu}\, g_{1}^{\nu_{1}}\ldots g_{n}^{\nu_{n}}=f(g_{1},\ldots,g_{n})\]

Assume now $f\in{}{}^{b}\mathbb{C}[X_{1},\ldots,X_{n}]$. It is clear
that \emph{$u(f)\in{}{}^{b}\mathbb{C}[Y_{1},\ldots,Y_{m}]$} since
if \emph{$x\in{}{}^{b}\mathbb{C}^{n}$} then \emph{$(g_{1}(x),\ldots,g_{n}(x))\in{}{}^{b}\mathbb{C}^{n}$}
and \emph{$f(g_{1}(x),\ldots,g_{n}(x))\in{}{}^{b}\mathbb{C}$.} \\

As usual ${}{}^{\circ}g_{1},{\ldots,}^{\circ}g_{n}$ stands for the
standard parts of $g_{1},\ldots,g_{n}$. Let $h\in\mathcal{O}(\mathbb{C}^{n})$.
We define ${\mathrm{st}}(^{b}u)$ by \[
{\mathrm{st}}(^{b}u)(h):=h({}^{\circ}g_{1},{\ldots,}^{\circ}g_{n}).\]

It is immediate that ${\mathrm{st}}(^{b}u)$ defines a homomorphism
of Stein algebras between \emph{$\mathcal{O}(\mathbb{C}^{n})$} and
\emph{$\mathcal{O}(\mathbb{C}^{m})$.} Moreover, for every \emph{$N\in{}{}^{*}\mathbb{N}_{\infty}$}
we have ${\mathrm{st}}(^{b}u)(h)={}^{\circ}(h_{N}(g_{1},\ldots,g_{n}))$,
where $h_{N}$ denotes the Taylor expansion up to order $N$ of $h$.
In other words, we get ${\mathrm{st}}(^{b}u)\circ\,\mathrm{st}=\mathrm{st}\,\circ({}^{b}u)$
on $^{b}\mathbb{C}[X_{1},\ldots,X_{n}]$ which was to show. - \\

Now we shall define the category of {}``bounded'' polynomial algebras
(over $\mathbb{C}$). Objects are given by $^{b}\mathbb{C}$-algebras
of the form $^{b}\mathbb{C}[X_{1},\ldots,X_{n}]/I$ where $n\in\mathbb{N}$
and $I$ is an arbitrary ideal of $^{b}\mathbb{C}[X_{1},\ldots,X_{n}]$.
Let $A=\,{}^{b}\mathbb{C}[X_{1},\ldots,X_{n}]/I$ and $B=\,{}^{b}\mathbb{C}[Y_{1},\ldots,Y_{m}]/J$
be two bounded polynomial algebras. A morphism between $A$ and $B$
is given by $^{b}u:\,{}^{b}\mathbb{C}[X_{1},\ldots,X_{n}]\longrightarrow\,{}^{b}\mathbb{C}[Y_{1},\ldots,Y_{m}]$,
a morphism of $^{b}\mathbb{C}$-algebras induced by a n-uplets of
bounded internal polynomials as in proposition 2.5.1 which sends the
ideal $I$ to the ideal $J$. In this way we get a category of algebras
which we call the category of \emph{bounded} \emph{polynomial} \emph{algebras}.
We note that there exist coproducts in this category.\\

We define the full subcategory of bounded polynomial algebras of finite
type where objects are given by $\,{}^{b}\mathbb{C}[X_{1},\ldots,X_{n}]/I$
where $I$ is an ideal generated by \emph{finitely} many bounded polynomials.
Let $I_{int}:=i_{b}(I).\mathbb{C}[X_{1},\ldots,X_{n}]_{int}$ denote
the ideal associated to $I$. Since the ideal $I$ is of finite type
then $I_{int}$ is an \emph{internal} ideal of $\mathbb{C}[X_{1},\ldots,X_{n}]_{int}$.
By our construction, we have the following proposition\\

\textbf{Proposition 2.5.2} \emph{There is a natural covariant functor
$\mathcal{F}$ from the} \emph{category of bounded polynomial algebras
of finite type (over $\mathbb{C}$) to the category of {*}-algebras
of finite type over $\mathbb{C}$. This functor is compatible with
coproducts.}

~

Now we have developed the necessary tools to prove the main result
of this section which is the construction of the standard part functor
from the category of bounded polynomial algebras to the category of
Stein algebras (over $\mathbb{C}$) of finite embedding dimension.
\\

\textbf{Theorem 2.5.3} \emph{There is an essentially surjective functor}
\emph{$\mathcal{ST}$, called the standard part functor, from the}
\emph{category of bounded polynomial algebras to the category of Stein
algebras (over $\mathbb{C}$) of finite embedding dimension}\[
\mathrm{\mathcal{ST}}:(\mbox{Bounded polynomial algebras)}\longrightarrow(\mbox{Stein algebras)}.\]

\emph{Proof}\textbf{.} Let $A=\,{}^{b}\mathbb{C}[X_{1},\ldots,X_{n}]/I$
be a bounded polynomial algebra and $\textrm{st}:\,^{\mathrm{b}}\mathbb{C}[X_{1},\ldots X_{n}]\longrightarrow\mathcal{O}(\mathbb{C}^{n})$
the ring epimorphism defined in proposition 2.4.4. Then $\,^{\circ}I:=\mathrm{\mathrm{st}(I)}$
is an ideal of $\mathcal{O}(\mathbb{C}^{n})$ and its closure $\overline{^{\circ}I}$
in $\mathcal{O}(\mathbb{C}^{n})$ gives us a Stein algebra \[
\mathrm{\mathcal{ST}}({}^{b}\mathbb{C}[X_{1},\ldots,X_{n}]/I):=\mathcal{O}(\mathbb{C}^{n})/\,\overline{^{\circ}I}.\]

Now, let $B={}^{b}\mathbb{C}[Y_{1},\ldots,Y_{m}]/J$ be another bounded
polynomial algebra and $^{b}u$ the morphism between $A$ and $B$
which is given by n-uplets of bounded internal polynomials. By proposition
2.5.1, ${\mathrm{st}}({}^{b}u)$ gives a homomorphism of Stein algebras
between \emph{$\mathcal{O}(\mathbb{C}^{n})$} and \emph{$\mathcal{O}(\mathbb{C}^{m})$,}
satisfying ${\mathrm{st}}({}^{b}u)\circ{}\mathrm{st}=\mathrm{st}{}\circ({}^{b}u)$.
Since ${\mathrm{st}}({}^{b}u)$ is continuous, we get

\[
{\mathrm{st}}({}^{b}u)(\overline{^{°}I})\subset\mathrm{st}(\overline{^{°}J}).\]
 Setting

\[
\mathrm{\mathcal{ST}}({}^{b}u):\mathcal{O}(\mathbb{C}^{n})/\,\overline{^{°}I}\longrightarrow\mathcal{O}(\mathbb{C}^{m})/\,\overline{^{°}J}\]
 which is given by ${\mathrm{st}}(^{b}u)$ modulo $\overline{^{\circ}I}$,
we defined $\mathcal{ST}$ on morphisms.

Now, let $C=\mathcal{O}(\mathbb{C}^{n})/\,\mathfrak{a}$ be a Stein
algebra. Then $\mathfrak{a}$ is a closed ideal of $\mathcal{O}(\mathbb{C}^{n})$
generated by a family $(g_{i})_{i\in I}$ of entire holomorphic functions
over $\mathbb{C}^{n}$\[
\mathfrak{a}=\sum_{i\in I}g_{i}\mathcal{O}(\mathbb{C}^{n}).\]

For each $i\in I$, let $f_{i}\in\,{}^{b}\mathbb{C}[X_{1},\ldots,X_{n}]$
be a bounded internal polynomial such that ${}{}^{\circ}f_{i}=g_{i}$.
Define $\mathfrak{a}_{b}$ the ideal generated by the family $(f_{i})_{i\in I}$
over $\,{}^{b}\mathbb{C}[X_{1},\ldots,X_{n}]$, that is \[
\mathfrak{a}_{b}=\sum_{i\in I}f_{i}\,{}^{b}\mathbb{C}[X_{1},\ldots,X_{n}].\]
Clearly we have ${}{}^{\circ}\mathfrak{a}_{b}=\mathfrak{a}$ which
implies that $\overline{^{\circ}\mathfrak{a}_{b}}=\mathfrak{a}$.
Hence 

~

$\mathrm{\mathcal{ST}}({}^{b}\mathbb{C}[X_{1},\ldots,X_{n}]/\mathfrak{a}_{b})=\mathcal{O}(\mathbb{C}^{n})/\,\mathfrak{a}$

~ \\
which finishes the proof. -

~\\
2.5.4 If one wants to treat also the case Stein algebras which are
\emph{not} necessarily of finite embedding dimension, one can proceed
in the following way: Let $X$ be a Stein complex space. Then there
exists an increasing sequence of natural numbers $n_{0}<n_{1}<...$
and a projective system of holomorphic maps 

$(f_{k}:X\rightarrow\mathbb{C}^{n_{k}})_{k\in\mathbb{N}}$ \\
such that

~

a) the categorical image $Y_{k}$ of $f_{k}$ is a closed complex
subspace of $\mathbb{C}^{n_{k}}$ and so we obtain a projective system
$(Y_{k})_{k}$ of Stein subspaces of $(\mathbb{C}^{n_{k}})_{k}$,

b) for each compact $K\subset X$, there is a $k$ such that $f_{k}$
is an embedding of $K$; so we may consider in particular the map
$X\rightarrow\mathrm{lim}_{\leftarrow}Y_{k}$ as a monomorphism and
$\mathrm{lim}_{\rightarrow}\mathcal{O}(Y_{k})$ as a dense subalgebra
of $\mathcal{O}(X)$ . 

~\\
To these date, we can associate an algebra of internal bounded polynomials
as follows: The sequence $(n_{0},n_{1,},...)$ defines an infinite
natural number $N$ and soforth the ring of internal polynomials $\,{}^{*}\mathbb{C}[X_{1},\ldots,X_{N}]$.
It contains $\mathrm{lim}_{\rightarrow k}{}^{b}\mathbb{C}[X_{1},\ldots,X_{n_{k}}]$
as a subring. The projective system $(Y_{k})_{k}$ defines in an obvious
(but non unique) way an ideal in this subring by fixing in addition
an infinite natural number $N'$, serving as a degree for replacing
holomorphic functions by internal polynomials. So we are able to {}``replace''
the complex space $X$ by a certain type of algebra of internal polynomials.
It is not difficult (but technical) to give a description in categorical
terms of a standard part functor (which is essentially surjective),
defined on this type of algebras and taking values in the category
of Stein algebras. Since we do not need this construction in the sequel,
we omit the details.

\subsection{Completions and enlargements}

The aim of this section is to compare the notions of completions and
enlargements of topological rings. Let $(A,\mathfrak{I})$ be a topological
ring and $\mathfrak{I}$ is an ideal of A which makes $A$ separated
and complete for the $\mathfrak{I}$ -adic topology, that is $\bigcap_{n>0}\,\mathfrak{I}^{n}=0$
and $A\cong\widehat{A}:=\lim_{\leftarrow}A/\,\mathfrak{I}^{n+1}$.
Let $(\,^{*}A,\,^{*}\mathfrak{I})$ be an enlargement of the couple
$(A,\mathfrak{I})$. First, we prove that the ring $\,^{*}A$ equipped
with the $\,^{*}\mathfrak{I}$-adic topology is in general not separated
since $\mu(0):=\bigcap_{n>0}\,^{*}\mathfrak{I}^{n}$ , the halo of
$0$, is not reduced to $\{0\}$, which occurs only if $0$ is isolated
in $A$, in other words if the ring $A$ is discrete. Let $\,^{(*)}A:=\,^{*}A/\mu(0)$
denote the separated ring associated with $\,^{*}A$. We show that
$\,^{(*)}A$ is complete for the $\mathfrak{^{*}I}$-adic topology.\\

Before giving the proof, we shall prove first some results in a more
general context and deduce from those the fact mentioned above.\\
 \\
 Let $A$ be a ring and $\,^{*}A:=\prod_{\mathcal{U}}A$ be an enlargement
of $A$, where $\mathcal{U}$ is a nonprincipal ultrafilter on $\mathbb{N}$.
Let $\mathfrak{a}$ be an \emph{internal} ideal of $^{*}A$. \\

\textbf{Theorem 2.6.1} \emph{The canonical homomorphism of rings}\[
\theta:\:^{*}A\longrightarrow\lim_{\longleftarrow}\:^{*}A/\,\mathfrak{a}^{n+1}\]
 \emph{is surjective and its kernel is $\cap{}_{n>0}\mathfrak{a}^{n}$.}\\

\emph{Proof}\textbf{.} The canonical homomorphism $\theta$ is defined
by \[
\theta(x)=(x_{k}:=x\mbox{\, mod\,}\mathfrak{a}^{k+1})\]

Let $(x_{k})$ be a sequence of elements of $^{*}A$, such that $x_{k+1}-x_{k}\in\mathfrak{a}^{k+1}$.
Since the enlargement $^{*}A$ is comprehensive, the sequence $(x_{k})$
extends to an internal sequence $(x_{k})_{k\in\,^{*}\mathbb{N}}$,
indexed by $^{*}\mathbb{N}$. Let \[
C=\{ k\in\,^{*}\mathbb{N}\,\mid\, x_{k+1}-x_{k}\in\mathfrak{a}^{k+1}\}.\]
 Then $C$ is an internal subset of $^{*}\mathbb{N}$ which contains
$\mathbb{N}$. Hence, by permanence, there exists $\omega_{0}\in\,^{*}\mathbb{N}_{\infty}$
such that \[
x_{k+1}-x_{k}\in\mathfrak{a}^{k+1}\,\mbox{ for every }\, k\in[\![1,\ldots,\omega_{0}]\!].\]
 We put $y:=x_{\omega_{0}+1}\in\,^{*}A$ and so we have \[
y-x_{k}=\sum_{l=k}^{\omega_{0}}(x_{l+1}-x_{l})\in\mathfrak{a}^{k+1}\mbox{ for every }k\in\mathbb{N}\]
 Then $\theta(y)=y\mbox{\, mod\,}\mathfrak{a}^{k+1}=x_{k}\mbox{\, mod\,}\mathfrak{a}^{k+1}$
which proves that $\theta$ is surjective. Clearly, $\mathrm{Ker}(\theta)=\cap{}_{n>0}\mathfrak{a}^{n}$.

Finally

\[
^{*}A/\cap_{n>0}\mathfrak{a}^{n}\cong\lim_{\longleftarrow}\:^{*}A/\,\mathfrak{a}^{n+1}.\]
 This shows that the separated space associated to $^{*}A$ for the
$\mathfrak{a}$-adic topology (i.e $^{*}A/\cap{}_{n>0}\mathfrak{a}^{n}$)
is complete for this topology, q.e.d.\\

\noun{Remark} 2.6.2 By the permanence \emph{}principle\emph{,} we
have \emph{$\cap{}_{k>0}\mathfrak{a}^{k}=\cup{}_{k\in\,^{*}\mathbb{N}_{\infty}}\mathfrak{a}^{k}$}
and we conclude that \emph{$\cap{}_{k>0}\mathfrak{a}^{k}\not=0$,}
if and only if\emph{, $\mathfrak{a}^{k+1}\not=0$} for each $k\in\mathbb{N}$\emph{,}
which is the case if and only if \emph{$\mathfrak{a}^{k+1}\not=0$}
for some $k\in\,^{*}\mathbb{N}_{\infty}$.

~ \\
Using elementary proprieties of projective limits, we get\\

\textbf{Corollary 2.6.3} \emph{Let $A$ be an $\mathfrak{I}$-adic
ring. Then we have }

i) $\lim_{\longleftarrow}\:^{*}A/\,^{*}\mathfrak{I}^{n+1}\cong{}^{*}A/\mu(0)$,

ii) $\widehat{A}\hookrightarrow{}^{*}A/\mu(0)$,

iii) $\widehat{A}\,\cong{}^{*}A/\mu(0)$ \emph{if and only if $A/\,\mathfrak{I}^{n+1}$
is finite for every $n\in\mathbb{N}$.}\\

\emph{Proof}\textbf{.} The first assertion is a direct consequence
of the theorem 2.6.1. For the second one, by transfer, we have $\,^{*}(A/\mathfrak{I}^{k})=\,^{*}A/\,^{*}\mathfrak{I}^{k}$.
Consider the following sequence of projective systems (with surjective
transition maps)\[
0\longrightarrow(A/\,\mathfrak{I}^{k})_{k}\longrightarrow({}^{*}A/\,^{*}\mathfrak{I}^{k})_{k}\longrightarrow(\mathrm{Coker}(*_{k}))_{k}\longrightarrow0.\]
 Taking the projective limit, we obtain the exact sequence

\[
0\longrightarrow\lim_{\longleftarrow}\: A/\,\mathfrak{I}^{k+1}\longrightarrow\lim_{\longleftarrow}\:^{*}A/\,^{*}\mathfrak{I}^{k+1}\longrightarrow\lim_{\longleftarrow}\mathrm{Coker}(*_{k+1})\longrightarrow0.\]
 Therefore we get an injective homomorphism of rings $\widehat{A}\hookrightarrow{}^{*}A/\mu(0)$
since $\widehat{A}\cong\lim_{\longleftarrow}A/\,\mathfrak{I}^{n+1}$
and $\lim_{\longleftarrow}\:^{*}A/\,^{*}\mathfrak{I}^{n+1}\cong{}^{*}A/\mu(0)$.
Furthermore, $\widehat{A}$ is (via this map) isomorphic to $^{*}A/\mu(0)$,
if and only if, $\lim_{\longleftarrow}\mathrm{Coker}(*_{k+1})=0$,
if and only if, $*_{k+1}\,:\, A/\,\mathfrak{I}^{k+1}\longrightarrow\:^{*}A/\,^{*}\mathfrak{I}^{k+1}$
is an isomorphism for each $k\in\mathbb{N}$, i.e. if and only if,
\emph{$A/\,\mathfrak{I}^{k+1}$} is finite for every \emph{$k\in\mathbb{N}$,}
q.e.d.\\

\noun{Example} 2.6.4 Let $K$ be a field and $A=K[X_{1},\ldots,X_{n}]$
be the ring of polynomials with coefficients in the field $K$ and
$\mathfrak{M}$ be the maximal ideal generated by $(X_{1},\ldots,X_{n})$.
Let $^{*}A\,=K[X_{1},\ldots,X_{n}]_{int}$ denote the ring of internal
polynomials, so $\mu(0)=\cap{}_{k>0}(X_{1},\ldots,X_{n})^{k}\,{}^{*}A$.
The ring $^{*}A/\mu(0)$ is isomorphic to a ring of power series in
the standard sense $\Sigma{}_{k\geq0}a_{\nu}\, X^{\nu}$ where $a_{\nu}\in\,^{*}K$
and the sequence $(a_{\nu})_{\nu\in\mathbb{N}^{n}}$ forms an initial
segment of a hyperfinite sequence in $\,^{*}K$. Hence, it is evident
that $K[\![X_{1},\ldots,X_{n}]\!]$ , the ring of power series, is
included in $^{*}A/\mu(0)$. Furthermore, if the field $K$ is \emph{finite}
then $K[\![X_{1},\ldots,X_{n}]\!]$ and $^{*}A/\mu(0)$ are isomorphic.

\section{Affine {*}-schemes and {*}-bounded schemes}

\subsection{Affine {*}-schemes}

We construct now the category of affine nonstandard schemes and later
that of {}``convergent'' affine nonstandard schemes (which well
call {*}-bounded schemes), more directly related to complex spaces.
Our approach is selfcontained and independent of the paper \cite{Brunjes2}.
There the authors defined the functor $^{*}\mathrm{Spec}$ from the
category $^{*}(\mathrm{Rings)}$ the category of internal rings to
$^{*}(\mathrm{Aff.\: Sch)},$ the category of internal affine schemes,
and {*}-affine schemes as the essential image of $^{*}\mathrm{Spec}$.\\

In this section, we equip in particular these objects by a topology
and a canonical sheaf structure. As a consequence {*}-affine schemes
form a subcategory of the category of locally ringed spaces.\\

Let $I$ be an infinite set and $\mathcal{U}$ be a nonprincipal ultrafilter
on $I$.

\subsubsection*{The internal spectrum of an internal ring}

Since the index set $I$ will be fixed in the sequel, we will write
$^{*}X$ instead of $^{*}X_{I}$ for the enlargement of any set $X$.\\

Let $^{*}R=\prod_{\mathcal{U}}R_{i}$ be an {*}-commutative {*}-ring.
Then $^{*}R$ is a commutative ring and $\mathrm{Spec}({}^{*}R)$
denotes the (usual) spectrum of $^{*}R$. We define

\[
\mathrm{Specint}({}^{*}R)=\left\{ \mathfrak{J}\in\mathrm{Spec}({}^{*}R)\,\mid\,\mathfrak{J}\,\,\mathrm{is\, an\, internal\, ideal\, of\,}{}^{*}R\right\} .\]

By transfer, we easily prove

i) $^{*}R$ is an integral domain if and only if $\left\{ i\in I\,\mid\, R_{i}\mathrm{\,\, is\, an\, integral\, domain}\right\} \in\mathcal{U}$,

ii) $^{*}R$ is a field if and only if $\left\{ i\in I\,\mid\, R_{i}\mathrm{\,\, is\, a\, field}\right\} \in\mathcal{U}$.

As a consequence, we get\\

\textbf{Proposition 3.1.1} \emph{With the above notations, we have
natural bijections}

i) $\mathrm{Specint}({}^{*}R)\cong\prod_{\mathcal{U}}\mathrm{Spec}(R_{i})$,

ii) $\mathrm{Specmaxint}({}^{*}R)\cong\prod_{\mathcal{U}}\mathrm{Specmax}(R_{i})$.
\\

In particular, let $k$ be an algebraic closed field, then by the
Hilbert Nullstellensatz, we have $\mathrm{Specmaxint}(k[T_{1},\ldots,T_{n}]_{int})\simeq{}^{*}k^{n}$.\\

Let $j:\,\mathrm{Specint}({}^{*}R)\,\longrightarrow\,\mathrm{Spec}({}^{*}R)$
denote the inclusion map. \\

We endow $\mathrm{Specint}({}^{*}R)$ with the induced Zariski topology,
defined on $\mathrm{Spec}({}^{*}R)$. Hence closed subsets of $\mathrm{Specint}({}^{*}R)$
are given by

\[
\mathcal{V}(\mathfrak{a})=\left\{ \mathfrak{J}\in\mathrm{Specint}({}^{*}R)\,\mid\,\mathfrak{J}\supset\mathfrak{a}\right\} \]
 where $\mathfrak{a}$ is an ideal of $^{*}R$ which may be external.\\

Let $\mathfrak{J}$ be an internal ideal of $^{*}R$. We set

\[
\sqrt[{int}]{\mathfrak{J}}=\left\{ f\in\,^{*}R\,\mid\,\exists\, n\in\,^{*}\mathbb{N}\,,\, f^{n}\in\mathfrak{J}\right\} \]
 Then $\sqrt[{int}]{\mathfrak{J}}$ is an ideal of $^{*}R$, containing
$\mathfrak{J}$.\\

If $\mathfrak{J}$ is an internal ideal of $^{*}R$, so $\mathfrak{J}=\prod_{\mathcal{U}}\mathfrak{J}_{i}$,
then $\sqrt[{int}]{\mathfrak{J}}$ is internal too and given by $\sqrt[{int}]{\mathfrak{J}}=\prod_{\mathcal{U}}\sqrt{\mathfrak{J}_{i}}$.
Again by transfer, we conclude that $\sqrt[{int}]{\mathfrak{J}}$
is the intersection of all internal prime ideals containing $\mathfrak{J}$,
i.e. \emph{\begin{eqnarray*}
\sqrt[{int}]{\mathfrak{J}} & = & {\displaystyle {\displaystyle \bigcap_{\mathfrak{p}\in\mathrm{Specint}({}^{*}R),\,\mathfrak{p}\supset\mathfrak{J}}}\mathfrak{p}}.\end{eqnarray*}
 }

We have\\

\textbf{Proposition 3.1.2} \emph{i) Let $\mathfrak{a}=\Pi{}_{\mathcal{U}}\mathfrak{a}_{i}$
be an internal ideal of $^{*}R$. Then $\mathcal{V}(\mathfrak{a})$
is an internal subset of $\mathrm{Specint}({}^{*}R)$ and $\mathfrak{\mathcal{V}(a)}=\Pi{}_{\mathcal{U}}V(\mathfrak{a}_{i})$
where $V(\mathfrak{a}_{i})=\left\{ \mathfrak{J}\in\mathrm{Spec}(R_{i})\,\mid\,\mathfrak{J}\supset\mathfrak{a}_{i}\right\} $.}

\emph{ii) Let $\mathfrak{a},\mathfrak{b}$ be two internal ideals
in $^{*}R$. Then $\sqrt[{int}]{\mathfrak{b}}\subset\sqrt[{int}]{\mathfrak{a}}$
if and only if $\mathcal{V}(\mathfrak{a})\subset\mathcal{V}(\mathfrak{b})$.}\\

\emph{Proof.} i) Let \emph{$\mathfrak{J}=\Pi{}_{\mathcal{U}}\mathfrak{J}_{i}$}
be an internal ideal of \emph{$^{*}R$}. By transfer, we have \emph{$\mathfrak{J}\supset\mathfrak{a}$}
if and only if $\left\{ i\in I\,\mid\,\mathfrak{J}_{i}\supset\mathfrak{a}_{i}\right\} \in\mathcal{U}$
which is equivalent to $\left\{ i\in I\,\mid\,\mathfrak{J}_{i}\in V(\mathfrak{a}_{i})\right\} \in\mathcal{U}$.

ii) Let $\mathfrak{p}$ be an internal ideal containing $\mathfrak{a}$.
Then $\mathfrak{p}$ contains $\sqrt[{int}]{\mathfrak{a}}$, hence
$\mathfrak{p}\supset\sqrt[{int}]{\mathfrak{b}}\supset\mathfrak{b}$.
The converse is, via transfer, an easy consequence, since \emph{$\mathcal{V}(\mathfrak{a})\subset\mathcal{V}(\mathfrak{b})$}
if an only if $\left\{ i\in I\,\mid\, V(\mathfrak{a}_{i})\subset V(\mathfrak{b}_{i})\right\} \in\mathcal{U}$
if and only if $\left\{ i\in I\,\mid\,\sqrt{\mathfrak{b}_{i}}\subset\sqrt[{}]{\mathfrak{a}_{i}}\right\} \in\mathcal{U}$
which means \emph{$\sqrt[{int}]{\mathfrak{b}}\subset\sqrt[{int}]{\mathfrak{a}}$},
q.e.d.\\

We fix $f=\Pi{}_{\mathcal{U}}f_{i}\in\,^{*}R$, and put \[
\mathfrak{D}(f)=\left\{ \mathfrak{p}\in\mathrm{Specint}({}^{*}R)\,\mid\, f\not\in\mathfrak{p}\right\} .\]

Then $\mathfrak{D}(f)=\Pi{}_{\mathcal{U}}D(f_{i})$, where $D(f_{i})=\left\{ \mathfrak{p}\in\mathrm{Spec}(R_{i})\,\mid\, f_{i}\not\in\mathfrak{p}\right\} $.
The sets $\left\{ \mathfrak{D}(f)\,\mid\, f\in\,^{*}R\right\} $ form
an internal open base for the induced Zariski topology, induced by
$\mathrm{Spec}({}^{*}R)$. For any ideal $\mathfrak{I}$ in $^{*}R$,
we have $\mathrm{Specint}(^{*}R)\setminus\mathcal{V}(\mathfrak{I})=\cup_{f\in\mathfrak{I}}\mathfrak{D}(f)$.

Since $\mathfrak{D}(f)$ are internal subsets of $\mathrm{Specint}({}^{*}R)$,
we have\\

\textbf{Proposition 3.1.3} \emph{Let $f=\Pi{}_{\mathcal{U}}f_{i},\,$
and $g=\Pi{}_{\mathcal{U}}g_{i}$ be two elements of $^{*}R$. Then }

\emph{i) $\mathfrak{D}(f)=\emptyset$ if and only if there exists
$n\in\,^{*}\mathbb{N}$ such that $f^{n}=0$,}

\emph{ii) $\mathfrak{D}(f)\cap\mathfrak{D}(g)=\mathfrak{D}(fg)$ and
for each $n\in\,^{*}\mathbb{N}$ positive, we have $\mathfrak{D}(f^{n})=\mathfrak{D}(f)$,}

\emph{iii) $\mathfrak{D}(f)\supset\mathfrak{D}(g)$ if and only if
$g\in\sqrt[{int}]{(f)}$, which is equivalent to}\\
\emph{$\left\{ i\in I\,\mid\, D(f_{i})\supset D(g_{i})\right\} \in\mathcal{U}$.}\\

Moreover, we have \\

\textbf{Proposition 3.1.4} \emph{Let $^{*}S=\Pi{}_{\mathcal{U}}S_{i}$
be an internal subset of} $^{*}R$. \emph{Then $^{*}S$ is a multiplicative
subset of $^{*}R$ if and only if} $\left\{ i\in I\,\mid\, S_{i}\,\mathrm{is\, multiplicative\, in\,}R_{i}\right\} \in\mathcal{U}$
and $^{*}S^{-1}\left(^{*}R\right)=\Pi{}_{\mathcal{U}}S_{i}^{-1}R_{i}$.\\

We consider two examples of internal multiplicative subsets of $^{*}R$\\

i) Let $\mathfrak{p}=\Pi{}_{\mathcal{U}}\mathfrak{p}_{i}\in\mathrm{Specint}({}^{*}R)$
be an internal prime ideal of $^{*}R$. Then $^{*}R\setminus\mathfrak{p}=\Pi{}_{\mathcal{U}}R_{i}\setminus\mathfrak{p}_{i}$
is an internal multiplicative subset of $^{*}R$ and $\left(^{*}R\right)_{\mathfrak{p}}=\,^{*}S^{-1}({}^{*}R)=\Pi{}_{\mathcal{U}}\left(R_{i}\right)_{\mathfrak{p}_{i}}$.

ii) Let \emph{$f=\Pi{}_{\mathcal{U}}f_{i}\in\,^{*}R$} and \emph{$^{*}S=\{1,f,f^{2,},\ldots,f^{N},\ldots,\, N\in\,^{*}\mathbb{N}\}$.}
Clearly $^{*}S$ is an internal multiplicative subset of $^{*}R$
and $^{*}S=\Pi{}_{\mathcal{U}}S_{i}$, where $S_{i}=\{1,f_{i},f_{i}^{2,},\ldots,f_{i}^{n},\ldots,\, n\in\mathbb{N}\}$.
We denote by\[
^{*}R_{[f]}=\,^{*}S^{-1}\left(^{*}R\right)=\Pi{}_{\mathcal{U}}\left(R_{i}\right)_{f_{i}}\]
 the localization of $^{*}R$ with respect to $^{*}S$ which is an
internal ring. Its internal prime spectrum is given by\\

\textbf{Proposition 3.1.5} \emph{Let $f\in\,^{*}R$. Then} $\mathrm{Specint}({}^{*}R_{[f]})=\mathfrak{D}(f)=\mathrm{Specint}({}^{*}R)\setminus\mathcal{V}(f$).

~

\emph{Proof.} We identify naturally

$\mathrm{Specint}({}^{*}R_{[f]})=\mathrm{Specint}(\Pi{}_{\mathcal{U}}(R_{i})_{f_{i}})=\Pi{}_{\mathcal{U}}\mathrm{Spec}((R_{i})_{f_{i}})=\Pi{}_{\mathcal{U}}D(f_{i})=\mathfrak{D}(f)$
. -\\

Let $^{*}R_{f}$ denote the localization of the ring $^{*}R$ with
respect to the multiplicative family \emph{$\{1,f,f^{2,},\ldots,f^{n},\ldots,\, n\in\,\mathbb{N}\}$}.
Then there is a natural morphism \[
^{*}R_{f}\,\longrightarrow\,^{*}R_{[f]}\]
 of rings, induced by the identity.

\subsubsection*{Structure sheaf of an internal prime spectrum}

Following the classical procedure, we define a sheaf of commutative
rings over $^{*}X:=\mathrm{Specint}({}^{*}R)$, the internal prime
spectrum, equipped with the internal topology of Zariski. We first
define sections and restriction maps on the sets $\mathfrak{D}(f)$,
$f\in\,^{*}R$, which form a base for the Zariski topology, induced
by that of $\mathrm{Spec}({}^{*}R)$. We set\[
\mathcal{A}_{^{*}X}\left(\mathfrak{D}(f)\right):={}^{*}R_{[f]}\]

Clearly, this defines a presheaf on $^{*}X$, where restriction maps
on elements of the base are given as follows: Let \emph{$\mathfrak{D}(f)\supset\mathfrak{D}(g)$}
which is equivalent to $J:=\{ i\in I\,\mid\, D(f_{i})\supset D(g_{i})\}\in\mathcal{U}$.
There is an internal homomorphism of internal rings \[
^{*}\rho_{\mathfrak{D}(g),\mathfrak{D}(f)}:\,{}^{*}R_{[f]}\longrightarrow\,{}^{*}R_{[g]}\]
 induced by the restriction maps $\rho_{D(g_{i}),D(f_{i})}:\,(R_{i})_{f_{i}}\longrightarrow(R_{i})_{g_{i}}$
for every $i\in J$. Trivially, we have $^{*}\rho_{\mathfrak{D}(f),\mathfrak{D}(f)}=id$
and $^{*}\rho_{\mathfrak{D}(h),\mathfrak{D}(g)}\circ{}^{*}\rho_{\mathfrak{D}(g),\mathfrak{D}(f)}={}^{*}\rho_{\mathfrak{D}(h),\mathfrak{D}(f)}$,
for $\mathfrak{D}(f)\supset\mathfrak{D}(g)\supset\mathfrak{D}(h)$.
Consider the collection of internal open sets $\mathfrak{D}(f)$,
containing $\mathfrak{p}\in\mathrm{Specint}({}^{*}R)$. We put $\mathcal{W}_{\mathfrak{p}}=\{\mathfrak{D}(f),\,\mathfrak{p}\in\mathfrak{D}(f)\}$.
It is a directed set.

Let $^{*}\rho_{\mathfrak{p}}^{f}:\,{}^{*}R_{[f]}\longrightarrow\,{}^{*}R_{\mathfrak{p}}$
be the canonical internal homomorphism for each $\mathfrak{p}\in\mathfrak{D}(f)$.
The following proposition is immediate \\

\textbf{Proposition 3.1.6} \emph{For an internal prime ideal $\mathfrak{p}$
of an internal ring $^{*}R$, there is a natural isomorphism of rings}

\[
\lim_{\rightarrow}\,^{*}R_{[f]}\rightarrow\,^{*}R_{\mathfrak{p}}\]
 \emph{where} \emph{the} \emph{limit} \emph{is} \emph{taken} \emph{over}
\emph{all} $f$ \emph{such} \emph{that} $\mathfrak{p}\in\mathfrak{D}(f)$.

~

\emph{Proof}\noun{.} We can use transfer to prove this by constructing
an isomorphism between $\lim_{\rightarrow}\,^{*}R_{[f]}$ and $\Pi{}_{\mathcal{U}}\lim_{\rightarrow}\,\left(R_{i}\right)_{f_{i}}$.
The last one is internally isomorphic to $\Pi{}_{\mathcal{U}}\left(R_{i}\right)_{\mathfrak{p}_{i}}=\,^{*}R_{\mathfrak{p}}$,
where $\mathfrak{p}=\Pi{}_{\mathcal{U}}\mathfrak{p}_{i}\in\mathfrak{D}(f)$,
so we are done.\\

We show now, using transfer, that our sheaf is already separated:\\

\textbf{Lemma 3.1.7} \emph{Let $f\in\,^{*}R$ and} $\mathfrak{D}(f)=\cup_{\alpha\in A}\mathfrak{D}(f_{\alpha})$\emph{.
Suppose that for} $a\in\,^{*}R_{[f]}$ \[
^{*}\rho_{\mathfrak{D}(f_{\alpha}),\mathfrak{D}(f)}(a)=0\,\mathrm{\quad\mathit{for\, each\,}}\alpha\in A.\]
 \emph{Then} $a=0$.

~

\emph{Proof.} The element \emph{$a\in\,^{*}R_{[f]}$} can be expressed
in the form \emph{$a=g/f^{m}$}, where \emph{$m$} is a hyperinteger.
If \emph{$\mathfrak{J=}\left\{ h\in\,^{*}R\,\mid\, hg=0\right\} $},
then \emph{$\mathfrak{J}$} is an internal ideal of \emph{$^{*}R$}
and \emph{$a=0$} in \emph{$^{*}R_{[f]}$} if and only if \emph{$f\in\sqrt[{int}]{\mathfrak{J}}={\displaystyle \cap{\displaystyle _{\mathfrak{J}\subset p\in\mathrm{Specint}({}^{*}R)}}\mathfrak{p}}$}.

Suppose that $a\not=0$ in \emph{$^{*}R_{[f]}$}. Then there exists
an internal prime ideal $\mathfrak{p}\supset\mathfrak{J}$ with $f\not\in\mathfrak{p}$,
i.e. $\mathfrak{p}\in\mathfrak{D}(f)$. We take $\alpha\in A$ such
that $\mathfrak{p}\in\mathfrak{D}(f_{\alpha})$ and $^{*}\rho_{\mathfrak{p}}^{f_{\alpha}}\circ{}^{*}\rho_{\mathfrak{D}(f_{\alpha}),\mathfrak{D}(f)}={}^{*}\rho_{\mathfrak{p}}^{f}$.
By assumption, the image of $g=f^{m}a$ in $^{*}R_{\mathfrak{p}}$
is zero, which means, there exists $b\in R\setminus\mathfrak{p}$
such that $bg=0$. Hence $b\in\mathfrak{J}$, which contradicts $b\in R\setminus\mathfrak{p}$,
since $\mathfrak{p}\supset\mathfrak{J}$, q.e.d.\\

\textbf{Definition 3.1.8} The \emph{}sheaf \emph{}of \emph{}rings
\emph{}associated \emph{}to \emph{}the \emph{}presheaf $\mathcal{A}_{^{*}X}$,
given \emph{}by $\mathcal{A}_{^{*}X}\left(\mathfrak{D}(f)\right):={}^{*}R_{[f]}$
on \emph{}the \emph{}basis $\left\{ \mathfrak{D}(f)\,\mid\, f\in\,{}^{*}X\right\} $
of $^{*}X$\emph{,} will \emph{}be \emph{}called \emph{}the \emph{structure}
\emph{sheaf} on \emph{}the \emph{}internal \emph{}spectrum \emph{$^{*}X=\mathrm{Specint}({}^{*}R)$,}
and \emph{}denoted \emph{}by $\mathcal{O}_{^{*}X}$.\\

Let $\theta_{f}:\,{}^{*}R_{[f]}\longrightarrow\Gamma(\mathfrak{D}(f),\mathcal{O}_{^{*}X})$
be the canonical map. It is given as follows: For $\mathfrak{p}\in\mathfrak{D}(f)$
and $\xi\in\,{}^{*}R_{[f]}$, we have $\left(\theta_{f}\left(\xi\right)\right)_{\mathfrak{p}}={}^{*}\rho_{\mathfrak{p}}^{f}\left(\xi\right)$.
So by (3.1.7), we get\\

\textbf{Corollary 3.1.9} \emph{For every $f\in\,^{*}R$, the homomorphism}
$\theta_{f}:\,{}^{*}R_{[f]}\longrightarrow\Gamma(\mathfrak{D}(f),\mathcal{O}_{^{*}X})$
\emph{is injective}. \\

Let $^{*}R_{f}\,\longrightarrow\,\Gamma(\mathfrak{D}(f),\mathcal{O}_{^{*}X})$
be the homomorphism of rings, given by composing $^{*}R_{f}\,\longrightarrow\,{}^{*}R_{[f]}$
and $\theta_{f}$. These homomorphisms induce a morphism of sheaves\\
 \[
\mathcal{O}_{\mathrm{Spec}({}^{*}R)}\,\longrightarrow\, j_{*}\mathcal{O}_{\mathrm{Specint}(^{*}R)}\]
 and so we obtain $j$ as a morphism of ringed spaces. \\

\textbf{Proposition 3.1.10} \emph{The sheaves} $\mathcal{O}_{^{*}X}$
\emph{and} $j^{-1}(\mathcal{O}_{\mathrm{Spec}({}^{*}R)})$ \emph{are
naturally isomorphic.}\\

\emph{Proof.} Applying the functor $j^{-1}$ to $\mathcal{O}_{\mathrm{Spec}({}^{*}R)}\,\longrightarrow\, j_{*}\mathcal{O}_{^{*}X}$,
we get a morphism of sheaves $j^{-1}\mathcal{O}_{\mathrm{Spec}({}^{*}R)}\,\longrightarrow\,\mathcal{O}_{^{*}X}$
on $\mathrm{Specint}({}^{*}R)$. This morphism is an isomorphism since
it gives the identity on stalks at every internal prime ideal of $^{*}R$,
q.e.d.\\

The interplay between the divers spectra can be summarized as follows
: \\
 Let $R$ be a ring and we denote by $^{*}R=\Pi_{\mathcal{U}}R$ its
ultrapower. The morphism of rings $*:\, R\longrightarrow\,^{*}R$
induces a continuous map

~

$\mathrm{Spec}(*):\,\mathrm{Spec}({}^{*}R)\longrightarrow\mathrm{Spec}(R)$

~ \\
 where $\mathrm{Spec}({}^{*}R)$ and $\mathrm{Spec}(R)$ are equipped
with the Zariski topology. We denote by

~

$\mathrm{Specint}(*):\,\mathrm{Specint}\left(^{*}R\right)\longrightarrow\mathrm{Spec}(R)$

~ \\
 the restriction of $\mathrm{Spec}(*)$ to $\mathrm{Specint}({}^{*}R)$
which carries the induced Zariski topology of $\mathrm{Spec}({}^{*}R)$.
Clearly, $\mathrm{Specint}(*)$ is continuous. We should mention that
the induced topology on $\mathrm{Specint}({}^{*}R)=\,^{*}\mathrm{Spec}(R$)
coincides with the so-called Q-topology since they have the same basis,
given by the {*}-open subsets $\mathfrak{D}(f)=\Pi{}_{\mathcal{U}}D(f_{i})$.
The enlargement construction gives an injective map

~

$*:\,\mathrm{Spec}(R)\longrightarrow\,^{*}\mathrm{Spec}(R)$

~\\
 $\mathfrak{p}\mapsto*\mathfrak{p}$ . This map is \emph{not continuous}
since the restriction of the Q-topology does not coincides with the
initial one. So, it is more natural to work with $\mathrm{Specint}(*)$
on the level of topological spaces.

Now, we will define the category of affine {*}-schemes. An object
is given by a locally ringed space ${({}}^{*}X,\mathcal{O}_{^{*}X})$
with the underlying topological space $^{*}X:=\mathrm{Specint{({}}^{*}R)}$
where $^{*}R$ is an internal ring and $\mathcal{O}_{^{*}X}=j_{*}\mathcal{O}_{\mathrm{Spec}{({}}^{*}R)}$
(the mapping $j:\mathrm{Specint}{({}}^{*}R)\rightarrow\mathrm{Spec}{({}}^{*}R)$
denotes the inclusion). Let ${({}}^{*}Y,\mathcal{O}_{^{*}Y})$ be
another object in the category of {*}-affine schemes defined by an
internal ring $^{*}S$, that is, $^{*}Y=\mathrm{Specint}{({}}^{*}S)$.
A morphism between ${({}}^{*}X,\mathcal{O}_{^{*}X})$ and ${({}}^{*}Y,\mathcal{O}_{^{*}Y})$
is represented by an internal morphism of rings $\varphi:\,^{*}S{\rightarrow{}}^{*}R$
which induces a morphism of locally ringed space as follows: the restriction
of $\mathrm{Spec}(\varphi)$ on $\mathrm{Specint}{({}}^{*}R)$ gives
a continuous mapping\[
\mathrm{Specint}(\varphi):\mathrm{Specint}{({}}^{*}R)\rightarrow\mathrm{Specint}{({}}^{*}S)\]
and, by localization, $\varphi$ gives a morphism of sheaves \[
\mathcal{O}_{^{*}Y}\rightarrow\mathrm{Specint}(\varphi)_{*}(\mathcal{O}_{^{*}X}).\]

\subsection{Affine {*}-bounded schemes}

It is well known that there is a natural functor between the category
of affine schemes of finite type over $\mathbb{C}$ and the category
of Stein spaces. This functor associates for each algebra of finite
type $A=\mathbb{C}[X_{1},\ldots,X_{n}]/(f_{1},\ldots,f_{q})$ the
Stein algebra $\mathcal{O}(\mathbb{C}^{n})/\sum_{i=1}^{q}f_{i}\mathcal{O}(\mathbb{C}^{n})$,
where each polynomial $f_{i}$ is considered as an entire holomorphic
function over $\mathbb{C}^{n}$. \\

In this section, we will construct the category of affine bounded
{*}-schemes as an {}``intermediate'' category between the category
of affine schemes of finite type over $\mathbb{C}$ and the category
of Stein spaces. This category will be a subcategory of locally ringed
spaces. \\

Before giving this construction, we recall some facts and notations
on the topology of $^{b}\mathbb{C}^{n}$, the space of bounded points
of $^{*}\mathbb{C}^{n}$. The space $^{b}\mathbb{C}^{n}$ is a S-closed
and S-open subspace of $^{*}\mathbb{C}^{n}$. A basis for its topology
is given by the S-balls $S(p,r)=\{ q\in{}^{b}\mathbb{C}^{n}\,|\,{}°|q-p|<r\}$,
where $p\in{}{}^{b}\mathbb{C}^{n}$ and $r$ is a positive real standard
number. As a consequence, the standard part mapping \[
\mathrm{st}:{}{}^{b}\mathbb{C}^{n}\rightarrow{}\mathbb{C}^{n}\]
 is continuous and open. Furthermore, we note that the inclusion mapping
\[
i:{}\mathbb{C}^{n}\rightarrow{}^{b}\mathbb{C}^{n}\]
 is continuous. If fact, the S-topology coincides with the initial
topology which makes the standard mapping $\mathrm{st}:{}{}^{b}\mathbb{C}^{n}\rightarrow{}\mathbb{C}^{n}$
continuous. Finally, we note that if $U$ is S-open, then ${}{}^{\circ}U$
is an open subset of $\mathbb{C}^{n}$ and ${}{}^{\circ}U=U\cap\mathbb{C}^{n}$.
\\

Let $\mathfrak{a}$ be an ideal of ${}{}^{b}\mathbb{C}[X_{1},\ldots,X_{n}]$.
Define the \emph{variety} \emph{of} $\mathfrak{a}$ by\[
V(\mathfrak{a})=\{ x\in{}{}^{b}\mathbb{C}^{n}\,|\, f(x)\approx0\,\mbox{ for each }f\in\mathfrak{a}\}.\]
 Clearly, $V(\mathfrak{a})$ is S-closed, since each bounded polynomial
is S-continuous. Again, by S-continuity, we have $V(\mathfrak{a})=\{ x\in{}{}^{b}\mathbb{C}^{n}\,|\, f(y)\approx0\,\mbox{ for each }f\in\mathfrak{a}\,,\, y\in\mu(x)\}$.\\

Let $\mathfrak{b}$ be a subset of $\mathcal{O}(\mathbb{C}^{n})$.
Define the \emph{zero} \emph{set} \emph{of} $\mathfrak{b}$ by \[
Z(\mathfrak{b})=\{ x\in\mathbb{C}^{n}\,|\, f(x)=0\,\mbox{ for each }f\in\mathfrak{b}\}.\]
 For each ideal $\mathfrak{a}$ of $^{b}\mathbb{C}[X_{1},\ldots,X_{n}]$,
we set $°\mathfrak{a}=\mathrm{st}(\mathfrak{a})$ which is an ideal
of $\mathcal{O}(\mathbb{C}^{n})$, its standard part. \\

We have the following rules

~

\textbf{Proposition 3.2.1} \emph{Let $(\mathfrak{a}_{i})_{i\in I}$
be a family of ideals and} $\mathfrak{a}$\emph{,}$\mathfrak{b}$
\emph{two ideals all in} \emph{${}{}^{b}\mathbb{C}[X]$. Moreover,
we fix $f\in{}{}^{b}\mathbb{C}[X]$. Then }

\emph{i)} $\cap_{i\in I}V(\mathfrak{a}_{i})=V(\sum_{i\in I}\mathfrak{a}_{i})$,

\emph{ii)} $V(\mathfrak{a})\cup V(\mathfrak{b})=V(\mathfrak{a}\cap\mathfrak{b})=V(\mathfrak{a}\mathfrak{b})$,

\emph{iii)} $^{\circ}V(\mathfrak{a})=Z{({}}^{\circ}\mathfrak{a})$
\emph{and} $i^{-1}(V(\mathfrak{a}))=Z({}{}^{\circ}\mathfrak{a})$,

\emph{iv)} $V(\mathfrak{a}+{}^{i}\mathbb{C}[X])=V(\mathfrak{a})$,

\emph{v)} $V(\mathfrak{a})\subset V(\mathfrak{b})$ \emph{if and only
if} $Z{({}}^{\circ}\mathfrak{a})\subset Z{({}}^{\circ}\mathfrak{b})$,

\emph{vi) $V(f)=\emptyset$ if and only if ${}{}^{\circ}f(x)\not=0$
for every $x\in\mathbb{C}^{n}$,}

\emph{vii) $V(f)={}{}^{b}\mathbb{C}^{n}$ if and only if $f{\in{}}^{i}\mathbb{C}[X]$.}\\

\emph{Proof}. The assertions (i), (ii) , (iv), (vi) and (vii) are
straightforward. For the assertion (iii), the inclusion $Z(°\mathfrak{a})\subset{}°V(\mathfrak{a})$
is obvious. Conversely, since each element $f\in\mathfrak{a}$ is
S-continuous, we have $°V(\mathfrak{a})\subset Z(°\mathfrak{a})$.
The same argument about S-continuity together with (iii) prove (v).
\\

\noun{Remark} 3.2.2 The sets $V(\mathfrak{a})$ where $\mathfrak{a}$
is an arbitrary ideal in \emph{${}{}^{b}\mathbb{C}[X]$,} form the
closed subsets for a Zariski-topology on ${}{}^{b}\mathbb{C}^{n}$.
From the assertion (iii) of the proposition 3.2.1, we deduce that
the inclusion $i:{}\mathbb{C}^{n}\rightarrow{}^{b}\mathbb{C}^{n}$
is continuous if both spaces are equipped with the Zariski-topology.
\-\\

Now, we construct the structure sheaf of $^{b}\mathbb{C}^{n}$. Let
$U$ be a nonempty subset of $^{b}\mathbb{C}^{n}$. Define\[
B(U):=\{ f\in\mathbb{C}[X_{1},\ldots,X_{n}]_{int}\,|\, f(x)\in{}{}^{b}\mathbb{C}\mbox{ for each }x\in U\}\]
 and

\[
S(U):=\{ f\in\mathbb{C}[X_{1},\ldots,X_{n}]_{int}\,|\, f(x)\in{}{}^{a}\mathbb{C}\mbox{ for each }x\in U\}\]
 where $^{a}\mathbb{C}$ denotes the multiplicative set of appreciable
elements of $^{b}\mathbb{C}$, i.e.

\[
^{a}\mathbb{C}=\{ x\in{}{}^{b}\mathbb{C}\,|\, x\not\approx0\}.\]

Trivially, $B(U)$ is an integral domain. We will write $S^{-1}B(U)$
instead of $S^{-1}(U)B(U)$. The correspondence $U\mapsto S^{-1}B(U)$
gives a presheaf of rings over $^{b}\mathbb{C}^{n}$, where $U$ runs
over all S-open subsets in $^{b}\mathbb{C}^{n}$. The restriction
maps are evident. Such a presheaf is clearly separated. We will denote
by $\mathcal{O}_{^{b}\mathbb{C}^{n}}$ its associated sheaf.\\

Let $p\in{}^{b}\mathbb{C}^{n}$ and $U$ be a S-open subset of \emph{$^{b}\mathbb{C}^{n}$,}
containing $p$. Hence $U$ also contains $\mu(p)$. We have a canonical
morphism \[
\rho_{U}^{p}{}:\, S^{-1}B(U)\,\rightarrow\, S^{-1}B(\mu(p))\]
 defined by restriction on $\mu(p)$. 

~

We shall describe the stalks of our structure sheaf\\

\textbf{Lemma 3.2.3} \emph{Let $p\in{}{}^{b}\mathbb{C}^{n}$ and $g\in{}^{*}\mathcal{O}(\mathbb{C}^{n})$.}

\emph{i) If $g(\mu(p))\subset{}^{b}\mathbb{C}$, there exists a S-open
subset $V$, containing $p$, such that $g(V)\subset{}{}^{b}\mathbb{C}$.}

\emph{ii) If $g(\mu(p))\subset{}^{a}\mathbb{C}$, there exists a S-open
subset $V$, containing $p$, such that $g(V)\subset{}{}^{a}\mathbb{C}$.}\\

\emph{Proof}. The first assertion is an easy consequence of permanence
principle. The second assertion is deduced from theorem 2.2.1 which
asserts that there exists a S-open subset $V$ such that $g$ is S-continuous
on $V$ and has a standard part $^{\circ}g$, which is actually holomorphic
and $^{\circ}g(z)\not=0$ for every $z\in{}{}^{\circ}V$. In particular,
we conclude that $g$ takes appreciable values on $V$, q.e.d.\\

Applying this lemma for internal polynomials, we prove\\

\textbf{Proposition 3.2.4} \emph{For every $p\in{}{}^{b}\mathbb{C}^{n}$,
there is a natural isomorphism }

\emph{\[
\mathcal{O}_{^{b}\mathbb{C}^{n},\, p}\,\rightarrow\, S^{-1}B(\mu(p))\]
 between the stalk of} $\mathcal{O}_{^{b}\mathbb{C}^{n}}$ \emph{in
$p$ and $S^{-1}B(\mu(p))$.}\\

\emph{Proof.} For each S-open $U$ containing $p$, the morphisms
$\rho_{U}^{p}{}$ are compatible with restriction maps and induce,
by taking the inductive limit, the morphism \[
\lim_{\rightarrow}\, S^{-1}B(U)\,\rightarrow S^{-1}B(\mu(p))\]
 which is actually an isomorphism: the permanence principle gives
the injectivity and surjectivity is a consequence of lemma 3.2.3.
- \-\\

\noun{Remark 3.2.5} i) Trivially, the assertion in lemma 3.2.3 is
false for $^{i}\mathbb{C}$, so we have to distinguish between \[
\mathfrak{m}_{p}:=\{\frac{f}{g}\in S^{-1}B(\mu(p))\,|\, f(p)\in{}{}^{i}\mathbb{C}^{n}\}\]
 maximal ideal of $S^{-1}B(\mu(p))$ and the ideal of infinitesimal
elements\[
\mathcal{I}nf_{p}:=\{\frac{f}{g}\in S^{-1}B(\mu(p))\,|\,\exists\mathrm{\,}V\mathrm{,\,\, S-open\,:\,}p\mathrm{\in V,}\, f(z)\in{}{}^{i}\mathbb{C}^{n},\,\forall z\in V\}.\]
 We have $S^{-1}B(\mu(p))/\mathcal{I}nf_{p}\cong\mathcal{O}_{\mathbb{C}^{n},°p}$.

ii) Let $p{\in{}}^{b}\mathbb{C}^{n}$ and $f=\sum_{|\nu|\leq d}b_{\nu}(X-p)^{\nu}$
be an internal polynomial. If $b_{0}$ and $|b_{\nu}|^{\frac{1}{|\nu|}}$
are bounded for each $\nu$ such that $0<|\nu|\leq d$, then $f\in B(\mu(p))$.
This is a consequence of the estimate, used in the proof of 1.4.3.\-\\

Let $U$ be a S-open. From theorem 2.2.1 we deduce that each element
$f\in B(U)$ (resp. $f\in S(U)$ ) has a standard part mapping ${}{}^{\circ}f\in\mathcal{O}_{\mathbb{C}^{n}}{({}}^{\circ}U)$
(resp. ${}{}^{\circ}f\in\mathcal{O}_{\mathbb{C}^{n}}{({}}^{\circ}U)$,
${}{}^{\circ}f(x)\not=0$ for each $x\in{}°U$). Using the fact that
${}{}^{\circ}U=U\cap\mathbb{C}^{n}$, we get a morphism of presheaves
$\mathrm{st_{U}:}\, S^{-1}B(U)\,\rightarrow\,\mathcal{O}_{\mathbb{C}^{n}}(U\cap\mathbb{C}^{n})$,
given by $\mathrm{st}_{U}(f/h)=°f/\,°h$ and which induces a morphism
of sheaves \[
\mathrm{st_{U}:}\,\mathcal{O}_{^{b}\mathbb{C}^{n}}(U)\,\rightarrow\,\mathcal{O}_{\mathbb{C}^{n}}(U\cap\mathbb{C}^{n})\]
 Hence, we have the following proposition\\

\textbf{Proposition 3.2.6} \emph{The inclusion map $i:\mathbb{C}^{n}\rightarrow{}{}^{b}\mathbb{C}^{n}$
induces an epimorphism of sheaves}\[
\mathcal{O}_{^{b}\mathbb{C}^{n}}\,\rightarrow\, i_{*}(\mathcal{O}_{\mathbb{C}^{n}})\]
 \emph{and so} \emph{we} \emph{obtain} $i$ \emph{as a morphism of
ringed spaces}.\\

\textbf{Proposition 3.2.7} \emph{Let $\mathfrak{a}$ be an ideal of
${}{}^{b}\mathbb{C}[X_{1},\ldots,X_{n}]$. Then}

\emph{\[
\mathrm{Supp}(\mathcal{O}_{^{b}\mathbb{C}^{n}}/\mathfrak{a\mathcal{O}_{^{b}\mathbb{C}^{n}}})=V(\mathfrak{a}).\]
 }

\emph{Proof}. It is clear that $V(\mathfrak{a})\subset\mathrm{Supp}(\mathcal{O}_{^{b}\mathbb{C}^{n}}/\mathfrak{a\mathcal{O}_{^{b}\mathbb{C}^{n}}})$.
Conversely, if $x\not\in V(\mathfrak{a})$ then there exists $f\in\mathfrak{a}$
such that $f(x){\in{}}^{a}\mathbb{C}$. Since $f$ is S-continuous,
we have $f(y){\in{}}^{a}\mathbb{C}$ for each $y\in\mu(x)$, which
implies that $1/f\in\mathcal{O}_{^{b}\mathbb{C}^{n},\, x}$ and as
a consequence we get $x\not\in\mathrm{Supp}(\mathcal{O}_{^{b}\mathbb{C}^{n}}/\mathfrak{a\mathcal{O}_{^{b}\mathbb{C}^{n}}})$,
q.e.d.\\

Now, we are able to define the category of affine {*}-bounded schemes.
An object is given by a locally ringed space $({}^{b}X,\mathcal{O}_{^{b}X})$.
The underlying topological space is defined by $^{b}X:=V(\mathfrak{a}){\subset{}}^{b}\mathbb{C}^{n}$
where $\mathfrak{a}$ is an ideal of \emph{${}{}^{b}\mathbb{C}[X_{1},\ldots,X_{n}]$}
and $\mathcal{O}_{^{b}X}:=j_{b}^{-1}(\mathcal{O}_{{}{}^{b}\mathbb{C}^{n}}/\mathfrak{a}\mathcal{O}_{{}{}^{b}\mathbb{C}^{n}})$.
The mapping $j_{b:}:{}^{b}X\rightarrow{}^{b}\mathbb{C}^{n}$ denotes
the inclusion.

Let $({}^{b}Y,\mathcal{O}_{^{b}Y})$ be another object in the category
of affine {*}-bounded schemes where $^{b}Y:=V(\mathfrak{b}){\subset{}}^{b}\mathbb{C}^{m}$.
A morphism between $({}^{b}X,\mathcal{O}_{^{b}X})$ and $({}^{b}Y,\mathcal{O}_{^{b}Y})$
is represented by a morphism of $^{b}\mathbb{C}$-algebras, \[
u_{b}:{}{}^{b}\mathbb{C}[Y_{1},\ldots,Y_{m}]\rightarrow{}^{b}\mathbb{C}[X_{1},\ldots,X_{n}]\]
 sending the ideal $\mathfrak{b}$ into $\mathfrak{a}$. We demand
that the morphism $u_{b}$ should lift to an internal morphism $u:{}\mathbb{C}[Y_{1},\ldots,Y_{m}]_{int}\rightarrow\mathbb{C}[X_{1},\ldots,X_{n}]_{int}$
defined by $u(Y_{i})=g_{i}\in{}{}^{b}\mathbb{C}[X_{1},\ldots,X_{n}]$
for each $i\in[\![1..m]\!]$. Then $u_{b}$ induces a morphism of
locally ringed spaces between $({}^{b}X,\mathcal{O}_{^{b}X})$ and
$({}^{b}Y,\mathcal{O}_{^{b}Y})$ as follows: the morphism $u_{b}$
gives a S-continuous map $u_{b}'{}:{}{}^{b}\mathbb{C}^{n}\rightarrow{}^{b}\mathbb{C}^{m}$
sending ${}{}^{b}X$ into ${}{}^{b}Y$, defined by ${}u_{b}'(x)=(g_{1}(x),\ldots,g_{m}(x))$.
We will denote by \[
\alpha_{b}:{}{}^{b}X\rightarrow{}{}^{b}Y\]
 the restriction of $u_{b}'$ on $^{b}X$. The mapping $u_{b}'$ defines
a morphism of sheaves $\mathcal{O}_{^{b}\mathbb{C}^{m}}\rightarrow(u_{b}')_{*}(\mathcal{O}_{^{b}\mathbb{C}^{n}})$.
Finally, we obtain a morphism of sheaves of rings\[
\mathcal{O}_{^{b}Y}\rightarrow{}\alpha_{b*}(\mathcal{O}_{^{b}X})\]
 induced by $\alpha_{b}$.\\

We state the main result of this section\\

\textbf{Theorem 3.2.8} \emph{There is a natural essentially surjective
functor \[
\mathrm{st}:(\mathrm{aff\,}*b-\mathrm{sch})\longrightarrow(\mathrm{Stein}\,\mathrm{spaces})\]
 from the category of affine {*}-bounded schemes to the category of
Stein spaces of finite embedding dimension.}

\emph{Moreover, if} $({}^{b}X,\mathcal{O}_{^{b}X})$ \emph{is} \emph{an
affine {*}-bounded scheme and} $(X,\mathcal{O}_{X}):=\mathrm{st}({}^{b}X,\mathcal{O}_{^{b}X})$\emph{,
then} \emph{the inclusion map $i_{b}:X{\rightarrow{}}^{b}X$ is a
monomorphism.}\\

\emph{Proof}. Let $({}^{b}X,\mathcal{O}_{^{b}X})$ be an object in
the category of affine {*}-bounded schemes, where $^{b}X:=V(\mathfrak{a}){\subset{}}^{b}\mathbb{C}^{n}$.
We set\[
\mathrm{st}({}^{b}X,\mathcal{O}_{^{b}X}):=(X,\mathcal{O}_{X})=(Z(°\mathfrak{a}),\, i^{-1}(\mathcal{O}_{\mathbb{C}^{n}}/\,°\mathfrak{a}\mathcal{O}_{\mathbb{C}^{n}}))\]
 where $i\,:\, Z(°\mathfrak{a})\rightarrow\mathbb{C}^{n}$ denotes
the inclusion.

By taking the standard part, the inclusion mapping $i_{b}:\, X\rightarrow{}{}^{b}X$
induces a surjective morphism of sheaves\[
\mathcal{O}_{^{b}X}\rightarrow{}i_{b*}(\mathcal{O}_{X}).\]

Let $({}^{b}Y,\mathcal{O}_{^{b}Y})$ be another object in the category
of affine {*}-bounded schemes where $^{b}Y:=V(\mathfrak{b}){\subset{}}^{b}\mathbb{C}^{m}$.
A morphism $u$ between $({}^{b}X,\mathcal{O}_{^{b}X})$ and $({}^{b}Y,\mathcal{O}_{^{b}Y})$
gives a morphism between the two associated Stein algebras \[
\mathrm{\mathrm{st}}(u):\mathcal{\, O}(\mathbb{C}^{m})/\,\overline{^{°}\mathfrak{b}}\longrightarrow\mathcal{O}(\mathbb{C}^{n})/\,\overline{^{°}\mathfrak{a}}.\]
 Such a morphism induces a morphism between the Stein spaces $(X,\mathcal{O}_{X})$
and $(Y,\mathcal{O}_{Y})$.

It remains to prove that the functor $"\mathrm{st}"$ is essentially
surjective. Let $X$ be a Stein space of finite embedding dimension.
Then there exists an ideal $\mathfrak{a\subset}\mathcal{O}(\mathbb{C}^{n})$,
generated by a family of entire holomorphic functions $(g_{i})_{i\in I}$
, such that $X=Z(\mathfrak{a})$ and $\mathcal{O}_{X}=\mathcal{O}_{\mathbb{C}^{n}}/\,\mathfrak{a}\mathcal{O}_{\mathbb{C}^{n}}$.
For each $i\in I$, let $h_{i}$ be a bounded internal polynomial
such that ${}{}^{\circ}h_{i}=g_{i}$. We put $\mathfrak{a}_{b}=\sum_{i\in I}h_{i}\,^{b}\mathbb{C}[X_{1},\ldots,X_{n}]$,
$^{b}X=V(\mathfrak{a}_{b})$ and $\mathcal{O}_{^{b}X}:=\mathcal{O}_{{}{}^{b}\mathbb{C}^{n}}/\mathfrak{a}_{b}\mathcal{O}_{{}{}^{b}\mathbb{C}^{n}}$.
Clearly, we obtain \emph{$\mathrm{st}({}^{b}X,\mathcal{O}_{^{b}X}):=(X,\mathcal{O}_{X})$},
q.e.d.\\

\textbf{Theorem 3.2.9} \emph{There is a natural functor between the
category of affine schemes of finite type over $\mathbb{C}$ and the
category of affine {*}-bounded schemes}.\\

\emph{Proof.} Each polynomial over $\mathbb{C}$ can be regarded as
an internal bounded polynomial. We consider the functor which associates
to each algebra of finite type over $\mathbb{C}$, say $A=\mathbb{C}[X_{1},\ldots,X_{n}]/(f_{1},\ldots,f_{q})$,
the affine {*}-bounded scheme $^{b}X:{={}}^{*}Z(f_{1},\ldots,f_{q}){\cap\,}^{b}\mathbb{C}^{n}$
and the structural sheaf $\mathcal{O}_{{}{}^{b}\mathbb{C}^{n}}/\sum_{i=1}^{q}f_{i}\mathcal{O}_{{}{}^{b}\mathbb{C}^{n}}$.
This construction is functorial and gives us our functor. -

~

\section{Global {*}-bounded schemes and the standard part functor}

\subsection{The basic functors}

4.1.1 In this section we want to define a category of global schemes
where the affine pieces are those described in section 3. This is
done in the usual way via locally ringed spaces and we obtain immediately
two new categories, the category of general {*}-schemes denoted by

~

$(*-\mathrm{sch})$

~\\
 and

~

$(*b-\mathrm{sch})$

~ \\
 that of so-called \emph{{*}-bounded} schemes which are more related
to complex geometry. The first one is of very general nature (i.e.
a category of schemes over $^{*}\mathbb{Z}$) whilst the second one
is a category of schemes over the ring $^{b}\mathbb{C}$ of bounded
complex numbers. We note that the structure sheaf of any {*}-bounded
scheme $X$ contains an intrinsic ideal sheaf $\mathcal{I}nf_{X}$
of infinitesimal sections which is locally described in section (3.2).

Our definition of a {*}-scheme is apriorily more general than that
of \cite{Brunjes2} where only finite coverings by affine pieces are
allowed. In order to define the notion of an internal subspace, we
would have to make restrictions on the cardinality of coverings, since
\emph{internal} is a global property.

~

We have an evident functor {}``associated {*}-scheme'' from the
category of algebraic $\mathbb{C}$-schemes, locally of finite type

\[
\mathrm{(\mathbb{C}-schemes\: l.f.t.)}\longrightarrow(*b-\mathrm{sch})\]
 which associates to a locally algebraic $\mathbb{C}$-scheme $X$
its {*}-bounded version $^{b}X$ and similar on the level of morphisms.

~

\noun{Example} 4.1.2 If $\mathbb{P}_{\mathbb{C}}^{n}$ denotes the
n-dimensional projective space over the complex numbers, then $^{b}\mathbb{P}_{\mathbb{C}}^{n}$
can be constructed as the quotient of $^{b}\mathbb{C}^{n+1}\setminus\mu(0)$
by the action of the multiplicative system of all appreciable complex
numbers. Note that the natural map $^{b}\mathbb{P}_{\mathbb{C}}^{n}\rightarrow{}^{*}\mathbb{P}_{\mathbb{C}}^{n}$
is bijective ($n$ is standard).

~

Recall that to every affine bounded scheme we associated in section
3 a Stein complex space (of finite embedding dimension) in a functorial
way. By a straightforward generalisation, we obtain a functor

\[
\mathrm{st}:(*b-\mathrm{sch})\longrightarrow(\mathrm{complex}\,\,\mathrm{spaces})\]
 which we call the \emph{standard part} \emph{functor}. We sometimes
write also $°X=\mathrm{st}(X)$ for the standard part of a bounded
scheme $X$. We may regard $°X$ as the {}``subspace'' of $X$,
defined by the ideal sheaf $\mathcal{I}nf_{X}$ all infinitesimal
sections and so we have a natural morphism of ringed spaces

~

$°X\rightarrow X$. 

~\\
 The conormal sheaf of this embedding is of particular interest, since
it gives us a non standard interpretation of classical one differential
forms on $°X$ (see section 6).

~

\textbf{Theorem 4.1.3} \emph{The standard part functor $"st"$ possesses
a left adjoint functor $^{b}:(\mathrm{complex\,\, spaces})\rightarrow(*\mathrm{b-sch})$
which associates to every complex space $X$ a natural {*}-bounded
scheme $^{b}X$ with locally no nontrivial infinitesimal elements.
Moreover, the adjunction morphism $id\rightarrow\mathrm{st\circ({}^{b})}$
is an isomorphism.}

~

\emph{Proof}. We first describe the functor $X\mapsto{}^{b}X$ locally,
i.e. for $X$ a finite dimensional Stein spaces (as in section 2).
As we have seen, there is an affine {*}-bounded scheme $\mathfrak{X}$
such that $X\cong{}^{st}\mathfrak{X}$. We can make now a natural
minimal choice of $\mathfrak{X}$ where there are \emph{no} nontrivial
infinitesimal elements in the local rings. In this case, the natural
homomorphism $\mathcal{O}_{\mathfrak{X}}\rightarrow(i_{b})_{*}\mathcal{O}_{X}$
will be an isomorphism. This construction is evidently functorial
in $X$, so it globalizes to complex spaces and we get our desired
functor. Intuitively speaking, we just enlarge $X$ (locally) by its
bounded points and conserve at the same time its structure sheaf.

In order to verify the adjunction property, let $\mathfrak{X}$ be
a {*}-bounded scheme and $Y$ a complex space. By applying the functor
{}``st'' , we get a functorial map\[
\mathrm{Hom}({}^{b}Y,\mathfrak{X})\longrightarrow\mathrm{Hom}(Y,\mathrm{st}(\mathfrak{X}))\]
 which is immediately seen to be bijective, since there are locally
no nontrivial infinitesimal elements in the structure sheaf of $^{b}Y$.
The last assertion is an obvious consequence of our construction,
q.e.d.

~

\subsection{DG- algebra resolutions of complex spaces via bounded polynomial
algebras}

Let us now consider a differential graded algebra $R=\oplus_{k\leq0}R_{k}$
with differential $s$ such that $R_{0}$ is an algebra of bounded
polynomials and $R$ is freely generated over $R_{0}$. Then we can
associate to it a Stein algebra in the following way: We may write
$R_{0}={}^{b}\mathbb{C}[T_{1},...,T_{n}]$ . Then $\mathcal{ST}\mathrm{(Coker}(R_{-1}\rightarrow R_{0}))$is,
by 2.2.3, a Stein algebra. This definition is clearly functorial.
Since it is evidently possible to construct DG-resolutions of that
type for a given Stein algebra of finite embedding dimension, we obtain

~

\textbf{Proposition 4.2.1} \emph{The} \emph{above} \emph{defined}
\emph{functor from the category of free bounded DG-algebras which
are exact in negative degrees, to the category of Stein algebras of
finite embedding dimension is essentially surjective.}

~

We generalize this fact to the simplicial case. Let $M$ be a totally
ordered set and $\mathcal{N}$ the category of non empty finite subsets
of $M$. An object $\alpha$ of $\mathcal{N}$ (i.e. a non empty finite
subset of $M$) is called a simplex. Its dimension is, by definition,
$\mathrm{card}(\alpha)-1$. The set $\mathrm{Mor}_{\mathcal{N}}(\alpha,\beta)$
is of cardinality $\leq1$ corresponding to the condition if $\alpha\subset\beta$
or not. A simplicial object in a category $\mathcal{C}$ is a \emph{contravariant}
functor $\mathcal{N}\rightarrow\mathcal{C}$. The category of all
these functors (or $\mathcal{N}$-objects) is denoted by $\mathcal{C}_{\mathcal{N}}$.
We can show

~

\textbf{Proposition 4.2.2} \emph{The} \emph{above defined functor
{}``standard part in degree zero cohomology'' from the category
of free bounded DG-$\mathcal{N}$-algebras which are exact in negative
degrees, to the category of Stein $\mathcal{N}$-algebras of simplicialwise
finite embedding dimension, is} \emph{essentially} \emph{surjective}.

~

For a sketch of \emph{proof}, given a simplicial Stein algebra $A$
(simplicialwise of finite embedding dimension), we proceed in the
usual way by induction on the dimension of simplicies to construct
a free DG-algebra resolution $R_{\alpha}$ such that $°\mathrm{H}^{0}(R_{\alpha})=A_{\alpha}$.
If $\mathrm{dim}\alpha=0$, this is has been done. For $\mathrm{dim}\alpha>0$
we know already the values of the differential on all free generators
which come from strictly lower dimensional simplicies, so that we
just add some new generators and construct the differential $s$ by
descending induction on the degree of generators in order to obtain
exactness in negative degrees and such that $\mathcal{ST}(\mathrm{H}^{0}(R_{\alpha}))=A_{\alpha}$
. -

~

4.2.3 We want to describe another more subtle construction of such
resolutions, by imposing additional assumptions on the resolution
type. For this, we fix a complex space $X$, an embedding $X_{\cdot}\hookrightarrow P_{\cdot}$
into a simplicial (free) polydisc $P_{\cdot}$ and a free simplicial
DG-algebra resolution $(\mathcal{R},s)$ of $\mathcal{O}_{X_{\cdot}}$
of the form $\mathcal{R}=\mathcal{O}_{P_{\cdot}}[e_{i}]_{i\in I}$
with free generators $e_{i}$ of strictly negative degrees. The graded
simplicial algebra $\mathcal{R}$ admits an evident lifting to the
category of simplicial affine bounded algebras. But we cannot extend
directly the identity $s^{2}=0$ to this algebra. Neverless, we will
obtain it on a suitable (even maximal) affine subspace $\mathfrak{Z}_{\cdot}\hookrightarrow{}^{b}P_{\cdot}$,
after fixing an infinite natural number $N$. The definition of $\mathfrak{Z}_{\cdot}$
is the following: we first lift $s$ to bounded derivation $s_{N}$
on $^{b}\mathcal{R}=\mathcal{O}_{^{b}P_{\cdot}}[e_{i}]_{i\in I}$.
Then $s_{N}^{2}(e_{i})$ is infinitesimal by construction. Dividing
out these infinitesimal polynomials (over each simplex), we obtain
a subspace $\mathfrak{Z}_{\cdot}\hookrightarrow{}^{b}P_{\cdot}$ over
which the class of $s_{N}^{2}$ is zero. In other words, we get a
simplicial DG-algebra in the category of bounded algebras, conserving
all original free generators of negative degree (and not adding any
further ones).

~

\noun{Remark} 4.2.4 The construction in 4.2.3 is inspired by deformation
theory of complex spaces: describing all small deformations of $X$
by varying the differential of a fixed DG-algebra resolution of $\mathcal{O}_{X}$.
Moreover, our construction leads us to an obvious definition of what
should be a deformation of $X$ over local {*}-bounded algebras, for
example algebras of the form $^{b}\mathbb{C}/\mathfrak{a}$ where
$\mathfrak{a}$ is an ideal generated by some infinitesimal complex
numbers. Such \emph{nonstandard} or \emph{Leibniz} \emph{deformations}
give an interesting alternative approach to deformation theory where
the meaning of {}``infinitesimal deformation'' becomes a metric
one.

\subsection{Special features of the standard part construction}

We come to some remarkable properties of our standard part construction.
Let $H$ be a standard hyperplane in the {*}-scheme ${}^{*}\mathbb{P}_{\mathbb{C}}^{n}$.
Then the natural inclusion ${}^{*}\mathbb{C}^{n}\hookrightarrow{}^{*}\mathbb{P}_{\mathbb{C}}^{n}$
induced by $H$, allows us to define the notion of a \emph{bounded}
\emph{point} \emph{in} ${}^{*}\mathbb{P}_{\mathbb{C}}^{n}\setminus H$.
The following theorem shows that we can {}``compactify'' analytic
subsets of $\mathbb{C}^{n}$ in a (non-canonical) way to projective
{*}-schemes

~

\textbf{Theorem 4.3.1} \emph{Let $X\subset\mathbb{C}^{n}$ be an analytic
subset, defined by a finite number of entire functions on $\mathbb{C}^{n}$
. Then there is a standard hyperplane $H\subset{}^{*}\mathbb{P}_{\mathbb{C}}^{n}$
and a {*}-projective subvariety $\mathfrak{X}\subset{}^{*}\mathbb{P}_{\mathbb{C}}^{n}$
such that $X$ is the standard part of $\mathfrak{X}\setminus H$.}

~

\emph{Proof}. Let $X$ be the zero set of a finite family $(f_{j})_{j\in J}$
of holomorphic functions on $\mathbb{C}^{n}$. We fix an infinite
natural number $N$ and extensions $F_{j}$ of each $f_{j}$ to a
finite family of {*}-bounded polynomials of degree $\leq N$. Next,
we homogenize each $F_{j}$ to a homogeneous internal (but in general
not {*}-bounded !) polynomial $\tilde{F}_{j}$ of degree $\leq N$
in $n+1$ variables. Let $\mathfrak{X}\subset{}^{*}\mathbb{P}_{\mathbb{C}}^{n}$
be defined by the $\tilde{F}_{j}$'s. Then, by our construction, $\mathfrak{X}$
has the desired property, q.e.d.

~

\textbf{Theorem 4.3.2} \emph{Let $X$ be an affine algebraic $\mathbb{C}$-scheme
such that $\mathrm{H}^{2}(X^{an},\mathbb{Z})=0$. Then every holomorphic
map $f:X^{an}\rightarrow\mathbb{P}_{\mathbb{C}}^{n}$ from $X^{an}$
to a projective space is the standard part of a {*}-bounded morphism
$F:{}^{b}X\rightarrow^{b}\mathbb{P}_{\mathbb{C}}^{n}$.}

~

\emph{Proof}. By our assumption, we can lift $f$ to a holomorphic
map $\tilde{f}:X^{an}\rightarrow\mathbb{C}^{n+1}\setminus\{0\}$.
By the results of section 2.5, the map $\tilde{f}$ is the standard
part of a {*}-bounded map $F:{}^{b}X\rightarrow{}^{b}\mathbb{C}^{n+1}$
which avoids necessarily the halo of $\{0\}.$ By passing to the quotient
(as in 4.1.2), we get our result, q.e.d.

~

\noun{Example} 4.3.3. Let $f:\mathbb{C}\rightarrow\mathbb{P}_{\mathbb{C}}^{n}$
be a holomorphic map. By fixing an infinite natural number $N$, we
may extend $f$ to a morphism $f_{N}:^{b}\mathbb{C}\rightarrow^{b}\mathbb{P}_{\mathbb{C}}^{n}$
of {*}-bounded schemes such that $\mathrm{st}(f_{N})=f$. Moreover,
applying homogenization, we can lift $f_{N}$ to a {*}-algebraic map
$\tilde{f_{N}}:{}^{*}\mathbb{P}^{1}\rightarrow{}^{*}\mathbb{P}_{\mathbb{C}}^{n}$
which will be a {*}-\emph{bounded} morphism if and only if $f$ is
algebraic. In fact, we can define $\tilde{f_{N}}:{}^{b}\mathbb{P}^{1}\rightarrow{}^{b}\mathbb{P}_{\mathbb{C}}^{n}$
just as a map, but it will be not continuous at the point $(1:0)$
if $f$ is not algebraic.

~\\
 We close this section by mentioning a simple fact concerning the
standard points of internal subspaces of nonstandard complex spaces

~

\textbf{Theorem 4.3.4} \emph{Let $X$ be a complex space and $\mathfrak{Y}\subset{}^{*}X$
an internal complex subspace of (finite) codimension $k$. Then the
intersection $\mathfrak{Y}\cap X$ is contained in a countable union
of complex subspaces $Y_{n}$, $n\in\mathbb{N}$, of $X$ of codimension
$k$. The same conclusion holds in the algebraic context.}

\section{The Nullstellensatz and nonstandard generic points for complex spaces}

In this section we prove in particular that any maximal ideal of a
Stein algebra $A$ is the vanishing ideal of an eventually nonstandard
point. Moreover, we can show that any prime ideal of $A$ is determined
by its nonstandard zero set if and only if it satisfies a Hilbert
Nullstellensatz. We give a large class of examples of prime ideals
which are not determined by their nonstandard zero set and with necessarily
empty standard zero set. Closed ideals in Stein algebras are extensively
treated in Forster's paper \cite{Forster}. 

Our notion of a nonstandard zero set of an ideal allows us to interpret
generic points of irreducible complex spaces in a natural geometric
way.

\subsection{Nonstandard zeros of holomorphic functions and the Nullstellensatz}

5.1.1 Let $X$ be a Stein complex space and $A:=\Gamma(X,\mathcal{O}_{X})$
its algebra of global holomorphic sections. For any ideal $I\subset A$,
we define the \emph{nonstandard} zero set of $I$ by

\[
\mathrm{V}(I):=\{ x\in{}^{*}X\,\mid\,{}^{*}f(x)=0,\,\forall f\in I\}.\]

Forster proved that the (standard) zero set of a proper closed ideal
in a Stein algebra is non empty. Siu\cite{Siu} showed that the closedness
is essential and gave an example of an ideal whose its variety is
empty. We prove that the \emph{nonstandard} zero set of any proper
ideal is nonempty. This indicates that our definition is an adequate
notion of a zero set in the context of Stein algebras. Let $I\subsetneq A$
be a proper ideal. We define the binary relation $P$ on $X\times I$
by $P<x,f>$ if $f(x)=0$. The relation $P$ is concurrent: if $f_{1},\ldots,f_{r}\in I$
then there exists $x\in X$ such that $f_{i}(x)=0$ for each $i=1,\ldots,r$.
Since $^{*}X$ is an enlargement of $X$, there exists $x\in\,^{*}X$
such that $^{*}f(x)=0$ for all $f\in I$ and hence $\mathrm{V}(I)\not=\emptyset$
(see appendix B).

~

There is also the dual construction: Let $M\subset{}^{*}X$ be any
subset. We put\[
\mathrm{Id}(M):=\{ f\in A\,\mid\,{}^{*}f(x)=0,\,\forall x\in M\}\]
 the \emph{nonstandard} ideal of $M$ in $A$. These constructions
transform an inclusion into the opposite one. Moreover, we have the
following rules: Let $(I_{r}){}_{r\in R}$ be a family of ideals in
$A$ and $(M_{s}){}_{s\in S}$ a family of subsets of $^{*}X$. Then

\[
\mathrm{V}(\Sigma_{r\in R}I_{r})=\cap_{r\in R}\mathrm{V}(I_{r}),\]
 and

\[
\mathrm{V}(\cap{}_{r\in R}I_{r})=\mathrm{V}(\Pi_{r\in R}I_{r})=\cup{}_{r\in R}\mathrm{V}(I_{r})\]
 if $R$ is finite, dually, for {}``$\mathrm{Id}$''\[
\mathrm{Id}(\cup{}_{s\in S}M_{s})=\cap_{s\in S}\mathrm{Id}(M_{s}),\]
 and \[
\mathrm{Id}(\cap{}_{s\in S}M_{s})\supset\Sigma_{s\in S}\mathrm{Id}(M_{s})\]
 Moreover, for any subset $M\subset{}^{*}X$ and any ideal $I\subset A$,
we have clearly

\[
\mathrm{V}(\mathrm{Id}(M))\supset M,\]

\[
\mathrm{Id}(\mathrm{V}(I))\supset\sqrt{I}.\]
 More precisely, for the first one, we show

~

\textbf{Proposition 5.1.2} \emph{We} \emph{have} $\mathrm{V}(\mathrm{Id}(M))=\cap_{Y}\mathrm{V}(\mathrm{Id}(Y))$,
\emph{where} \emph{the} \emph{intersection} \emph{is} \emph{taken}
\emph{over} \emph{all} $Y\subset X$ \emph{closed analytic} \emph{subsets
such} \emph{that} $M\subset{}^{*}Y$.

~

\emph{Proof}. If $Y\subset X$ is analytic, we get $\mathrm{Id}({}^{*}Y)=\mathrm{Id}(Y)$
and therefore $\mathrm{V}(\mathrm{Id}({}^{*}Y))=\mathrm{V}(\mathrm{Id}(Y))\supset{}^{*}Y$.
We note that $\mathrm{Id}(M)=\Sigma_{Y}\mathrm{Id}(Y)$, where the
sum is taken over all $Y\subset X$ analytic such that $M\subset{}^{*}Y$,
which gives us the desired identity. -

~

We want to show that every prime ideal of $A$ is the ideal of a (nonstandard)
point of $X$ if and only if it satisfies a Hilbert Nullstellensatz. 

~

\textbf{Theorem 5.1.3} i) \emph{Let $\mathfrak{m}$ be a maximal or
a minimal prime ideal of $A$, then $\mathfrak{m}$ is of the form}
$\mathrm{Id}(\{ x\})$,

ii) \emph{Let $\mathfrak{p}$ be a prime ideal of $A$, then $\mathfrak{p}$
is of the form} $\mathrm{Id}(\{ x\})$, \emph{if} \emph{and} \emph{only}
\emph{if it} \emph{satisfies} \emph{the} \emph{Nullstellensatz, i.e.}
$\mathrm{Id}(\mathrm{V(\mathfrak{p})})=\mathfrak{p}$.

~

\emph{Proof}. i) We first treat the case of a maximal ideal. By the
Nullstellensatz, we have\[
\forall n\in\mathbb{N},\forall f_{1},...,f_{n}\in\mathfrak{m}\Rightarrow\exists\, y\in X:f_{i}(y)=0\,\forall i,1\leq i\leq n.\]
 By the concurrence principle (see the appendix), we conclude that
there is a point $x\in{}^{*}X$ such that $\mathfrak{m}\subseteq\mathrm{Id(\{ x\})}$.
Since $\mathfrak{m}$ is maximal, we have equality.

ii) First, let $\mathfrak{p}$ be of the form $\mathrm{Id}(\{ x\})$.
Then $\mathrm{V}(\mathfrak{p})\supset\{ x\}$and so $\mathrm{Id}(\mathrm{V}(\mathfrak{p}))\subset\mathrm{Id}(\{ x\})=\mathfrak{p}$,
which finally gives equality. For the converse, we note that it suffices
to show

\[
\forall n\in\mathbb{N},\forall f_{1},...,f_{n}\in\mathfrak{p},\forall g_{1},...,g_{n}\notin\mathfrak{p}\Rightarrow\exists y\in X:f_{i}(y)=0\,\forall i,\, g_{1}(y)\cdot...\cdot g_{n}(y)\neq0\]
 and since $g_{1}\cdot...\cdot g_{n}\notin\mathfrak{p}$, we may take
$n=1$ for the $g_{i}'s$ for the verification of this implication.
By applying the concurrence principle, we obtain\[
\exists\, x\in{}^{*}X:\:\forall f\in\mathfrak{p},\forall g\notin\mathfrak{p}\Rightarrow{}^{*}f(x)=0,\,{}^{*}g(x)\neq0\]
 which means precisely that $\mathfrak{p}=\mathrm{Id}_{X}(x)$. So
let now be $f_{1},...,f_{n}\in\mathfrak{p}$ and $g\notin\mathfrak{p}$.
Assume that for the usual zero sets we have $\mathrm{Z(}f_{1},...,f_{n})\subset\mathrm{Z}(g)$.
Then we get also $\mathrm{V(}f_{1},...,f_{n})\subset\mathrm{V}(g)$
and so for the corresponding ideals $\mathrm{Id}(\mathrm{V(}f_{1},...,f_{n}))\supset\mathrm{Id}(\mathrm{V}(g))\ni g$.
But, $\mathrm{Id}(\mathrm{V(}f_{1},...,f_{n}))\subset\mathrm{Id}(\mathrm{V}(\mathfrak{p}))=\mathfrak{p}$
and we would obtain $g\in\mathfrak{p}$, which is a contradiction.
This shows (ii). 

iii) It remains to show that any minimal prime ideal is that of a
point. Let $\mathfrak{p}$ be a prime ideal. In contrast to the formula
above, the following implication is always true

~\[
\forall n\in\mathbb{N},\forall g_{1},...,g_{n}\notin\mathfrak{p}\Rightarrow\exists y\in X:\, g_{1}(y)\cdot...\cdot g_{n}(y)\neq0\]
 since it is true for $n=1$ and $\mathfrak{p}$ is prime. By the
concurrence principle, we conclude that there is a point $x\in{}^{*}X$
such that\[
\forall g\notin\mathfrak{p}\Rightarrow g(x)\neq0\]
 which means $\mathfrak{p}\supset\mathrm{Id}_{X}(x)$, so finally
$\mathfrak{p}=\mathrm{Id}_{X}(x)$by minimality of $\mathfrak{p}$,
q.e.d.

~

\noun{Examples} 5.1.4 i) Regarding the inclusion $A\hookrightarrow{}^{*}A$,
we obtain certain prime ideals of $A$ by intersecting just with those
of $^{*}A$. In particular every internal prime ideal of $^{*}A$
gives us one of $A$ (for example, fixing $(\mathfrak{p}_{i})_{i}$,
a countable sequence of closed prime ideals of $A$).

ii) Assume $A$ to be an integral domain and that we have an {}``order
function'' $\mathbb{\omega}:A\setminus\{0\}\rightarrow{}^{*}\mathbb{N}$,
i.e. $\omega$ satisfies $\omega(fg)=\omega(f)+\omega(g)$ and $\omega(f+g)\geq\mathrm{inf}\{\omega(f),\omega(g)\}$
as well as $\omega(1)=0$. Then the subset of $A$ formed by zero
and all elements of \emph{infinite} $\omega$-order is a prime ideal
of $A$.

iii) Combining the constructions of (i) and (ii), we can immediately
construct explicitly many prime ideals of the ring of entire holomorphic
functions on $\mathbb{C}^{n}$ which are \emph{not} ideals of nonstandard
points.

~

We can determine the possible \emph{residue} \emph{fields} of a maximal
ideal in a Stein algebra 

~

\textbf{Theorem 5.1.5} \emph{Let $A$ be a Stein algebra (of finite
embedding dimension) and $\mathfrak{m}=\mathrm{Id}(x)$ any maximal
ideal of $A$. We suppose that the underlying ultrafilter is $\delta-$stable
(see appendix A).}

\emph{Then the residue field of $\mathfrak{m}$ is $^{*}\mathbb{C}$
if $\{ x\}$ is an infinite point and $\mathbb{C}$ if $\{ x\}$ is
bounded. In the last case, $\{ x\}$ is standard and so $\mathfrak{m}$
is a closed ideal.}

~

\emph{Proof}. We may assume that $A=\Gamma(\mathbb{C}^{n},\mathcal{O}_{\mathbb{C}^{n}})$.
If $\{ x\}$ is an infinite point, then at least one coordinate of
$\{ x\}$, say $\{ x_{k}\}$, is infinite. By projecting to this coordinate,
we may take first $n=1$. But in this case, the point can be represented
(using $\delta$-stability) by a sequence in $\mathbb{C}$ which tends
to $\infty$. Clearly, here the residue field must be $^{*}\mathbb{C}$,
by classical function theory in one complex variable. In this way,
we get homomorphisms of fields $\kappa(x_{k})\rightarrow\kappa(x)\rightarrow{}^{*}\mathbb{C}$,
where the second arrow is induced by evaluation in $\{ x\}$. The
composition is, by construction, evaluation in $^{*}\mathbb{C}$ and
therefore bijective.

If the point $\{ x\}$ is bounded in $^{*}\mathbb{C}^{n}$, it has
a standard part $a$. Since every holomorphic function which vanishes
in $\{ x\}$, must also vanish in $a$, we obtain, by maximality of
$\mathfrak{m}$, our result, q.e.d

~

For the rest of this section we will always assume that our ultrafilter
is $\delta$- stable.

\subsection{On the spectrum of a Stein algebra}

Here is the main result of this section

~

\textbf{Theorem 5.2.1} \emph{Let $X$ be a Stein complex space. Then
the image of the natural map}

\[
\mathrm{Id_{X}}:{}^{*}X\longrightarrow\mathrm{Spec}(\Gamma(X,\mathcal{O}_{X}))\]
 \emph{which} \emph{associates} \emph{to} \emph{a point} $x\in{}^{*}X$
\emph{its ideal} $\mathrm{Id}(\{ x\})$, \emph{consists of all prime
ideals satisfying the Nullstellensatz. Moreover, this image contains
all maximal, minimal and all closed prime ideals of} $\Gamma(X,\mathcal{O}_{X})$.

~

\emph{Proof}. First we remark that this map is well-defined, by considering
the evaluation homomorphism $\Gamma(X,\mathcal{O}_{X})\rightarrow{}^{*}\mathbb{C}$,
given by a point $x\in{}^{*}X$. Let $\mathfrak{p}$ be a prime ideal
in $\Gamma(X,\mathcal{O}_{X})$. By the results of the last subsection,
we only need to treat the case of closed prime ideals. But such a
prime ideal satisfies the usual Nullstellensatz (see \cite{Forster})
and so in particular the nonstandard one, i.e. $\mathrm{Id}(\mathrm{V}(\mathfrak{p}))=\mathfrak{p}$.
Again, 5.1.3 gives us our statement, q.e.d.

~

\textbf{Corollary 5.2.2} \emph{Let $X$ be an irreducible Stein complex
space. Then there exists a point} $x\in{}^{*}X$ \emph{such} \emph{that}
\emph{the} \emph{evaluation} \emph{in} $x$

\[
\chi_{x}:\Gamma(X,\mathcal{O}_{X})\longrightarrow{}^{*}\mathbb{C}\]
 \emph{is injective. For the field of meromorphic functions, we get}\[
\mathcal{M}(X)\longrightarrow{}^{*}\mathbb{C},\]
 \emph{the induced homomorphism of fields}.

~

Now we come to study some topological properties of the map $\mathrm{Id_{X}}$.
One can define a topology $\mathcal{T}$ on $^{*}X$ where $V(I)$
are closed sets for this topology (see section 5.1). For each $f\in A=\Gamma(X,\mathcal{O}_{X})$,
the open set \[
\mathcal{D}(f)=\{ x\in\:^{*}X,\:^{*}f(x)\not=0\}=\,^{*}\{ x\in X,\, f(x)\not=0\}=\,^{*}D(f)\]
 is a distinguished open set, the family of distinguished open sets
is a basic for the topology $\mathcal{T}$. Hence $\mathcal{T}$ is
the S-topology on $^{*}X$ when $X$ is equipped with the Zariski
topology where $\{ D(f),\, f\in A\}$ forms an open base, called the
\emph{S-Zariski topology on $^{*}X$} (see appendix C).

~

\textbf{Proposition} \textbf{5.2.3} Let \emph{Let $X$ be a Stein
complex space. Then the natural map}

\[
\mathrm{Id_{X}}:{}^{*}X\longrightarrow\mathrm{Spec}(\Gamma(X,\mathcal{O}_{X}))\]
\emph{is continuous for the S-Zariski topology on $^{*}X$ and the
Zariski topology on $\mathrm{Spec}(\Gamma(X,\mathcal{O}_{X}))$.}

~

\emph{Proof}. For each $x\in\,^{*}X$ and $f\in\,\mathrm{Spec}(\Gamma(X,\mathcal{O}_{X}))$,
we have $f\in\,\mathrm{Id_{X}}(x)$ if and only if, $^{*}f(x)=0$.
Hence, we get $\mathrm{Id_{X}}^{-1}(D(f))=\mathcal{D}(f)$ which shows
that our map is continuous. -

~

\textbf{Definition 5.2.4} Let $X$ be a complex space and $x\in{}^{*}X$
a point. The \emph{Zariski} \emph{closure} $\overline{x}$ of $x$
in $X$ is the smallest analytic subset $Y$ of $X$ such that $x\in{}^{*}V$
for every Zariski-open $V\subset Y$. If $Y\subset X$ is an analytic
subset and $x\in{}^{*}Y$, we call $x$ a \emph{generic} \emph{point}
\emph{of} $Y$ if $\overline{x}=Y$.

~

\noun{Remark} 5.2.5 i) The Zariski closure $\overline{x}$ of $x$
in $X$ is always irreducible.

ii) Note that $\overline{x}$ may be sometimes the empty set, for
example if $X=\mathbb{C}$ and $x\in{}^{*}\mathbb{C}$ is given be
a discrete sequence converging to infinity.

~

\textbf{Theorem 5.2.6} \emph{Let $X$ be an irreducible complex space.
Then} $X$ \emph{has always a generic point. Moreover, for every standard
point $x\in X$, there is a generic point of $X$ in the halo of $x$
(for the usual topology).}

~

\emph{Proof}. Let $U_{1},...,U_{n}$ be a finite family of non-empty
Zariski open subsets of $X$. Then, by the irreducibility of $X$,
the intersection $U_{1},\cap...\cap U_{n}$ is also non-empty. The
concurrence principle allows us to conclude that there is a point
$\eta\in{}^{*}X$ such that $\eta\in{}^{*}U$ for every non-empty
Zariski open subset $U$ of $X$. This point will be generic for $X$:
Let $Y\subset X$ be any strictly smaller analytic subset such that
$\eta\in{}^{*}Y$. Since we have $\eta\in{}^{*}(X\setminus Y)$ too,
we get a contradiction. This shows the first part of (5.2.6).

We refine slightly our argument in order to obtain the second part.
Let $x$ be a standard point of $X$. We decompose $X$ \emph{locally}
around $x\in X$ into irreducible components $X_{1}\cup...\cup X_{r}$
and fix one of them, say $X_{1}.$ In $X_{1}$, we consider the subsets
of the form $V\setminus Y$ where $V$ is a (usual) open neighborhood
of $x$ in $X_{1}$ and $Y\subset V$ a strictly smaller closed analytic
subset, defined by finitely many holomorphic functions on $V$. Clearly,
every finite intersection of such subsets is non-empty. By the concurrence
principle, there is a point $\eta\in{}^{*}(V\setminus Y)$ for every
one of our subsets $V\setminus Y$. Evidently, $\eta$ is bounded
and its standard part is $x$ by construction. We claim that $\eta$
is a generic point for $X$. Let $Z\subset X$ be a strictly smaller
analytic subset. By irreducibility and purity of dimension of $X$,
we cannot have $Z=X_{1}$ locally around $x$. So $Z$ is strictly
smaller than $X_{1}$ in some neighborhood of $x$ too. We therefore
get $\eta\in{}^{*}(X\setminus Z)$ which means that $\eta$ is a generic
point for $X$, q.e.d.

~

~

\noun{Remark} 5.2.7 We should mention in this context the paper
of G.Wallet \cite{Wallet} in which he shows algebraic versions of
the existence of nonstandard generic points by a more direct method.
The case of complex space space germs has already been treated by
A.Robinson in \cite{Robinsion2}.

\section{Differential forms seen in a modern nonstandard way}

In this section, we show how differential forms find a natural description
in the context of our {*}-bounded algebras. We prove in particular
that they coincide (in all relevant cases) with the associated analytic
ones.

~\\
Let $D_{n}$ be the subring of the internal polynomials in $2n$-
variables $\mathbb{C}[X,dX]_{int}=\mathbb{C}[X_{1},\ldots,X_{n},dX_{1},\ldots,dX_{n}]_{int}$
defined by

\[
D_{n}:=\{ P\in\mathbb{C}[X,dX]_{int}\,|\, P(x,\xi)\in{}{}^{b}\mathbb{C}\,\mbox{ for each }(x,\xi)\in{}{}^{b}\mathbb{C}^{n}\times{}^{i}\mathbb{C}^{n}\}.\]
\\
By theorem 2.2.1, it is immediate that the ring $D_{n}$ is invariant
under derivations. Furthermore, each $P\in D_{n}$ can be written
as\[
P(X,dX)=P_{0}(X)+\sum_{i=1}^{n}Q_{i}(X,dX)dX_{i}\]
 where $P_{0}(X)=P(X,0)\in{}{}^{b}\mathbb{C}[X_{1},\ldots,X_{n}]$
and $Q_{i}\in D_{n}$, since $Q_{i}(x,\xi)=\int_{0}^{1}\frac{\partial P}{\partial\xi_{i}}(x,t\xi)\, dt$.

More generally, let $A={}^{b}\mathbb{C}[X]/\mathfrak{a}$ be a {*}-bounded
algebra, using the abbreviation $X=(X_{1},...,X_{n})$. We define
an infinitesimal version of the diagonal algebra of $A$ by setting

\[
D(A):=D_{n}/\mathfrak{a}D_{n}+\delta(\mathfrak{a})D_{n}\]
 where $\delta(f):=f(X+dX)-f(X)$ for any $f\in{}^{b}\mathbb{C}[X]$.
Evidently, there is a natural injective homomorphism of rings $A\rightarrow D(A)$
and $D(..)$ is a covariant functor which conserves epimorphisms.
The standard part construction gives us an epimorphism

\[
\mathrm{st}:D(A)\longrightarrow A^{an}\]
 to the Stein algebra $A^{an}$, associated to $A$. Let $I_{A}$
be the kernel of this map. If $A={}^{b}\mathbb{C}[X]$, we simply
write $I_{n}$. Clearly, $I_{A}$ contains each $D(A)dX_{i}$ and
also $Inf(A)D(A)$ where $Inf(A)$ denotes, by definition, the kernel
of the standard part map $A\rightarrow A^{an}$. 

~

\textbf{Definition 6.1} We call the ideal $\mathfrak{a}$ \emph{saturated}
if the standard part map $\mathfrak{a}$$\rightarrow\overline{^{°}\mathfrak{a}}$
is surjective. 

~

\noun{Remark} 6.2 i) $\mathfrak{a}$ is saturated if $^{°}\mathfrak{a}$
is a closed ideal in $\mathcal{O}(\mathbb{C}^{n})$, for example if
$\mathfrak{a}$ is finitely generated.

ii) If $\mathfrak{a}$ is a closed subset of ${}^{b}\mathbb{C}[X]$
with respect to the S-topology (inherited and constructed from the
compact-open topology of $\mathcal{C}(\mathbb{C}^{n},\mathbb{C})$),
then $\mathfrak{a}$ is saturated.

iii) If $\mathfrak{a}$ is saturated, then the natural mappings $^{i}\mathbb{C}[X]\rightarrow Inf(A)$
and $I_{n}\rightarrow I_{A}$ are surjective . 

~

We can show now\\

\textbf{Proposition 6.3} \emph{i) $I_{A}=Inf(A)+\sum_{i=1}^{n}D(A)\, dX{}_{i}$
if $\mathfrak{a}$ is saturated,}

\emph{ii) for $f\in{}{}^{b}\mathbb{C}[X]$ we have, $f(X+dX)-f(X)\in I_{n}$,}

\emph{iii) for $f\in{}{}^{i}\mathbb{C}[X]$ we have, $f(X+dX)-f(X)\in I_{n}^{2}$.}\\

\emph{Proof.} The assertion (i) is clear for $A={}^{b}\mathbb{C}[X]$,
since $P\in I_{n}$ if and only if $P_{0}\in\,{}^{i}\mathbb{C}[X]$.
By saturation of $\mathfrak{a}$, we can reduce to the case $A={}^{b}\mathbb{C}[X]$.
For the second one, we already know that $f$ is S-continuous at each
$x\in{}{}^{b}\mathbb{C}^{n}$. Hence, for every $x\in{}{}^{b}\mathbb{C}^{n}$
and $\xi\in{}{}^{i}\mathbb{C}^{n}$, we have $f(x+\xi)-f(x)$$\approx0$.
The third assertion is a consequence of the stability of $\,{}^{i}\mathbb{C}[X]$
under partial derivations. - \\

\textbf{Lemma 6.4} \emph{If $\mathfrak{a}$ is saturated, then we}
\emph{have} $Inf(A)^{2}=Inf(A)$.\\

\emph{Proof.} By saturation, it is sufficient to treat the case $A={}^{b}\mathbb{C}[X]$.
Let $P=\sum a_{\nu}X^{\nu}\in{}{}^{i}\mathbb{C}[X]$. Then $a_{0}{\in\,}^{i}\mathbb{C}$
and $|a_{\nu}|^{\frac{1}{|\nu|}}{\in\,}^{i}\mathbb{C}$, for each
$0<|\nu|\leq d$ which implies in particular that $a_{\nu}{\in{}}^{i}\mathbb{C}$.
Setting $\varepsilon=\max_{0\leq|\nu|\leq d}\,|a_{\nu}|^{\frac{1}{2}}$,
we have $\varepsilon{\in{}\,}^{i}\mathbb{C}$, Define now\[
Q=\sum_{0\leq|\nu|\leq d}b_{\nu}X^{\nu}\quad\mbox{where\,\,\,}b_{\nu}:=\frac{a\nu}{\varepsilon}.\]
 It is easy to prove that $b_{0}{\in\,}^{i}\mathbb{C}$ and $|b_{\nu}|^{\frac{1}{|\nu|}}{\in\,}^{i}\mathbb{C}$,
hence $Q\in{}{}^{i}\mathbb{C}[X]$ and $P=\varepsilon Q\in({}^{i}\mathbb{C}[X])^{2}$.
- \\

Next, we consider the canonical map $\delta_{A}:A\rightarrow D(A)$
which associates to each class $[F]$ modulo $\mathfrak{a}$ of a
bounded polynomial $F\in{}^{b}\mathbb{C}[X]$, the class $[\delta(F)]$
in $D(A)$.\\

\textbf{Proposition 6.5} \emph{The induced mapping $\delta_{A}:A\rightarrow I_{A}/I_{A}^{2}$
is a derivation.}

~

\emph{Proof.} First, we note that we can reduce immediately to the
case $A={}^{b}\mathbb{C}[X]$. But here, it is the usual standard
verification. -

~

We remark that $I_{A}/I_{A}^{2}$ carries a natural $A^{an}$-module
structure (since annihilated by $Inf(A)\subset I_{A}$). So, we make
the following definition

~

\textbf{Definition 6.6} The {*}-\emph{bounded} \emph{module} \emph{of}
\emph{1-forms} is denoted by

\[
^{b}\Omega_{A}:=I_{A}/I_{A}^{2}.\]
 If $A={}^{b}\mathbb{C}[X]$, we simply write $^{b}\Omega_{n}$. 

~\\
By a standard verification, we get

~

\textbf{Proposition 6.7} \emph{If $\mathfrak{a}$ is saturated, the
natural sequence of $A$-modules}

\[
\begin{array}{ccccccc}
 & \delta\\
\mathfrak{a}/\mathfrak{a}^{2} & \longrightarrow & ^{b}\Omega_{n}/\mathfrak{a}\cdot{}^{b}\Omega_{n} & \longrightarrow & ^{b}\Omega_{A} & \rightarrow & 0\end{array}\]
 \emph{is exact} \emph{and all three terms are} \emph{$A^{an}$-modules}.

~\\
We are able to show the comparison theorem

~

\textbf{Theorem 6.8} \emph{The natural homomorphism} $\gamma_{A}:{}^{b}\Omega_{A}\rightarrow\Omega_{A^{an}}$
\emph{is always surjective and bijective if the ideal $\mathfrak{a}$
is saturated.}

\emph{~}

\emph{Proof.} We first treat the case $A={}^{b}\mathbb{C}[X]$ and
put $I:=I_{n}$. Consider the map $\varphi:\, I\rightarrow{\Gamma(\mathbb{C}}^{n},\Omega_{\mathbb{C}^{n}})=\oplus_{i}\mathcal{O}(\mathbb{C}^{n})dx_{i}$
defined by \[
\varphi(P_{0}+\sum_{i=1}^{n}Q_{i}\, dX_{i})=\sum_{i=1}^{n}\mathrm{st}(Q_{i})dx_{i}.\]
 This mapping is well defined since $Q_{i}\in D_{n}$ and we know
that $\varphi(I^{2})=0$. We shall prove that $\mathrm{Ker}(\varphi)=I^{2}$
and that $\varphi$ is surjective.

Let $P=P_{0}+\sum_{i=1}^{n}Q_{i}\, dX_{i}\in I$ such that $\varphi(P)=0$,
so $Q_{i}\in I$. Hence, $Q_{i}=Q_{i,0}+\sum_{j=1}^{n}Q_{ij}\, dX_{j}$
where $Q_{i,0}\in{}{}^{i}\mathbb{C}[X]$. As a consequence, the polynomial
$P$ can be expressed in the form $P=P_{0}+\Sigma_{i=1}^{n}Q_{i,0}dX_{i}+\sum_{i,j=1}^{n}Q_{ij}\, dX_{i}\, dX_{j}$.
By lemma 6.4, we deduce that $P\in I^{2}$. 

Now let $f(x)dx_{i}$ be an element of ${\Gamma(\mathbb{C}}^{n},\Omega_{\mathbb{C}^{n}})$
and $f_{N}\in{}{}^{b}\mathbb{C}[X]$ such that ${}{}^{\circ}f_{N}=f$.
Consider $F_{N}\in{}{}^{b}\mathbb{C}[X]$ such that $\partial^{i}F_{N}=f_{N}$
. If we set \[
P(X,dX):=F_{N}(X+dX_{i})-F_{N}(X)\]
 we obtain $P\in I$ and \[
P(X,dX)-f_{N}(X)dX_{i}=F_{N}(X+dX_{i})-F_{N}(X)-f_{N}(X)dX_{i}\in I^{2}.\]
 In particular we get $\varphi(P)=\varphi(f_{N}(X)dX_{i})=\mathrm{st}(f_{N})\, dx_{i}=f\ (x)dx_{i}$.

The general case will be a consequence of the functorial properties
of both $\Omega$-constructions: The surjectivity of $\gamma_{A}$
is immediate. If $\mathfrak{a}$ is separated, then we can use the
sequence in 6.7 and bijectivity for $A={}^{b}\mathbb{C}[X]$ to conclude
that of $\gamma_{A}$, q.e.d.\\

\textbf{Theorem 6.9} \emph{Let $\mathfrak{a}$ be saturated. Then
the exact sequence of $A^{an}$-modules}\[
0\rightarrow I_{A}/I_{A}^{m+1}\rightarrow D(A)/I_{A}^{m+1}\rightarrow A^{an}\rightarrow0\]
 \emph{splits by a natural section} $s\,:\, A^{an}\rightarrow D(A)/I_{A}^{m+1}$
\emph{for any $m\in\mathbb{N}$.}

\emph{~}

\emph{Proof}. Consider the exact sequence of $A^{an}$-modules \[
0\rightarrow I_{A}/I_{A}^{m+1}\rightarrow D(A)/I_{A}^{m+1}\rightarrow A^{an}\rightarrow0.\]

Let $s\,:\, A^{an}\rightarrow D(A)/I_{A}^{m+1}$ be the morphism,
defined by $s(f)=F\,\mathrm{mod\,}I_{A}^{m+1}$, where $F\in A$ such
that $\mathrm{st}(F)=f$. This map is well defined: Let $F$ and $G$
be two elements in $A$ satisfying $\mathrm{st}(F)=\mathrm{st}(G)=f$.
Then $F-G\in Inf(A)=Inf(A)^{m+1}\subset I_{A}^{m+1}$. Clearly, $s$
is a natural section of the quotient map $D(A)/I_{A}^{m+1}\rightarrow A^{an}$
and this implies that the above sequence of $A^{an}$-modules is in
fact split exact, q.e.d.

~

\noun{Remark} 6.10 It is possible to globalize our $\Omega$-construction
for {*}-bounded schemes. But this approach is only satisfactory if
the affine pieces are define by saturated ideals.

\section*{Appendix}

\addcontentsline{toc}{section}{Appendix}

\subsection*{A Filters and Ideals\protect \protect \\
  }

A filter is a special subset of a partially ordered set. We start
by introducing a filter of sets which is the most used special case,
the partially ordered set is the power set of some set. Filters play
an important role in many fields of mathematical like topology from
where they originated and also lattice theory. In this section, $I$
denotes a nonempty set and $\mathcal{P}(I)$ the set of subsets of
$I$.\\

\textbf{Definition.} A \emph{filter on $\textrm{I}$} is \emph{}a
\emph{}non \emph{}empty \emph{}collection \emph{$\mathcal{F}\subset\mathcal{P}(I)$}
of \emph{}subsets \emph{}of \emph{$I$,} satisfying \emph{}the \emph{}following
\emph{}conditions

i) if $A$,$B$$\in\mathcal{F}$, then $A\cap B\in\mathcal{F}$,

ii) if $A\in\mathcal{F}$ and $A\subset B\subset I$, then $B\in\mathcal{F}$.

We say that $\mathcal{F}$ is \emph{proper} if $\emptyset\not\in\mathcal{F}$
that is $\mathcal{F}\subsetneq\mathcal{P}(I)$ . \\

A trivial example of a filter on $I$ is the collection $\mathcal{F}=\{ I\}$
that consists only of the set $I$ itself.

Let $A$ be a nonempty subset of $I$. Then the collection $\mathcal{F}=\{ X\subset I\,\mid\, X\supset A\}$
is a filter on $I$. It is called the \emph{principal filter} \emph{on}
$I$ \emph{generated by} $A$. In general, if $\mathcal{H}\subset\mathcal{P}(I)$
has the finite intersection property, that is, the intersection of
every nonempty finite sub-collection of $\mathcal{H}$ in nonempty,
then there exist a smallest proper filter containing $\mathcal{H}$,
the filter \emph{generated} \emph{by} $\mathcal{H}$.\\

An \emph{ultrafilter} is a proper filter which satisfies : for any
$A\subset I$, one has either $A\in\mathcal{F}$ or $I\setminus A\in\mathcal{F}$.
It is an easy exercise to prove that $\mathcal{F}$ is an ultrafilter
on $I$ if and only if it is a maximal proper filter on $I$, that
is, there is no proper filter $\mathcal{F}'$on $I$ containing $\mathcal{F}$.
 If $A=\{ a\}$, the principal filter generated by $a$ is an ultrafilter.
If $I$ is finite, then every ultrafilter on $I$ is of this kind.
If $I$ is infinite, there exists nonprincipal filters on $I$, for
example the filter of all cofinite subsets of $I$ .\\

The dual notion of a filter is an (ordered) ideal. Let us recall the
definition of an ideal in a power set of some non empty set.\\

\textbf{Definition.} \emph{} An \emph{ideal on $\textrm{I}$} is \emph{}a
non \emph{}empty \emph{}collection \emph{$\mathfrak{J}\subset\mathcal{P}(I)$}
of \emph{}subsets \emph{}of \emph{$I$,} satisfying \emph{}the \emph{}following
\emph{}conditions

i) if $X,Y$$\in\mathfrak{J}$, then $X\cup Y\in\mathfrak{J}$,

ii) if $Y\in\mathfrak{J}$ and $X\subset Y$, then $X\in\mathfrak{J}$.

~

This definition is that of an ideal in the Boolean ring $\mathcal{P}(I)$,
expressed in the terms of the operations union and intersection in
the Boolean algebra $\mathcal{P}(I)$.\\
 An ideal $\mathfrak{J}$ on $I$ is called a \emph{prime} ideal if
for every $X\subset I$, either $X\in\mathfrak{J}$ or $I\setminus X\in\mathfrak{J}$.
In Boolean algebras, the terms prime ideal and maximal ideal coincide,
as do the terms prime filter and maximal filter.\\

We note that if $\mathcal{F}$ is a filter on $I$ then $\mathfrak{J}=\{ I\setminus X\:\mid\: X\in\mathcal{F}\}$
is an ideal, and vice versa, if $\mathfrak{J}$ is an ideal on $I$,
then $\mathcal{F}=\{ I\setminus X\,\mid\, X\in\mathfrak{J}\}$ is
a filter.\\

By application of Zorn's Lemma, one proves that every proper filter
on a set $I$ can be extended to an ultrafilter on $I$. So if $I$
is infinite, then there exists a non principal ultrafilter on $I$
which contains in particular the filter of cofinite sets.

~

\subparagraph*{Algebraic description of filters \protect \protect \\
 \protect \protect \\
 }

The following algebraic description of filters on a set $I$ is due
to Kochen \cite{Kochen}. Let $I$ be an index set, and $K_{i}$ (where
$i\in I$) be a family of commutative fields. We put $R:={\displaystyle \Pi{}_{i\in I}K_{i}}$.
Then $R$ is a von Neumann regular ring, that is, it satisfies: $\forall x\in R\,\exists y\in R$
such that $x=xyx$. Kochen \cite{Kochen} proved that the ideal structure
of the ring $R$ can be described by filters on the set $I$. We use
the following notations: for $f\in R$, let $Z(f)=\{ i\in I\,\mid\, f(i)=0\}$.
If $\mathfrak{a}$ is any ideal in $R$, we put $\mathcal{F}(\mathfrak{a}):=\{ Z(f)\,\mid\, f\in\mathfrak{a}\}$.
If $\mathcal{F}$ is a filter on $I$, we define $\mathrm{Id}(\mathcal{F}):=\{ f\in R\,\mid\, Z(f)\in\mathcal{F}\}$.
One verifies easily that $\mathrm{Id}(\mathcal{F})$ is an ideal of
$R$ and $\mathcal{F}(\mathfrak{a})$ is a filter on $I$.\\

\textbf{Theorem.} \emph{The construction above gives a one-one correspondence
between the family of ideals in $R$ and the family of filters on
$I$. Furthermore, prime ideals correspond to ultrafilters. If $\mathfrak{a}$
is an ideal, then $R/\mathfrak{a}\simeq R/\mathcal{F}(\mathfrak{a})$.
If $\mathfrak{a}$ is a prime ideal, then $R_{\mathfrak{a}}\simeq R/\mathfrak{a}$.}\\

We want to generalize this to arbitrary commutative rings, instead
of fields. We first consider the case of local rings $(R_{i},\mathfrak{m}_{i})_{i\in I}$,
each one different from the zero ring. Put $R:={\displaystyle \Pi{}_{i\in I}R_{i}}$.
For $f=(f_{i})_{i}\in R$, we define the \emph{variety} \emph{of}
$f$ by

\[
\mathrm{V}(f):=\{ i\in I\,\mid\, f_{i}\in\mathfrak{m}_{i}\}.\]
 If $\mathfrak{a}\subset R$ is an ideal, then $\mathcal{F}(\mathfrak{a}):=\{\mathrm{V}(f)\,\mid\, f\in\mathfrak{a}\}$
is a filter on $I$. This follows from the following facts: take $f,g\in R$,
then $\mathrm{V}(fg)=\mathrm{V}(f)\cup\mathrm{V}(g)$; there are $\varepsilon,\delta\in R$,
such that $\mathrm{V}(\varepsilon f+\delta g)=\mathrm{V}(f)\cap\mathrm{V}(g)$.
Moreover, if $\mathfrak{a}\subset R$ is a proper ideal, then $\mathcal{F}(\mathfrak{a})$
is proper too.

There is the dual construction: Let $\mathcal{F}$ be a filter on
$I$. We put \[
\mathrm{Id}(\mathcal{F}):=\{ f\in R\,\mid\,\mathrm{V}(f)\in\mathcal{F}\}.\]
 Then $\mathrm{Id}(\mathcal{F})$ is an ideal of $R$. If $\mathcal{F}$
is proper, then $\mathrm{Id}(\mathcal{F})$ is proper too. We note
that both constructions conserve inclusions and we have always

\[
\mathfrak{a}\subset\mathrm{Id}(\mathcal{F}(\mathfrak{a})),\]

\[
\mathcal{F}=\mathcal{F}(\mathrm{Id}(\mathcal{F})).\]
 One can show

~

\textbf{Proposition.} i) \emph{Let} $\mathfrak{m}\subset R$ \emph{be}
\emph{an} \emph{ideal}. \emph{If} $\mathfrak{m}$ \emph{is} \emph{maximal},
\emph{then} $\mathcal{F}(\mathfrak{m})$ \emph{is} \emph{an} \emph{ultrafilter
on $I$.}

ii) \emph{Let} $\mathcal{F}$ \emph{be a} \emph{proper filter} \emph{on}
$I$. \emph{Then} $\mathrm{Id}(\mathcal{F})$ \emph{is a} \emph{maximal
ideal} \emph{of $R$ if and} \emph{only if} $\mathcal{F}$ \emph{is
an} \emph{ultrafilter on $I$.}

~

\textbf{Corollary.} \emph{The} \emph{mapping}

\[
\mathcal{F}:\mathrm{Specmax}(R)\longrightarrow\mathrm{ufil}(I)\]
 \emph{which} \emph{associates to every maximal ideal of} $R$ \emph{its
ultrafilter}, \emph{is well-defined and bijective}. 

~\\
We come to the case of not necessarily local (but non zero) rings
$R_{i}$, $i\in I$. For a commutative ring $A$, we define

\[
\mathrm{W}(A):=\prod_{\mathfrak{m}\in\mathrm{Specmax}(A)}A_{\mathfrak{m}}.\]
 Then we prove 

~

\textbf{Theorem}. \emph{Every} \emph{maximal} \emph{ideal} \emph{of}
$R:={\displaystyle \Pi{}_{i\in I}R_{i}}$ \emph{is induced} \emph{by
a} \emph{(not necessarily} \emph{unique) maximal ideal of} ${\displaystyle \Pi{}_{i\in I}\mathrm{W}(R_{i})}$,
\emph{i.e. by an ultrafilter} \emph{on the disjoint union} $\sqcup_{i\in I}\mathrm{Specmax}(R_{i})$.

~

For the \emph{proof}, we take a maximal ideal $\mathfrak{m}$ of $R$.
Put $S:={\displaystyle \Pi{}_{i\in I}\mathrm{W}(R_{i})}$. Then we
note that $\mathfrak{m}S$ is a proper ideal of $S$. Otherwise $1_{S}$
would be a \emph{finite} linear combination of elements of $\mathfrak{m}$.
But this implies, by the local-global principle of commutative algebra,
that these elements generate each $R_{i}$, and so finally $\mathfrak{m}=R$
which is impossible. It follows that $\mathfrak{m}S$ is contained
in at least one maximal ideal of $S$. We obtain our result by applying
the above proposition.

\subsubsection*{$\delta$-stable ultrafilters\protect \protect \\
  }

There are some applications in analysis where a particular kind of
ultrafilter are very advantageous. In \cite{Choquet2} Choquet constructed
$\delta$-stable ultrafilters, previously called absolutely 1-simple
ultrafilters in \cite{Choquet1}.

~

\textbf{Definition \cite{Choquet2}} We say that an ultrafilter $\mathcal{U}$
is $\delta$-\emph{stable} if $(J_{n})_{n\in\mathbb{N}}$ is any sequence
of elements of $\mathcal{U}$, there is a set $J_{\infty}\in\mathcal{U}$,
almost contained in each $J_{n}$, so $J_{\infty}\setminus J_{n}$
is \emph{finite} for each $n$.

~

Choquet has shown (see \cite{Choquet1}, theorem 6 ) that under the
continuum hypothesis, there exist $\delta$-stable ultrafilters over
$\mathbb{N}$.

Stroyan and Luxemburg (see \cite{Stroyan Luxemburg}, theorem 7.1.1
p.175) proved that each infinite $\lambda$ natural number can be
represented by a sequence $\lim\lambda(j)=+\infty$. It is straightforward
that this is still valid for infinite positive reals 

~

\textbf{Theorem.} \emph{If} $\,{}^{*}\mathbb{R}$ \emph{is} \emph{a}
$\delta$-\emph{stable} \emph{ultrapower} \emph{of} $\mathbb{R}$
\emph{and} $\lambda{=[\lambda(j)]\in{}}^{*}\mathbb{R}_{\infty}^{+}$
\emph{is} \emph{infinite} \emph{positive real}, \emph{then there}
\emph{exists a} \emph{set} $J_{\infty}\in\mathcal{U}$ \emph{such
that} \[
\lim_{j\rightarrow\infty,\, j\in J_{\infty}}\lambda(j)=+\infty.\]

\subparagraph*{The standard part of hyperreals\protect \protect \\
 \protect \protect \\
 }

Let $\mathcal{U}$ be a nonprincipal ultrafilter on $\mathbb{N}$
and$^{*}\mathbb{R}$ the corresponding ultrapower of $\mathbb{R}$,
that is, $^{*}\mathbb{R}=\mathbb{R^{\mathbb{N}}}/\mathcal{U}$ with
$(x_{i})\sim(y_{i})$ if $\{ i\in\mathbb{N},\, x_{i}=y_{i}\}\in\mathcal{U}$.
Then $({}^{*}\mathbb{R},+,.,<)$ is an ordered field extension of
$(\mathbb{R},+,.,<)$. Elements of $^{*}\mathbb{R}$ are called (hyper)real
numbers.

A hyperreal $r\in{}{}^{*}\mathbb{R}$ is bounded or finite if $|r|<n$
for some $n\in\mathbb{N}$, and infinitesimal if $|r|<\frac{1}{n}$
for every $n\in\mathbb{N}$, $n\geq1$. Let $s\in{}{}^{*}\mathbb{R}$.
We say that $r$ and $s$ are infinitely close if $r-s$ is infinitesimal.
We write $r\approx s$ in this case. We denote by $^{b}\mathbb{R}$
the subring of $^{*}\mathbb{R}$ of bounded numbers and by $^{i}\mathbb{R}$
the ideal of $^{b}\mathbb{R}$ of all infinitesimals.\\

\textbf{Theorem.} \emph{If $\rho$ is bounded, then there exists a
unique real $r$ such that $\rho\approx r$.}\\

We call $r$ is the standard part of $\rho$ and write $r=°\rho$
or $r=\mathrm{st}(\rho)$. The map\[
\mathrm{st}:\,^{b}\mathbb{R}\longrightarrow\mathbb{R}\]
 is called the \emph{standard} \emph{part} \emph{map}. It is an order
preserving homomorphism from the ring $^{b}\mathbb{R}$ onto $\mathbb{R}$.
The kernel {}``st'' is $^{i}\mathbb{R}$ and the quotient ring $^{b}\mathbb{R}{/\,}^{i}\mathbb{R}$
is isomorphic to $\mathbb{R}$. \\
We restrict {}``st'', the standard part map, to $^{*}\mathbb{Q}$,
the set of hyperrationals, $\mathrm{st}_{^{*}\mathbb{Q}{\cap{}}^{b}\mathbb{R}}{:\,}^{*}\mathbb{Q}{\cap{}}^{b}\mathbb{R}\longrightarrow\mathbb{R}$.
By density of $\mathbb{Q}$ in $\mathbb{R}$, $\mathrm{st}_{^{*}\mathbb{Q}{\cap{}}^{b}\mathbb{R}}$
is surjective and the quotient $^{*}\mathbb{Q}{\cap{}}^{b}\mathbb{R}/{}^{*}\mathbb{Q}{\cap{}}^{i}\mathbb{R}$
is isomorphic to $\mathbb{R}$ (the completion of $\mathbb{Q}$).

Brunjes and Serpé \cite{Brunjes2} treated the case of a non trivially
valued filed $(K,|.|)$ with locally compact completion $(\hat{K},|.|)$.
They proved that $^{*}K{}^{{\rm fin}}/{}^{*}K{}^{{\rm inf}}$ is isomorphic
to $\hat{K}$, where $^{*}K{}^{{\rm fin}}$ denotes the set of finite
elements of $^{*}K$ and $^{*}K{}^{{\rm inf}}$ the set of infinitesimal
elements of $^{*}K$. \\

\subsection*{B Principles of Nonstandard Analysis}

\subparagraph*{Transfer Principle\protect \protect \\
 \protect \protect \\
 }

The fundamental property of ultraproducts is the following:\\

\textbf{Theorem} (\L \v{o}s) \emph{Let $L$ be a first order language
and $A_{i}$, $i\in I$ be structures for $L$. Let $\mathcal{U}$
be an ultrafilter on $I$ and $^{*}A=(\Pi{}_{i\in I}A_{i})/\mathcal{U}$.
Then for any first order formula $\varphi(x_{1,},\ldots,x_{n})$ with
$x_{1},\ldots,x_{n}$ its only free variables and $[a_{1}],\ldots,[a_{n}]\in\,^{*}A$,
we have that $^{*}\varphi([a_{1}],\ldots,[a_{n}])$ is true in $^{*}A$
if and only if $\{ i\in I\,:\,\varphi(a_{1}(i),\ldots,a_{n}(i))\mbox{ is true in }A_{i}\}\in\mathcal{U}$.
In particular if $\varphi$ is a sentence, then $\varphi$ is true
if and only if $^{*}\varphi$ is true.}

~\\
We present some interesting applications of \L \v{o}s' Theorem.

\begin{enumerate}
\item Let $(k_{i})_{i\in I}$ be a family of (algebraically closed) fields.
Then $^{*}k$, their ultraproduct, is again a (algebraically closed)
field. \\
If for each prime $p$, only finitely many $k_{i}$ have characteristic
$p$, then $^{*}k$ has characteristic zero. In particular, if $P$
is an \emph{infinite} prime in $^{*}\mathbb{Z}$ then $^{*}\mathbb{Z}/P{}^{*}\mathbb{Z}$
is a field of characteristic zero. Its algebraic closure is (non-canonically)
isomorphic to $\mathbb{C}$, because each field is algebraically closed
with the same cardinality which is equal to that of the continuum.
In the same manner, one can show that $^{*}\mathbb{C}$ is (non-canonically)
isomorphic to $\mathbb{C}$. 
\item Let $(A_{i})_{i\in I}$be a family of local rings with maximal ideal
$\mathfrak{m}_{i}$ and residue filed $k_{i}=A_{i}/\mathfrak{m}_{i}$.
Then $^{*}A$, their ultraproduct, is a local ring with maximal ideal
$^{*}\mathfrak{m}$ and residue field $^{*}k={}^{*}A/{}^{*}\mathfrak{m}$.
In fact, a ring is local if and only if the sum of two non-units is
a non unit. 
\end{enumerate}

\subparagraph*{Permanence Principle\protect \protect \\
 \protect \protect \\
 }

This principle asserts that certain functions can be extended to larger
domains than those over which they are originally defined.

~

\textbf{Theorem} \emph{Let $\varphi(x)$ be an internal formula with
the only free variable $x$. Then}

\emph{i) If there exists $k\in\mathbb{N}$ such that $\varphi(n)$
is true for all $n\in\mathbb{N}$ with $k\leq n$, then there exists
$K\in\,^{*}\mathbb{N}\setminus\mathbb{N}$ such that $\varphi(n)$
is true for all $n\in\,^{*}\mathbb{N}$ with $k\leq n\leq K$.}

\emph{ii) If there exists $K\in\,^{*}\mathbb{N}\setminus\mathbb{N}$
such that $\varphi(n)$ is true for all $n\in\,^{*}\mathbb{N}\setminus\mathbb{N}$
with $n\leq K$, then there exists $k\in\mathbb{N}$ such that $\varphi(n)$
is true for all $n\in\,^{*}\mathbb{N}$ with $k\leq n\leq K$.}

\emph{iii) If $\varphi(x)$ holds for each infinitesimal $x$, then
there is a standard $r>0$ in $\mathbb{R}$ so that $\varphi(b)$
holds for all $b$ with $|b|\leq r$ in $^{*}\mathbb{R}$.}

\emph{~}

\textbf{Corollary~(Spillover Principle)} \emph{Let $A$ be an internal
subset of $^{*}\mathbb{R}$.}

\emph{i) If $A$ contains all standard natural numbers, then $A$
contains an infinite natural number.}

\emph{ii) If $A$ contains all infinite natural numbers, then $A$
contains a standard natural number.}

\subparagraph*{Concurrence or Saturation Principle\protect \protect \\
 \protect \protect \\
 }

This principle provides us with new objects in nonstandard extensions.
Let $X$ be an infinite set and $R$ be a binary relation on $X$.
$R$ is called $concurrent$ if for any finite subset $\{ x_{1},\ldots,x_{n}\}$
of the domain of $R$ there exists an element $y$ with $x_{i}\, R\, y$
for all i between $1$ and $n$.

~

\textbf{Theorem}\emph{. If} $R$ \emph{is concurrent}, \emph{then}
\emph{there} \emph{exists} $b{\in{}}^{*}X$ \emph{such} \emph{that}
$^{*}x(^{*}R)b$ \emph{for} \emph{all} $x\in{\rm dom}\, R$.

~

This is also called the \emph{enlargement property}: $^{*}X$ is an
\emph{enlargement} of $X$ if every concurrent relation verifies the
enlargement property. Let us consider some examples:

1. $R(x,y):(x\in\mathbb{R})\wedge(y\in\mathbb{R})(x\leq y)$ $R$
is a concurrent relation on $\mathbb{R}$, so there exists $b\in\,^{*}\mathbb{R}$
such that $x\leq b$ for any $x\in\mathbb{R}$. Evidently, $b$ is
an \emph{infinite} number.

2. Let $X$ be a topological space, $x\in X$ and $^{*}X$ an enlargement
of $X$. Then there exists an internal open set $W$ in $^{*}X$ containing
$x$, such that $W\subset\mu(x)$. \\
 We denote by $\Omega_{x}$ the system of open neighborhoods of $x$.
The binary relation $R$ defined on $\Omega_{x}\times\Omega_{x}$
defined by $R(U,V)$ if $U\supset V$ is concurrent. So, the enlargement
property guarantees the existence of $W{\in{}}^{*}\Omega_{x}$ such
that $\,^{*}U\supset W$, for every $U\in\Omega_{x}$. In a metric
space, one can take $W$ an $*$-ball with $x$ as center and positive
infinitesimal radius.

\subsection*{C Nonstandard Topologies \protect \protect \\
 }

\subparagraph*{Topological spaces\protect \protect \\
 \protect \protect \\
 }

Let $(X,\mathcal{T})$ be a topological space. In the literature,
usually two topologies on the nonstandard extension of $^{*}X$ are
considered. The first one, called $Q$-topology, introduced by Robinson
\cite{Robinson1}, has a basis, consisting of elements in $^{*}\mathcal{T}$.
In fact, by transfer, elements of $^{*}\mathcal{T}$ are stable under
$*$-finite intersection, hence under finite intersection and under
internal union. Furthermore, $\emptyset$ ,$^{*}X$$\in\,$$^{*}\mathcal{T}$.
Elements in the $Q$-topology are called $Q$-open and elements in
the base $^{*}\mathcal{T}$ are called $*$-open subsets of $^{*}X$.
A $*$-open set is clearly $Q$-open .The converse is false in general:
there are external sets which are open in this topology. Robinson
proved that if $Y$ is an internal subset of $^{*}X$, then $Y$ is
$Q$-open if and only if $Y$ is $*$-open. The second one, called
$S$-topology by Robinson, has as basis $\mathcal{B=}\,\{^{*}U\,\mid\, U\in\mathcal{T}\,\}$.
The $Q$-topology is finer than the $S$-topology. \\

Let $x\in X$ be a point. The \emph{monad} of $x$ or \emph{halo}
is the subset 

~

${\displaystyle \mu(x)=\cap{}_{U\in\mathcal{T},x\in U}}{}^{*}U$ . 

~\\
A point $y\in\,^{*}X$ is, by definition, \emph{near-standard}, if
there exists $x\in X$ such that $y\in\mu(x)$. We write $y\approx x$
in this case. The set of near-standard points of $^{*}X$ is $\mbox{ns}{\rm (^{*}X)=}\cup_{x\in X}{\displaystyle \mu(x)}$.
Many topological properties of $X$ can be expressed via monads, so
the halos of $^{*}X$ which encode the topology of $X$. \\

For a Hausdorff topological space, we have a standard part map ${\rm \mathrm{st}:\,\mathrm{ns}{({}}^{*}X)}\rightarrow X$
defined as follows: we set ${\rm \mathrm{st}}(y){={}}^{\circ}y=x$,
for every $y\in\mu(x)$ and $x\in X$,. This map is well defined because
in the case of Hausdorff spaces, halos of standard points constitute
a disjoint partion of near standard points. Even for a non-Hausdorff
space, we can define the standard part of a set $B{\subset{}}^{*}X$,
by setting:

\[
\mathrm{st}(B){={}}^{\circ}B\,=\{ x\in X\,\mid\,\mu(x)\cap B\not=\emptyset\}.\]
 ~\\
Under certain cardinality restrictions (expressed via the notion {}``$\kappa$-saturated
enlargement''), standard parts of internal subsets turn out to be
always closed. The following result can be found in the book of Hurd
and Loeb (\cite{Hurd Loeb} p.117)

~

\textbf{Theorem.} \emph{Assume that $(X,\mathcal{T})$ is a topological
space and $^{*}\mathcal{T}$ is a $\kappa$-saturated enlargement
of $\mathcal{T}$ with $\kappa>\mathrm{card}\,\mathcal{T}$. }

\emph{i) If $B$ is internal subset of $^{*}X$, then $\mathrm{st}(B)$
is closed.}

\emph{ii) If $B$ is internal subset of near-standard points of $^{*}X$,
then $\mathrm{st}(B)$ is compact.}

~

If $X$ satisfies the first axiom of countability (for example if
$X$ is a metric space), one can use just $\aleph_{1}$-saturation.

\subparagraph*{Metric spaces\protect \protect \\
 \protect \protect \\
 }

Let $X$ be a metric space with distance function $\rho$ and $\Gamma$
be the set of all open balls B, where $B(x,r)=\{ y\in X\,\mid\,\rho(x,y)<r\}$
for a point $px\in X$ and a positive real $r$. We describe the monad
of a point $x{\in{}}^{*}X$ by $\mu(x)=\cap{}^{*}B(x,r)=\{ y{\in{}}^{*}X\,\mid\,{}^{*}\rho(x,y)\approx0\}$.
Obviously, $^{*}\Gamma$ forms a basis for the $Q$-topology of $^{*}X$.
Let $\  p{\in{}}^{*}X$ and $r$ a \emph{standard} positive number.
We put $S(p,r)=\{ q{\in{}}^{*}X\,\mid\,{}^{\circ}\rho(p,q)<r\}$.
The $S$-balls give us a topology in $^{*}X$, which we is the $S$-topology.
The space $^{*}X$, endowed with the $S$-topology in not Hausdorff.
In fact, if $x,y{\in{}}^{*}X$ such that $x\approx y$ and if $r>0$
is a standard positive real, we obtain $S(x,r)=S(y,r)$.\\

Let $(X,\rho),(Y,d)$ be two metric spaces and $f:\,^{*}X\rightarrow\,{}^{*}Y$
be a map from $^{*}X$ to $^{*}Y$. Let $x\in\,^{*}X$ be a point.
We say that $f$ is \emph{$S$-continuous} at $x$, if $f$ is continuous
at $x$ as a map from $^{*}X$ to $^{*}Y$ where both sets are equipped
with the $S$-topology, that is, for every standard $\varepsilon>0$,
there exists a standard $\delta>0$ such that $^{*}d(f(x),f(y))<\varepsilon$
for each $y\in\,^{*}X$, such that $^{*}\rho(x,y)<\delta$. There
is a simple characterization of \emph{$S$-}continuity of internal
maps. If $f$ is an \emph{internal} map, then $f$ is \emph{$S$-}continuous
at $x$ if and only if $f(\mu(x))\subset\mu(f(x)$).\\

Let $X$ be a metric space with distance function $\rho$. A point
$p{\in{}}^{*}X$ is called \emph{bounded} (or \emph{finite}) if there
exists a standard point $q\in X$ such that $^{*}\rho(p,q)$ is \emph{bounded},
i.e. there exists a standard positive real $m>0$ such that $p\in{}^{*}B(q,m)$.
We denote by $\mathrm{bd}{({}}^{*}X)$ the set of bounded points of
$^{*}X$. Robinson proved (see \cite{Robinson1} theorem 4.3.1 p.100)
that a metric space $X$ is bounded if and only if $^{*}X=\mathrm{bd}{({}}^{*}X)$.
Clearly, we always have $\mathrm{ns}({}^{*}X)\subset\mathrm{bd}({}^{*}X)$.\\

\subparagraph{Standard part of a map\protect \protect \\
 \protect \protect \\
}

Let $X$ and $Y$ be two topological spaces. $\mathcal{F}(X,Y)$ (resp.
$\mathcal{C}(X,Y)$) denotes the set of all (resp. continuous) maps
from $X$ to $Y$. As usual $\mathrm{ns}(\:^{*}X)$ is the set of
nearstandard points and ${\displaystyle \mathrm{cpt}(\:^{*}X)=\cup{}_{K\subset X,K\mathrm{\, compact}}\,^{*}K}$
is the set of {}``compact points''.

~

One of the most powerfull tools in nonstandard analysis is taking
standard parts of objects of the nonstandard universe. If $Y$ is
Hausdorff, we can define for every $f\in\,^{*}\mathcal{C}(X,Y)$ with
$f(\,^{*}x)\in\mathrm{ns}(\,^{*}X)$ for all $x\in X$, the standard
part function $\mathrm{st}\, f:X\rightarrow Y$ by $\mathrm{st}f(x):=\,^{\circ}\left(f(\:^{*}x)\right).$
Thus, we have a mapping\[
\mathrm{st}:\{ f\in^{*}\mathcal{C}(X,Y)\,\mid\, f(\,^{*}x)\in\mathrm{ns}(\,^{*}Y)\,\mathrm{\,\forall\,}x\in\,^{*}X\}\longrightarrow\mathcal{F}(X,Y)\]

It is known that $\mathrm{st}\, f$ is in general not a continuous
function.

~

\textbf{Theorem.} \emph{Let $B$ be a subset of $^{*}X$ and $f:\, B\rightarrow\,{}^{*}Y$
an internal map such that $f(B)\subset{\rm \mathrm{ns}}(^{*}Y)$ and
$f$ is S-continuous on $B.$ Then $\mathrm{st}\, f$ exists and is
continuous on $^{\circ}B$.}

~

We endow $\mathcal{C}(X,Y)$ with $\tau$, the compact open topology.
Let $f\in\mathcal{C}(X,Y)$ be a mapping. We denote by $\mu_{\tau}(f)$
the halo of $f$ with respect to the compact open topology. One easily
verifies that

\[
\mu_{\tau}(f)=\{ g\in\,^{*}\mathcal{C}(X,Y)\,\mid\, g(x)\approx\,^{*}f(x)\,\,\mathrm{\forall\,}x\in\mathrm{cpt}(\,^{*}X)\}\]

If $X$ is \emph{locally compact}, then $\mathrm{cpt}(\,^{*}X)=\mathrm{ns}(\,^{*}X)$,
hence\[
\mu_{\tau}(f)=\{ g\in\,^{*}\mathcal{C}(X,Y)\,\mid\, g(x)\approx\,^{*}f(x)\,\,\mathrm{\forall\,}x\in\mathrm{ns}(\,^{*}X)\}.\]

\textbf{Theorem.} \emph{Let $X,Y$ be two Hausdorff topological spaces
and $X$ locally compact. Then $\mathrm{ns}_{\tau}(\,^{*}\mathcal{C}(X,Y))$
consists of all internal maps sending nearstandard points of $^{*}X$
to nearstandard points of $^{*}Y$ and which are S-continuous on $\mathrm{ns}(\,^{*}X)$. }

\emph{In particular, nearstandard functions for the compact-open topology
are those which are bounded and S-continuous on $\mathrm{ns}(\,^{*}X)$.}

~

\[
\_\]
~\\
$\mathcal{AMS}$-classification: 32C15, 14A15, 32E10, 26E35, 58A10

~\\
\emph{Key}-\emph{words}: complex spaces, nonstandard schemes, Stein
algebras, analytic Nullstellensatz, generic points, internal polynomials,
ultraproducts of algebras

~\\
\emph{Author's} \emph{addresses}: 

Adel Khalfallah, Institut Préparatoire aux Etudes d'Ingénieur, Rue
Ibn El Jazzar, 5019 Monastir, Tunisie 

\emph{email}: adel.khalfallah@ipeim.rnu.tn 

Siegmund Kosarew, Institut Fourier Université de Grenoble 1, 100 rue
des Mathématiques, 38402 Saint Martin d'Hères, France 

\emph{email}: siegmund.kosarew@ujf-grenoble.fr

\end{document}